\documentclass[12pt, leqno]{article}
\usepackage{amsmath, amsthm, enumitem, thmtools}
\usepackage{amsfonts}
\usepackage{amssymb}
\usepackage[usenames,dvipsnames]{color}
\usepackage{graphicx}
\usepackage{bbm}
\usepackage[toc]{appendix}
\usepackage{bm}
\usepackage{listings}
\usepackage{nicefrac} 
\usepackage{xfrac}    
\usepackage[utf8]{inputenc}
\usepackage[T1]{fontenc}
\usepackage{hyperref}
\usepackage{pgffor}
\usepackage{xcolor}
\usepackage{soul}

\usepackage{orcidlink} 

\usepackage{color,calc}

\makeatletter
\def\namedlabel#1#2{\begingroup
    #2%
    \def\@currentlabel{#2}%
    \phantomsection\label{#1}\endgroup
}
\definecolor{shade}{gray}{0.8}
        {
          \raggedright
        \setlength{\rightmargin}{\leftmargin}
        \setlength{\itemsep}{-12pt}
        \setlength{\parsep}{20pt}
        \begin{lrbox}{\@tempboxa}%
        \begin{minipage}{\linewidth-2\fboxsep}
        }%
        {
        \end{minipage}%
        \end{lrbox}%
        \fcolorbox{black}{shade}{\usebox{\@tempboxa}}\newline
        }%
\makeatother

\newtheorem{theorem}{Theorem}
\newtheorem{lemma}{Lemma}

\renewcommand{\eqref}[1]{\hyperref[#1]{(\ref*{#1})}}

\newcommand{\dd}{\mathrm{d}}

\newcommand*{\norm}[1]{\lVert #1 \rVert}

\newtheorem{remark}{Remark}

\renewcommand{\ln}{\log}

\newcommand*{\pref}[1]{\hyperref[#1]{(\ref*{#1})}}
\newcommand*{\refpref}[2]{\hyperref[#2]{\ref*{#1}(\ref*{#2})}}

\newcommand{\N}{\mathbb{N}}

\newcommand{\sT}{\textsf{T}}
\newcommand{\sP}{\textsf{P}}

\newcommand{\sZ}{\textsf{Z}}

\newcommand{\sG}{\textsf{G}}
\newcommand{\sV}{\textsf{V}}
\newcommand{\sU}{\textsf{U}}

\newcommand{\sv}{\textsf{v}}
\newcommand{\su}{\textsf{u}}
\newcommand{\sA}{\textsf{A}}

\newcommand{\sm}{\textsf{m}}

\newcommand{\bP}{\textsf{P}}
\newcommand{\sJ}{\textsf{J}}

\definecolor{amethyst}{rgb}{0.6, 0.4, 0.8}
\definecolor{applegreen}{rgb}{0.55, 0.71, 0.0}
\definecolor{aqua}{rgb}{0.0, 1.0, 1.0}
\definecolor{asparagus}{rgb}{0.53, 0.66, 0.42}
\definecolor{amber(sae/ece)}{rgb}{1.0, 0.49, 0.0}
 	\definecolor{armygreen}{rgb}{0.29, 0.33, 0.13}
	\definecolor{shitbrown}{rgb}{0.43, 0.21, 0.1}
	\definecolor{brightpink}{rgb}{1.0, 0.0, 0.5}
	\definecolor{brightube}{rgb}{0.82, 0.62, 0.91}
	 	\definecolor{byzantine}{rgb}{0.74, 0.2, 0.64}
		\definecolor{chartreuse(web)}{rgb}{0.5, 1.0, 0.0}
\allowdisplaybreaks

\title{Critical non-local spatial branching processes with infinite variance conditioned on survival}

\author{Natalia  Cardona-Tob\'on\thanks{Universidad Nacional de Colombia, Departamento de Estadística,
Carrera 45 No 26-85, CP 111321 
Bogotá, Colombia. E-mail: \texttt{
ncardonat@unal.edu.co} }, Andreas E. Kyprianou\thanks{
Department of Statistics,
University of Warwick,
Coventry
CV4 7AL, UK. E-mail: \texttt{ \{andreas.kyprianou\}, \{pedro.martin-chavez.1\} @warwick.ac.uk}
}, \ Pedro Mart\'in-Ch\'avez\orcidlink{0000-0001-5530-3138}\footnotemark[2] 
}

\newcommand{\bra}[1]{\ensuremath{\left[ #1\right] }}

\begin{document}

\maketitle
\begin{abstract}
\noindent {
We consider the setting of either a general non-local branching particle process or a general non-local superprocess. Under the assumption that the mean semigroup has a Perron--Frobenious type behaviour in combination with a regularly varying assumption on the reproductive point process,  which permits infinite second moments, we consider  sufficient conditions that ensure limiting distributional stability when conditioned on survival  at criticality.  

We offer two main results. Under the aforesaid  conditions, our first main contribution establishes the polynomial decay in time of the  survival probability in the spirit of a classical Kolmogorov limit. The second main contribution pertains to the stability, when conditioning on survival, in the spirit of a Yaglom limit. In both cases our proofs work equally well for the analogous setting of superprocesses with non-local branching.

Our results complete a series of articles for various families of non spatial branching processes, cf. \cite{zolotarev1957, slack1968branching, vatutin1977, Goldstein_Hoppe}, as well spatial branching processes, cf. \cite{ren2020limit}. The generality of our results  is consistent with older work of \cite{HeringHoppe, AH3} who dealt with a similar context for branching particle systems. We add to this by providing a general framework for non-local superprocesses. 
}

\medskip

\medskip

\end{abstract}

\section{Introduction}
In this article, we  show that two classical results for  Bienaym\'e--Galton--Watson processes remain valid as universal principles in the setting of critical non-local general branching Markov processes and non-local superprocesses with infinite second moments. In particular, to the best of our knowledge, 
our results offer the first results at this level of generality for non-local superprocesses. Moreover, several known results appear as special cases of our general framework. 

In what follows, we focus on critical processes. In the setting of Bienaym\'e--Galton--Watson processes, the notion of criticality is dictated by the mean number of offspring. In the general setting we present in this work, the notion of criticality pertains to the value of an assumed lead eigenvalue for the mean semigroup. 

Our first main result extends the classical Kolmogorov's estimation for the survival probability of discrete-time Bienaym\'e--Galton--Watson processes where the offspring distribution is allowed to have infinite variance. Specifically, we show that, at criticality, the survival probability decays according to a regularly varying function. Our second main result concerns a Yaglom type limit: we prove that, in the critical case, a properly normalized version of the process, conditioned on non-extinction, converges in distribution to a non-degenerate random variable.

The asymptotic behaviour of critical branching processes with infinite variance was first studied by Zolotarev \cite{zolotarev1957} in the context of certain continuous-time processes. Later on, Slack \cite{slack1968branching} investigated discrete-time Bienaym\'e--Galton--Watson processes $(Z_n, n\ge 0)$ whose offspring generating function is given by
\begin{equation}\label{eq:genh}
h(s)= s + (1-s)^{1+\alpha}\ell(1-s), \quad s\ge 0,
\end{equation} 
where $\alpha\in (0,1]$ and $\ell$ is a slowly varying function at 0. The author showed that the survival probability decays as a regularly varying function at zero with index $-1/\alpha$. Furthermore, Slack established the Yaglom-type limit
\begin{equation*}
	\lim_{n\to \infty} \mathbb{E}\left[\exp\left(-\theta n^{-1/\alpha} \tilde{\ell}(n) Z_n \right) \Big|\Big. \  Z_n >0 \right] = 1- \frac{\theta}{(1+  \theta^{\alpha})^{1/\alpha}},\quad \theta \ge 0,
\end{equation*}
where $\tilde{\ell}$ is a slowly varying function at $\infty$. In \cite{slack1972further}, Slack also proved the converse: if the left-hand side of previous equality converges to a non-degenerate limit for all points of continuity, then the generating function must be of the form \eqref{eq:genh}. In addition, the case where the probability generating function is given by \eqref{eq:genh} with $\alpha = 0$ was studied in \cite{nagaev2007critical}.

These results have since been extended to critical multi-type Bienaym\'e--Galton--Watson processes with infinite variance by Vatutin \cite{vatutin1977} and by Goldstein and Hoppe \cite{Goldstein_Hoppe}. More recently, Vatutin et al. \cite{vatutin2021} generalized these findings to the case of critical Bienaym\'e--Galton--Watson processes with countably many types and infinite second moments. Additional {other} extensions can be found in \cite{smadi2015reduced}.

Asmussen and Hering \cite[Theorem 6.4.2]{AH3}, as well as Hering and Hoppe \cite{HeringHoppe}, obtained counterpart results for critical non-local branching Markov processes in general state spaces. 
Our approach differs  from theirs, as our proofs {are based} on a `spine decomposition'. This important difference enables us to extend the results to the superprocess setting. We should also note that one of our assumptions is not directly comparable with the corresponding condition in~\cite{HeringHoppe}. A more detailed discussion of this point is provided later in Remark \ref{rem-assumptionH5}.

In the context of superprocesses, similar questions have also been addressed, {for instance,} in 
critical continuous-state branching processes with infinite variance. Kyprianou and Pardo \cite{Kyprianou_Pardo_2008} showed analogous results assuming stable branching mechanisms, and their work was later extended by Ren, Yang, and Zhao \cite{RenYangZhao2014} to more branching mechanism given by a regularly varying function. More recently, Ren, Song, and Sun \cite{ren2020limit} investigated critical local superprocesses with spatial motion and spatially dependent stable branching. 

Therefore, one of our main contributions lies in providing a framework that allows the extension of the infinite variance Kolmogorov and Yaglom limits to general non-local superprocess setting. Furthermore, this framework being robust for non-local branching processes also allow us to recover several of the previous works in this area.
Further discussions and examples are provided in Section~\ref{sec:examples}.

\section{Non-local spatial branching processes}
Let us spend some time describing the general setting in which we wish to work.
Let $E$ be a Lusin space. Throughout, {we} write $B(E)$ for the Banach space of bounded measurable functions on $E$ with supremum norm $\norm{\cdot}$, $B^{+}(E)$ for the space of non-negative bounded measurable functions on $E$ and $B^{+}_1(E)$ for the subset of functions in $B^{+}(E)$ that are uniformly bounded by unity. We are interested in spatial branching processes that are defined in terms of a Markov process and a branching {\color{black}mechanism}, whether that be a branching particle system or a superprocess. We characterise Markov processes by a semigroup on $E$, denoted by $\sP=(\sP_t, t\geq0)$. We do not need $\bP$ to have the Feller property, and it is not necessary that $\bP$ is conservative. Indeed, in the case where it is non-conservative, we can append  a cemetery state $\{\dagger\}$ to $ E$, which is to be treated as an absorbing state, and regard $\bP$ as conservative on the extended space $E\cup\{\dagger\}$, which can also be treated as a Lusin space.
However, we must then alter the definition of $B(E)$ (and accordingly $B^+(E)$ and $B^+_1(E)$) to ensure that any function $f\in B(E)$ satisfies $f(\dagger)=0$.

\subsection{Non-local Branching Markov Processes}
Consider now a spatial branching process in which, given their point of creation, particles evolve independently according to a $\sP$-Markov process. In an event, which we refer to as `branching', particles positioned at $x$ die at rate {\color{black}$\beta(x)$, where }$\beta\in B^+(E)$, and instantaneously, new particles, say $N$, are created in $E$ according to a point process. The configurations of these offspring are described by the random counting measure
\[
\sZ(A) = \sum_{i = 1}^N \delta_{x_i}( A), 
\]
for Borel subsets $A$ of $E$, {\color{black} for which we also assume that $\sup_{x\in E}\mathcal{E}_x[N]<\infty$.} The law of the aforementioned point process depends on $x$, the point of death of the parent, and we denote it by $\mathcal{P}_x$, $x\in E$, with associated expectation operator given by $\mathcal{E}_x$, $x\in E$.  This information is captured in the so-called branching mechanism
\begin{equation}
  \sG[g](x) :=  \beta(x)\mathcal{E}_x\left[\prod_{i = 1}^N g(x_i) - g(x)\right], \qquad x\in E,
  \label{linearG}
\end{equation}
where $ g\in B^+_1(E)$. 
Without loss of generality, we can assume that $\mathcal{P}_x(N =1) = 0$ for all $x\in E$ by viewing a branching event with one offspring as an extra jump in the motion. On the other hand, we do allow for the possibility that $\mathcal{P}_x(N =0)>0$ for some or all $x\in E$. 

\medskip

Henceforth we refer to this spatial branching process as a $(\sP, \sG)$-branching Markov process (or $(\sP, \sG)$-BMP for short). It is well known that if the configuration of particles at time $t$ is denoted by $\{x_1(t), \ldots, x_{N_t}(t)\}$, then, on the event that the process has not become extinct or exploded, the branching Markov process can be described as the co-ordinate process $X= (X_t, t\geq0)$, given by
\[
X_t (\cdot) = \sum_{i =1}^{N_t}\delta_{x_i(t)}(\cdot), \qquad t\geq0,
\]
evolving in the space of {\color{black} counting} measures on $E$ with {\color{black} finite} total mass, which we denote $N(E)$.
In particular, $X$
is Markovian in  $N(E)$. Its probabilities will be denoted $\mathbb{P}: = (\mathbb{P}_\mu, \mu\in N(E))$. 
For convenience, we will write for any measure $\mu \in N(E)$ and function $f\in B^+(E)$,
\[
\langle f, \mu\rangle  = \int_E f(x)\mu(\dd x).
\]
In particular, 
\[
\langle f, X_t\rangle= \sum_{i = 1}^{N_t} f(x_i(t)), \qquad f\in B^+(E).
\]

With this notation in hand, it is worth noting that the independence that is inherent in the definition of the Markov branching property  implies that, if we define, 
\begin{equation*}
{\rm e}^{-\sv_t[f](x)} = \mathbb{E}_{\delta_x}\left[{\rm e}^{- \langle f, X_t\rangle}\right], \qquad t\geq 0, \, f\in B^+(E),\, x\in E,
\label{nonlin}
\end{equation*}
then for $\mu\in N(E)$, 
we have
\begin{equation}
\label{MBP}
\mathbb{E}_{\mu}\left[{\rm e}^{- \langle f, X_t\rangle}\right] = {\rm e}^{-\langle\sv_t[f], \mu\rangle }, \qquad t\geq 0.
\end{equation}
Moreover, for $f\in B^+(E)$ and  $x\in E$, 
\begin{equation}
{\rm e}^{-\sv_t[f](x)} = \sP_t[{\rm e}^{-f}](x) + \int_0^t \sP_s\left[ \sG[{\rm e}^{-\sv_{t-s}[f]}]\right](x)\dd s, \qquad t\geq0.
\label{nonlinv}
\end{equation}
 The above equation describes the evolution of the semigroup $\sv_t[\cdot]$ in terms of the action of  transport and branching. That is, either the initial particle has not branched {\color{black}and undergone a Markov transition (including the possibility of being  absorbed)} by time $t$ or at some time $s \le t$, the initial particle has branched, producing offspring according to $\sG$. We refer the reader to~\cite{SNTE-I, SNTE-II} for a proof.

Branching Markov processes enjoy a very long history in the literature, dating back as far as the late 1950s, \cite{Seva1958, Seva1961, Sk1964, INW1, INW2, INW3}, with a broad base of literature that is arguably too voluminous to give a fair summary of here. Most literature focuses on the setting of local branching. This corresponds to the setting that all offspring are positioned at their parent's point of death (i.e. $x_i = x$ in the definition of $\sG$). In that case, the branching mechanism reduces to 
\[
  \sG[s](x) = \beta(x)\left[\sum_{k =0}^\infty p_{k}(x)s^k - s\right], \qquad x\in E,
\]
where $s\in[0,1]$ and $(p_k(x), k\geq 0)$ is the offspring distribution when a parent branches at site $x\in E$.  
The branching mechanism  $\sG$ may otherwise be seen, in general, as a mixture of  local and non-local branching.

\subsection{Non-local Superprocesses} 
Superprocesses can be thought of as the high-density limit of a sequence of branching Markov processes, resulting in a new family of measure-valued Markov processes; see e.g.  \cite{ZL11, Dawson1993, Wat1968, Dynkin2, dawson2002nonlocal}. Just as branching Markov processes are Markovian in $N(E)$,  the former are Markovian in $M(E)$, the space of finite Borel measures on $E$ equipped with the {\color{black} topology of weak convergence}.  
There is a broad literature  for superprocesses, e.g. \cite{ZL11, Dawson1993, Wat1968, Alison, Janos}, with so-called local branching mechanisms, which has been  broadened to the more general setting of non-local branching mechanisms in \cite{dawson2002nonlocal, ZL11}. Let us now introduce these concepts with a self-contained  definition of what we mean by a non-local superprocess (although the reader will note that we largely conform to the presentation in \cite{ZL11}).

\medskip

A Markov process $X : = (X_t:t\geq 0)$ with state space $M(E)$  and probabilities $\mathbb{P} := (\mathbb{P}_\mu, \mu\in M(E))$ is called a $(\sP,\psi,\phi)$-superprocess  (or $(\sP,\psi,\phi)$-SP for short) if it has  semigroup $(\sV_t, t\geq 0)$ on $M(E)$ satisfying 
\begin{equation}\label{non-local-transition-semigroup}
\mathbb{E}_\mu\big[{\rm e}^{-\langle{f},{X_t}\rangle}\big] ={\rm e}^{-\langle{\sV}_t[f], \mu\rangle}, \qquad \mu\in M(E), f\in B^{+}(E),
\end{equation}
where  $({\sV}_t, t\geq 0)$ is characterised as the {unique} non-negative solution {(cf. \cite[Proposition 2.20]{ZL11})}
of the evolution equation
\begin{equation}\label{non-local-evolution-equation}
{\sV}_t[f](x)=\sP_t[f](x)-\int_{0}^{t}\sP_s\big[\psi(\cdot,{\sV}_{t-s}[f](\cdot))+\phi(\cdot,{\sV}_{t-s}[f])\big](x)\dd s.
\end{equation}
Here $\psi$ denotes the local branching mechanism
\begin{equation}\label{local-branching-mechanism}
\psi(x,\lambda)=-b(x) \lambda + c(x)\lambda^2 +\int_{0}^\infty( {\rm e}^{-\lambda y}-1+\lambda y)\nu(x,\dd y),\qquad \lambda\geq 0,
\end{equation}
where $b\in B(E)$, $c\in B^+(E)$ and 
$(y\wedge y^2)\nu(x, \dd y)$
is a uniformly (for $x\in E$) bounded kernel from $E$ to $(0,\infty)$,
and $\phi$ is the non-local branching mechanism
\begin{equation}\label{non-local-branching-mechanism}
\phi(x,f)=\beta(x)(f(x)-\eta(x,f)),
\end{equation}
where $\beta\in B^+(E)$ and $\eta$ has representation
\begin{equation*}\label{non-local-zeta-representation}
\eta(x,f)=\gamma(x,f)+\int_{M_{0}(E)}(1-{\rm e}^{-\langle{f},{\nu}\rangle})\Gamma(x,\dd \nu),
\end{equation*}
such that $\gamma(x,f)$ is a {\color{black}uniformly} bounded function on $E\times B^+(E)$ and $\langle{1},{\nu}\rangle\Gamma(x,\dd \nu)$ is a  {\color{black} uniformly (for $x\in E$)} bounded kernel from $E$ to $M_{0}(E):=M(E)\backslash \{0\}$ with
\begin{equation}
{\mathsf{m} [1](x)}:= \gamma(x,1)+\int_{M_{0}(E)}\langle{1},{\nu}\rangle\Gamma(x,\dd \nu)\leq 1.
\label{leq1}
\end{equation}

We refer the reader to \cite{dawson2002nonlocal, PY} for more details regarding the above formulae. 
Lemma 3.1 in \cite{dawson2002nonlocal} tells us that the functional $\eta(x,f)$ has the following equivalent representation
\begin{equation}\label{zeta-representation}
\eta(x,f)=\int_{M_0(E)}\bigg[\delta_\eta(x,\pi)\langle{f},{\pi}\rangle+\int_{0}^{\infty}(1-{\rm e}^{-u\langle{f},{\pi}\rangle})n_\eta(x,\pi,\dd u)\bigg]P_\eta(x,\dd\pi),
\end{equation}
where $M_0(E)$ denotes the set of probability measures on $E$, $P_\eta(x,\dd\pi)$ is a probability kernel from $E$ to $M_0(E)$, $\delta_\eta\geq 0$ is a bounded function on $E\times M_0(E)$, and $un_\eta(x,\pi,\dd u)$ is a bounded kernel from $E\times M_0(E)$ to $(0,\infty)$ with
\begin{equation*}
\delta_\eta(x,\pi)+\int_{0}^{\infty}un_\eta(x,\pi,\dd u)\leq 1.
\end{equation*}

The reader will note that we have deliberately used some of the same notation for both branching Markov processes and  superprocesses. In the sequel there should be no confusion and the motivation for this choice of repeated notation is that our main results are indifferent to which of the two processes we are talking about.

\section{Assumptions and main results}
Before stating our results, we first introduce some assumptions that will be crucial in analysing the models defined above. Unless a specific difference is indicated, the assumptions apply both to the setting that $X$ is either a non-local branching Markov process or a non-local superprocess.

\begin{description}

\item[\namedlabel{H1}{(H1)}]
In the particle setting, we assume first moments
\begin{equation*}
	\sup_{x\in E} \mathcal{E}_x[N]<\infty.
\end{equation*}

\end{description}

\begin{remark}\label{rem1} 
In the setting of superprocesses, since we have already assumed that the kernels $(y \wedge y^2)\nu(x, \dd y)$ and $\langle 1, \nu \rangle \Gamma(x, \dd \nu)$ are uniformly bounded (for  $x\in E$), and that condition \eqref{leq1} holds, no additional assumption is required. 
\end{remark}

\begin{description}
\item[\namedlabel{H2}{(H2)}]
Define the expectation semigroup $(\sT_t, t \ge 0)$ via $\sT_t[f](x) = \mathbb{E}_{\delta_x}[\langle f, X_t\rangle]$, $t\geq 0$, $f\in B^+(E)$. There exist 
a function $\varphi \in B^+(E)$ {strictly positive} and finite measure 
{$\tilde\varphi\in M(E)$} 
such that, for $f\in B^+(E)$,
\begin{equation}
\langle{\sT_t[f] },{\tilde\varphi} \rangle=\langle f, \tilde{\varphi}\rangle
\ \text{ and } \ 
\langle \sT_t[\varphi] , \mu\rangle = \langle{\varphi},{\mu}\rangle,
\label{criticaldef}
\end{equation}
for all $\mu\in N(E)$ (resp. $M(E)$) if $(X, \mathbb{P})$ is a branching Markov process (resp. a superprocess),  where
    \[
        \langle\textsf{T}_t[f], \mu\rangle = \int_{E}\mu(\dd x)\mathbb{E}_{\delta_x} \left[ \langle f, X_t\rangle\right]= \mathbb{E}_{\mu} \left[ \langle f, X_t\rangle\right], \qquad t\geq 0.
    \]
Further let us define
 \begin{equation}
 \Delta_t := \sup_{x\in E,\, f\in B^+(E),\, \|f\|=1}\left|{(\langle f,\tilde\varphi}\rangle\varphi(x))^{-1}\sT_t\bra{f}(x)-1\right| , \qquad {t>0}.
\label{Deltat}
 \end{equation}
 We suppose that 
 {$\Delta_t<\infty$ for all $t>0$ and}
$
\lim_{t\to\infty} \Delta_t=0.
$
Without loss of generality, we conveniently impose the normalisation $\langle\varphi,\tilde{\varphi}\rangle=1$.
 \end{description}

\begin{remark}
    The  requirement \eqref{criticaldef} necessarily identifies the process $X$ (whether a particle process or a superprocess) as critical as $\varphi$ and $\tilde\varphi$ are associated with a unit eigenvalue with respect to the mean semigroup $(\textsf{T}_t,t\geq0)$. For further details of this classification, see for example \cite{Horton2023, AH1, AH2, AH3}.
\end{remark}

\begin{remark} 
    Note that in the definition of $\Delta_t$, we take the supremum over functions $f \in B^+(E)$ with $\|f\| = 1$. This is not a restrictive condition, as the supremum involves the ratio $(\langle f,\tilde\varphi\rangle\varphi(x))^{-1}\sT_t\bra{f}(x)$, which is invariant under positive scalar multiplication of $f$. The normalization simply avoids the case $f \equiv \mathbf{0}$, for which the ratio is not well-defined. Equivalently, we could take the supremum over $B^{+}_1(E)\setminus \{\textbf{0}\}$.
\end{remark}

\begin{description}
\item[\namedlabel{H3}{(H3)}]
For each $x\in E$ 
\[
\mathbb{P}_{\delta_x}(\zeta < \infty) = 1,
\]
where $\zeta = \inf \{t>0 : \langle 1, X_t\rangle = 0\}.$

\medskip

\item[\namedlabel{H4}{(H4)}]
 For  $g\in B_1^+(E)$ and $x\in E$, define
\begin{equation*}
    \textsf{A}[g](x):= \beta(x)\mathcal{E}_x \left[\prod_{i=1}^N (1-g(x_i))-1 + \sum_{i=1}^N g(x_i)\right].
\end{equation*}
We assume that, 
\begin{equation*}
	\langle \textsf{A}[x\varphi], \tilde{\varphi} \rangle = x^{1+\alpha}\ell(x),\text{ as }x\to0,
\end{equation*}
where $\alpha\in (0,1]$ and $\ell$ is a slowly varying function at 0.
In the superprocess setting, define 
\begin{align}\label{J-H5}
\sJ[h](x) &:= c(x)h(x)^2 +\int_{(0,\infty)}( {\rm e}^{-h(x) y}-1+h(x) y)\nu(x,\dd y)\\
&\hspace{1cm} +\beta(x)\int_{M_{0}(E)}({\rm e}^{-\langle{h},{\nu}\rangle}-1+ \langle h, \nu\rangle)\Gamma(x,\dd \nu). \notag
\end{align} 
Then we assume instead 
\begin{equation*}
	\langle \textsf{J}[x\varphi], \tilde{\varphi} \rangle = x^{1+\alpha}\ell(x),\text{ as }x\to0,
\end{equation*}
with the same requirements of $\alpha$ and $\ell$. 
\end{description}

\begin{remark} 
In  \cite{horton2024stability}, for the the classical Kolmogorov and Yaglom limits in a similar setting of non-local spatial branching particle processes and superprocesses the key assumption in place of \ref{H4} was given as  
 $\sup_{x\in E}\mathcal{E}_x[N^2]<\infty$; equivalently, in the superprocess setting, 
\[
\sup_{x\in E}\left( \int_{0}^\infty y^2 \nu (x, \dd y)  + \beta(x)\int_{M(E)^\circ} \langle 1, \nu\rangle^2 \Gamma(x,\dd\nu)\right)<\infty.
\]
This corresponds to the setting that $\alpha = 1$, in which case $\langle \textsf{A}[x\varphi], \tilde{\varphi} \rangle = x^{2}\ell(x) = \langle \textsf{J}[x\varphi], \tilde{\varphi} \rangle $ as $x\to0$, where $\lim_{x\to 0}\ell(x)=\langle \frac12 \mathbb{V}[\varphi],\tilde{\varphi}\rangle$ {with}, for $f\in B^+(E)$
\begin{equation*}
    \mathbb{V}[f](x) = \beta(x)\mathcal{E}_x\left[\sum_{i, j= 1;\, i\neq j}^N f(x_i)f(x_j)\right], \qquad x\in E,
\end{equation*}
in the branching Markov process setting and
\begin{align*}
\mathbb{V}[f](x)  &= \psi''(x,0+)f(x)^2+ \beta(x)\int_{M(E)^\circ} \langle f, \nu\rangle^2 \Gamma(x,\dd\nu),\\ 
&= \left(2c(x) +  \int_{0}^\infty y^2 \nu (x, \dd y)\right)f(x)^2 + \beta(x)\int_{M(E)^\circ} \langle f, \nu\rangle^2 \Gamma(x,\dd\nu),
\end{align*}
in the superprocess setting. 

On the other hand, we believe that the case $\alpha=0$ requires a different approach from the one used in this paper. This is supported by the techniques observed in Theorem 2 of \cite{nagaev2007critical} for Bienaym\'e--Galton--Watson processes, where a different normalization appears in the Yaglom limit. To remain consistent with the universality of such asymptotic results, we expect a similar phenomenon to arise in the more general settings of non-local spatial branching processes and superprocesses.
 \end{remark}

\begin{description}
\item[\namedlabel{H5}{(H5)}]
There exists $\delta \in (0,1)$ such that
\[
\sup_{x\in E} \frac{\beta(x)}{\varphi(x)} \mathcal{E}_x \left[(N-1)^\delta\sum_{i=1}^N \varphi(x_i)  \right]
< \infty
\]
for branching Markov processes, and
\[
\sup_{x\in E} \left(\int_{(1,\infty)} y ^{1+\delta} \nu(x,\mathrm{d} y) + \frac{\beta(x)}{\varphi(x)} \int_{M_0(E)} (1\vee \langle1,\nu \rangle^\delta) \langle \varphi, \nu \rangle \Gamma (x, \mathrm{d} \nu)\right) < \infty
\]
for superprocesses.
\end{description}

\begin{remark} 
In some sense, the above condition can be interpreted as the existence of a uniform bound for a ($1+\delta$) moment related to the eigenfunction $\varphi$.
\end{remark}

\begin{remark}\label{rem-assumptionH5}
In Section 3 of  \cite{HeringHoppe}, a condition analogous to \ref{H4} is imposed in the context of non-local branching processes. In their setting, additional assumptions are required, which in our notation can be expressed as follows: there exist positive constants $c$ and $c^*$ such that 
\begin{equation}\label{eq:HeringHoppe}
\beta \mathsf{m}[\varphi] \le c \varphi \quad \text{and}\quad \langle \beta \mathsf{m}[f], \tilde{\varphi}\rangle \le c^* \langle f, \tilde{\varphi} \rangle, \quad \text{for any}\ f\in B^+(E),
\end{equation}
where $\beta$ is 
branching rate and
\[\mathsf{m}[f](x) := \mathcal{E}_x \left[\sum_{i=1}^N f(x_i)\right], \quad f \in B^+(E), \ x\in E.\]
The first inequality in \eqref{eq:HeringHoppe} can be rewritten as
\[\sup_{x\in E}\frac{\beta(x)}{\varphi(x)} \mathcal{E}_x \left[\sum_{i=1}^N \varphi(x_i)\right] \leq c,\]
which, on its own, is actually weaker than our assumption~\ref{H5}. However, the second condition in~\eqref{eq:HeringHoppe} does not have a direct analogue in our framework. Roughly speaking, it imposes an averaged boundedness under the  measure~$\tilde{\varphi}$, ensuring that the total expected offspring mass, weighted by~$\beta$, remains uniformly controlled. In contrast, our assumption~\ref{H5} replaces this averaged bound by a moment condition of order~$1+\delta$ for some $\delta \in (0,1)$. Therefore, instead of needing control on both the eigenfunction and the eigenmeasure as in \eqref{eq:HeringHoppe}, we require the $1+\delta$ condition only on the eigenfunction. 
\end{remark}

{
Finally, for the setting of superprocesses, we require one additional assumption to rigorously justify the existence of the Dynkin-Kuznetsov excursion measure $\mathbb{N}_x$ associated with $X$ that we will use later in Theorem \ref{superspine} (the spine decomposition). Let 
$\{Q_t(\mu, \cdot) := \mathbb{P}_\mu(X_t \in \cdot) : t \ge 0, \mu \in M(E)\}$ 
denote the transition kernel of the superprocess. 
From \eqref{non-local-transition-semigroup}, we have
\begin{equation*}
\int_{M(E)} {\rm e}^{-\langle f, \nu \rangle} Q_t(\mu, \dd\nu) = \mathrm{e}^{-\langle \sV_t[f], \mu \rangle}, \qquad \mu \in M(E), \quad t \ge 0.
\end{equation*}
Consequently, the branching property implies that $Q_t(\mu, \cdot)$ is an infinitely divisible probability measure on $M(E)$. Thus, the non-linear semigroup admits the representation
\begin{equation*}
\sV_t[f] = \lambda_t(x, f) + \int_{M_0(E)} (1 - {\rm e}^{-\langle f, \nu \rangle}) L_t(x, \dd\nu), \qquad t > 0, \quad x\in E,  \quad f \in B^+(E), 
\end{equation*}
where $\lambda_t(x, \dd y)$ is a bounded kernel on $E$ and $(1 \wedge \langle 1, \nu \rangle) L_t(x, \dd\nu)$ is a bounded kernel from $E$ to $M_0(E)$. Let us define the set where the linear term vanishes as
$$E_0 := \{x \in E : \lambda_t(x, E) = 0 \text{ for all } t > 0\}.$$
For $x \in E_0$, the family of measures $\{L_t(x, \cdot) : t > 0\}$ constitutes a valid entrance law, which ensures the existence of a unique $\sigma$-finite measure $\mathbb{N}_x$ (the Dynkin-Kuznetsov measure) on the c\`adl\`ag path space of $M(E)$ such that $\mathbb{N}_x(\{0\}) = 0$ and, for any $0<t_1<t_2<\cdots <t_n<\infty$ and $\nu_1, \dots, \nu_n \in M_0(E)$,
\[\mathbb{N}_x(\omega_{t_1}\in \dd \nu_1, \ldots, \omega_{t_n}\in \dd \nu_n)= L_{t_1}(x,\dd \nu_1)Q_{t_2 - t_1}(\nu_1, \dd \nu_2)\cdots Q_{t_n - t_{n-1}}(\nu_{n-1}, \dd \nu_n),\]
where  $\omega = (\omega_s, s\geq0)$ is the càdlàg path on   ${M}(E)$.
Further, it satisfies
\begin{equation*}
\mathbb{N}_x(1 - {\rm e}^{-\langle f, \omega_t \rangle}) = \sV_t [f](x), \qquad t > 0, \quad f \in B^+(E).
\end{equation*}
We refer the reader to \cite{dynkin2004measures} and Section 8.4 of \cite{ZL11} for details on the Dynkin--Kuznetsov measure. As we will see later, in the spine decomposition (Theorem \ref{superspine}), 
continuous immigration occurs along the spine at certain rate which includes a product of the coefficient $c$ and the Dynkin--Kuznetsov measure.
To guarantee that $\mathbb{N}_x$ is well-defined at all locations where this immigration rate is strictly positive, we require the following condition, analogous to Assumption 0 in \cite{RenSongYang2022}.
\begin{description}
		\item[\namedlabel{H6}{(H6)}] In the superprocess setting, we assume that $$E_c := \{x \in E : c(x) > 0\} \subset E_0 .$$
\end{description}
\begin{remark}\label{remH6}
In general, it is worth noting that \ref{H6} is automatically satisfied for a large class of branching mechanisms. 
Indeed, according to \cite[Corollary 5.33]{ZL11}, a sufficient condition for $E_0 = E$ (which trivially implies $E_c \subset E_0$) is the existence of a local branching mechanism 
$$
\psi_\ast(\lambda) = b_\ast \lambda + c_\ast \lambda^2 +\int_{0}^\infty ( {\rm e}^{-\lambda y}-1+\lambda y)\nu_\ast(\dd y),\qquad \lambda\geq 0,
$$
where $b_\ast \in \mathbb{R}$, $c_\ast \ge 0$, and $(y \wedge y^2)\nu_\ast (\dd y)$ is a bounded kernel on $(0, \infty)$, such that  
\begin{equation}\label{sufcondH6}
\lim_{\lambda \to \infty} \psi_\ast^\prime(\lambda) = \infty \qquad \text{and}\qquad 
\Phi_1(x,\lambda) \ge \psi_\ast(\lambda) \quad \text{for all } x \in E \text{ and } \lambda \geq 0,
\end{equation}
where $\Phi_1$ is the 
projection of the branching mechanism $\psi(x, f(x)) + \phi(x, f)$ and it is given by
$$\Phi_1(x,\lambda) = \psi(x, \lambda) + \beta(x) \left(1-\mathsf{m} [1](x) \right) \lambda + \beta(x) \int_{M_0(E)}  \left({\rm e}^{-\langle \lambda \mathbf{1}_{\{x\}}, \nu\rangle }-1+\langle \lambda \mathbf{1}_{\{x\}}, \nu\rangle \right) \Gamma(x, \dd \nu),$$
where $\mathsf{m} [1](x)$ is defined in \eqref{leq1}.
\end{remark}
}

We are now ready to state our main results. The reader will note that the results are stated for both branching Markov processes and superprocesses simultaneously. Moreover, as alluded to in the introduction, we give greater emphasis on giving the full details for the BMP proofs, offering only adjustments to the arguments to deal with the superprocess setting where necessary.
The first result  is a version Kolmogorov's asymptotic for the survival probability asymptotic holds, albeit now for processes with infinite variation. 

\begin{theorem}[Kolmogorov's estimation]\label{theo:kolmogorov}
Suppose that $(X, \mathbb{P})$ is a $(\sP, \sG)$-BMP (resp.  a $(\sP,\psi,\phi)$-SP) satisfying \ref{H1}-\ref{H5} {(and additionally \ref{H6} for $(\sP,\psi,\phi)$-SP)}. Then, for all $\mu \in N(E)$  (resp. $\mu\in M(E)$),
\begin{equation*}
	\lim_{t\to \infty} \frac{\mathbb{P}_\mu (\zeta>t)}{t^{-1/\alpha} \tilde{\ell}(t)} = \langle \varphi,\mu \rangle, 
\end{equation*}
where $\tilde{\ell}$ is a slowly varying function at $\infty$. More precisely, $\tilde{\ell}(t) = 1 / L^\ast (t^{1/\alpha})$ where $L^*$ is the Bruijn conjugate of $s\mapsto \alpha^{-1/\alpha} \ell(1/s)^{-1/\alpha}$.
\end{theorem}

The next result is the analogue of Yaglom's classical conditional critical convergence theorem, albeit now for the setting of infinite variance. 
The reader will note that Kolmogorov's estimate and Yaglom's theorem are asymptotically compatible: the survival probability decays like $t^{-1/\alpha} \tilde{\ell}(t)$, while the conditional population size grows as $t^{-1/\alpha} \tilde{\ell}(t)$ and the resulting distribution limit is heavy-tailed.  

\begin{theorem}[Yaglom's limit]\label{theo:Yaglom}
Suppose that $(X, \mathbb{P})$ is a $(\sP, \sG)$-BMP (resp. a $(\sP,\psi,\phi)$-SP) satisfying \ref{H1}-\ref{H5} {(and additionally \ref{H6} for $(\sP,\psi,\phi)$-SP)}.  Then, for all $\mu \in N(E)$ (resp. $\mu\in M(E)$),
\begin{equation*}
	\lim_{t\to \infty} \mathbb{E}_\mu\left[\exp\left(- t^{-1/\alpha} \tilde{\ell}(t) \langle f, X_t \rangle \right) \Big|\Big. \ \langle 1, X_t \rangle >0 \right] = 1- \frac{\langle f, \tilde{\varphi} \rangle }{(1+ \langle f, \tilde{\varphi} \rangle^{\alpha})^{1/\alpha}},
\end{equation*}
where 
$f\in B^+(E)$.
\end{theorem}

We would like to point out that our results are consistent with the finite second moment case given in Theorem 1 and 2 of \cite{horton2024stability}.  Since, in this case,  $\alpha=1$ and we choose $\ell$ such that $\lim_{t\to0} \ell(t) = \langle \frac12 \mathbb{V}[\varphi],\tilde{\varphi}\rangle$, it follows that $L(t) := \ell(1/t)^{-1} \to \langle \frac12 \mathbb{V}[\varphi],\tilde{\varphi}\rangle^{-1}$ as $t\to\infty$. Recalling that the Bruijn conjugate satisfies $\lim_{t\to\infty} L(t L^\ast(t))L^\ast (t) = 1$, we then conclude that $\lim_{t\to\infty} L^\ast(t) = \langle \frac12 \mathbb{V}[\varphi],\tilde{\varphi}\rangle$.

\begin{remark}
{From Theorem \ref{theo:Yaglom} and appealing to Theorems 1.18 and 1.19 in \cite{ZL11}, 
%
the random measure $t^{-1/\alpha} \tilde{\ell}(t) X_t$ converges weakly under $\mathbb{P}_\mu\left(\cdot \mid \langle 1, X_t\rangle >0\right)$ as $t\to\infty$ to the limiting measure $W\tilde\varphi$ , where $W$ is a random variable with 
a law whose Laplace transform is given by
		$$L_W(\theta) = \mathbb{E} \left[\mathrm{e}^{-\theta W}\right] = 1 - \frac{\theta}{(1+\theta^\alpha)^{1/\alpha}}, \qquad \theta \geq 0.$$
Clearly, to see if $W\tilde\varphi$ is infinitely divisible, it is enough to prove that $W$ is infinitely divisible. Indeed,
	by integration by parts, it is not difficult to see that
$$\int_0^\infty \mathrm{e}^{-\theta x} \mathbb{P}(W > x) \dd x = \frac{1 - L_W(\theta)}{\theta} = (1+\theta^\alpha)^{-1/\alpha}.$$
Thus $(1+\theta^\alpha)^{-1/\alpha}$ is the Laplace transform of the tail function $x\mapsto \mathbb P(W>x)$. We now show that this tail function is completely monotone. Since
$\theta\mapsto\theta$ is a complete Bernstein function, and the constant function $1$ is also a complete
Bernstein function, Corollary 7.15 in \cite{schilling2012bernstein}
implies that $\theta\mapsto (1+\theta^\alpha)^{1/\alpha}$
is a complete Bernstein function for $\alpha \in (0,1]$. Therefore, by Theorem 7.3 in
\cite{schilling2012bernstein},
\[
\theta\mapsto (1+\theta^\alpha)^{-1/\alpha}
\]
is a Stieltjes function. Now, appealing to Theorem 2.2(i) in \cite{schilling2012bernstein}, every Stieltjes
function is the Laplace transform of a measure whose density on
$(0,\infty)$ is completely monotone. Hence there exists a completely
monotone function $g$ such that
\[
(1+\theta^\alpha)^{-1/\alpha} =\int_0^\infty e^{-\theta x} g(x) \dd x, \qquad  \theta>0.
\]
By uniqueness of Laplace transforms, we obtain $g(x)=\mathbb P(W>x), \ x\ge0.$
Therefore the tail function $x\mapsto \mathbb P(W>x)$ is completely
monotone. Since the density $f_W=-g'$ of $W$ is the negative derivative of a completely monotone function, $f_W$ is completely monotone. Finally, by the Goldie-Steutel theorem (cf. \cite[Theorem 1]{Steutel}), every non-negative random variable
with a completely monotone density is infinitely divisible. Hence $W$ is
infinitely divisible, and therefore so is $W\widetilde{\varphi}$.}

{Finally, we give the Lévy-Khintchine representation of $W\tilde\varphi$. 
Its Lévy measure $\Lambda$ is supported strictly on the one-dimensional ray $\{ r \tilde\varphi : r > 0 \}$ within the space of measures, so that
$$\mathbb{E}\left[\mathrm{e}^{-\langle f, W\tilde\varphi\rangle }\right] = \exp\left( - \int_{M_0 (E) } (1 - \mathrm{e}^{-\langle f,\mu\rangle}) \Lambda(\dd\mu)\right)= \exp\left( - \int_0^\infty (1 - \mathrm{e}^{-r \langle f, \tilde\varphi\rangle}) \Pi(\dd r) \right),$$ 
where $\Pi(\dd r)$ is the one-dimensional Lévy measure of the scalar variable $W$ on $(0, \infty)$ and is uniquely determined by the Laplace exponent of $W,$ 
$$\int_0^\infty (1 - \mathrm{e}^{-\theta r}) \Pi(\dd r) = -\ln\left( 1 - (1+\theta^{-\alpha})^{-1/\alpha} \right).$$}
		\end{remark}

\section{Discussion of results {and examples}}\label{sec:examples}
In this section, we spend some time discussing the consistency of our results with the existing literature. Moreover, we also take the opportunity to discuss assumptions \ref{H1}-\ref{H6} in the setting of some specific processes.

 In \cite{horton2024stability} a full treatment of the universal nature of the Yaglom limit for non-local branching processes with and without immigration was treated under a natural second moment assumption together with a Perron--Frobenius type assumption similar to \ref{H2}. The starting point of this study has therefore been   to understand what the natural regular variation replacement assumption for second moments should be that is consistent with existing literature. Moreover, our  results and proofs have been achieved by taking a holistic view across a wide variety of literature spanning five decades. In this respect, \ref{H4} agrees, for example, with the assumptions of \cite{zolotarev1957, slack1968branching} in the setting of Galton--Watson processes, \cite{Goldstein_Hoppe, vatutin1977, vatutin2021} in the setting of multi-type Galton--Watson processes, \cite{RenYangZhao2014} for continuous-state branching processes, Chapter VI.4 of \cite{AH3} and \cite{HeringHoppe} for non-local branching Markov processes and \cite{ren2020limit} for a relatively particular class of superprocesses.

The approach we take adopts and adapts methodologies from across the aforementioned literature, in particular, we use a spine approach modified from \cite{ren2020limit} for the Kolmogorov's estimation (Theorem \ref{theo:kolmogorov}) and the one in Chapter VI.4 of \cite{AH3} for the Yaglom limit (Theorem \ref{theo:Yaglom}).  
In the setting of branching Markov processes, our results are consistent with those of \cite{HeringHoppe}, albeit  our proofs are based on a spine change of measure, which results in a  difference in the assumptions, as earlier mentioned   in Remark \ref{rem-assumptionH5}. In the setting of non-local superprocesses, we offer the first general result of this type, where the previous work of \cite{ren2020limit} focused on a specific category of  local branching mechanism.

Let us conclude this section by offering several concrete examples which our general results cover and, in particular, some of them not currently carried by the literature. 

\bigskip

\noindent {\bf Multi-type continuous-time Bienaym\'e--Galton--Watson processes:} Let us consider a continuous-time Bienaym\'e--Galton--Watson process with a countable set of types and infinite second moments for the offspring distribution of some or all types (the counterpart model in the discrete-time setting has been studied recently in \cite{vatutin2021}). The continuous-time processes can be viewed as a BMP, where the non-locality occurs in the set of types $E=\mathbb{N}$, but there is no spatial motion  ($\sP$ is the identity operator). 
The branching mechanism defined in \eqref{linearG} takes the form
\[
	\sG [g] (i) 
	= 
	\beta (i) \mathcal{E}_i \left[\prod_{j \in \N} g(j)^{N_j} - g(i)\right] 
	= 
	\beta (i) \left[ \sum_{\mathbf{k}\in \N_0^\N} \mathcal{P}_i (\mathbf{N}=\mathbf{k})  \prod_{j\in\N} g(j)^{k_j} - g(i)\right], \qquad i\in \N,
\]
where $\beta\in B^+(\N)$, $g\in B_1^+(\N)\simeq [0,1]^\N$ and $\mathbf{N}=(N_1, N_2, \ldots)$. Note that $N_j$ represents the particles of type $j$ produced in a branching event, for $j\in \mathbb{N}$. We refer the reader to \cite[Chapter V.7]{AN} for the setting with finite number of types to see the parallelism. Assuming $\sup_{i\in \N} \mathcal{E}_i \left[N\right] <\infty$ where $N= \sum_{j\in\N} N_j$, i.e. \ref{H1}, the offspring mean matrix
\[
 (\sm)_{i,j} = \mathcal{E}_i [N_j], \qquad i,j\in \N,
\]
is well defined. 
Finally, we define our generator in this framework, that is, the matrix
\begin{equation}\label{Lmatrix}
(L)_{i,j} = \beta(i) \mathcal{E}_i [N_j-\delta_{i=j}], \qquad i,j\in \N.
\end{equation}
Let $(\sT_t, t\ge 0)$ denote the mean semigroup, then 
$\sT_t[f] (i) = \mathrm{e}^{tL} [f] (i)$ according to \eqref{eq:BMP-linear} (see the equivalence with  \cite[p. 202]{AN} for the finite type case).
Now we assume, as in Theorem D of 
Vere-Jones \cite{VJ}, that  $\mathrm{e}^L$ is irreducible and 1-positive. This ensures that {there exist the right-eigenvector  $\varphi = (\varphi(j), j\in \mathbb{N})$ and left-eigenvector $\tilde\varphi= (\tilde \varphi(j), j\in \mathbb{N})$} of $\mathrm{e}^L$ associated to the maximal eigenvalue which is equal to 1. Furthermore, $\varphi$ and $\tilde\varphi$ are also eigenvectors for $L$ but in this case associated to the eigenvalue 0. 
{{Now note that} $\Delta_t$ can be expressed as 
\[\Delta_t = \sup_{i\in \mathbb{N}} \sup_{(\xi_1,\xi_2,\ldots)\in  \mathbb{R}_+^\N : \sup_{j\in \mathbb{N}} \xi_j = 1} \left|\frac{\sum_{j\in \mathbb{N}}M_{ij}(t) \xi_j}{\varphi(i)\sum_{j\in \mathbb{N}} \xi_j \tilde{\varphi}(j)} -1 \right|,\]
where
\[M_{ij}(t)= {\sT_t [\mathbf{1}_{\{j\}}] (i)} =\mathbb{E}_{\delta_i} \left[\sum_{k=1}^{N_t}\mathbf{1}_{\{j\}}(x_k(t))\right], \qquad i,j \in \mathbb{N}.\]
From Theorem D in \cite{VJ}, we also have that for each $i, j \in \mathbb{N}$
\[M_{ij}(t) \to \varphi(i)\tilde{\varphi}(j), \quad \text{as}\quad t\to \infty.\]
 See also the discussion around equations (1.10)--(1.11) in \cite{vatutin2021}. However, this pointwise convergence is not sufficient to verify condition \ref{H2}, which requires uniform convergence in order to ensure that $\Delta_t \to 0$ as $t\to\infty$. We therefore impose this as an additional condition, i.e. we assume
\[\sup_{i,j\in \mathbb{N}} \left|\frac{M_{ij}(t)}{\varphi(i)\tilde\varphi(j)} -1 \right| \to 0, \quad \text{as}\quad t\to \infty
.\] 
  Note that, in the discrete-time setting, Vatutin et al.~\cite{vatutin2021} also impose a kind of uniform convergence assumption on the offspring mean matrix (see \cite[Definition 2]{vatutin2021}).
 } 


Regarding Assumption \ref{H3}, it would follow from the counterpart result of \cite[Lemma 3]{vatutin2021} in the continuous-time setting. For simplicity, here we directly assume the claim of \ref{H3}. Hypotheses \ref{H4} and \ref{H5} are also required.

To conclude the example, we would like to point out that if we consider a finite number of types, i.e. $E=\{1,\dots,n\}$, {then the uniform convergence described above is guaranteed and {\ref{H2} is nothing more than the classical Perron-Frobenius asymptotic.}} Further, in this case, it is known that \ref{H3} also holds and we can show that \ref{H5} is a consequence of \ref{H4}. To see this, note first that \ref{H4} can be written as follows
\begin{equation*}
	\langle \textsf{A}[x\varphi], \tilde{\varphi} \rangle 
	= 
	\sum_{i=1}^n \tilde{\varphi}(i)  \beta(i)\mathcal{E}_i \left[\prod_{j=1}^n \big(1-x\varphi(j)\big)^{N_j} -1 +x \sum_{j=1}^n N_j  \varphi(j)\right]
	= x^{1+\alpha}\ell(x), \text{ as }x\to0.
\end{equation*}
Using the Taylor expansion of the logarithm and the exponential, we have
\begin{align*}
	&\mathcal{E}_i \left[\prod_{j=1}^n \big(1-x\varphi(j)\big)^{N_j}\right] 
	=
	\mathcal{E}_i \left[\mathrm{e}^{\sum_{j=1}^n N_j \ln (1-x\varphi(j))}\right] 
	=
	\mathcal{E}_i \left[\mathrm{e}^{\sum_{j=1}^n N_j  (-x\varphi(j)- x^2 \varphi(j)^2 / 2 - O(x^3))}\right] 
	\\
	&\quad =
	\mathcal{E}_i \left[\mathrm{e}^{-x\sum_{j=1}^n N_j  \varphi(j)}\left(1-O\big(x^2 \sum_{j=1}^n N_j  \varphi(j)^2\big)\right)\right] 
	=
	\mathcal{E}_i \left[\mathrm{e}^{-x\sum_{j=1}^n N_j  \varphi(j)}\right] - O(x^2).
\end{align*}
Thus, \ref{H4} implies
\begin{equation*}
	\sum_{i=1}^n \tilde{\varphi}(i)  \beta(i) \left(\mathcal{E}_i \left[\mathrm{e}^{-x\sum_{j=1}^n N_j  \varphi(j)}\right] - 1 + x  \mathcal{E}_i \left[\sum_{j=1}^n N_j  \varphi(j)\right]\right)
	=
	x^{1+\alpha}\ell(x), \text{ as }x\to0.
\end{equation*}
Applying a Tauberian Theorem, it is possible to derive that the tail $\mathcal{P}_i(\sum_{j=1}^n N_j  \varphi(j) > y)$ is a regularly varying function with index $-\alpha$ as $y\to\infty$, for each  $i\in\{1,\ldots,n\}$. This is enough to get
\[
\mathcal{E}_i \left[N^{1+\delta}  \right]<\infty , \qquad i\in\{1,\ldots,n\}, \quad \delta \in (0, \alpha),
\]
concluding the proof of \ref{H5}. 

\bigskip

\noindent{\bf Branching diffusion on a bounded domain: }In this model, particles move through the space $D\subset\mathbb{R}^n$ according to a diffusion whose generator $L = \frac{1}{2} \sum_{i,j=1}^n \partial_{x_j}(a^{ij} \partial_{x_i})$ is uniformly elliptic with coefficients $a^{i j} = a^{ ji}\in C^1 (\bar{D})$ for all $1\le i,j \le n$, and they branch independently at a constant rate $\beta > 0$ with the number of offspring being identically and independently distributed like $N$. Branching is local, meaning the offspring are positioned at the location where their parent dies. This process  
examined by \cite{Ellen}, occurs within a open connected bounded domain, such that when an individual first touches its boundary $\partial D$, it is killed and sent to a cemetery state. Powell demonstrated  in  \cite{Ellen} that \ref{H2} is valid as long as $\partial D$ is Lipschitz. 
Indeed, using the so-called many-to-one formula, we have
\[\mathsf{T}_t [f] (x) = \mathrm{e}^{(\mathcal{E}[N]-1)\beta t}\mathsf{P}_t[f](x) = \mathrm{e}^{(\mathcal{E}[N]-1)\beta t}\int_D f(y) p_t^D(x,y)\dd y,\]
where $p_t^D(x,y)$ is the transition density of the diffusion process. Let us denote by $0 < \lambda_1 < \lambda_2 \leq \lambda_3 \cdots$ the eigenvalues of $-L$ and by $\{\varphi_k\}_{k\geq 1}$ the associated eigenfunctions. 
At criticality, it is necessarily the case that $\lambda_1 = (\mathcal{E}[N]-1)\beta$, $\varphi=\varphi_1$ and \(\tilde\varphi(\mathrm{d}x) = \varphi_1(x)\mathrm{d}x\), where we suppose that $1<\mathcal{E}[N]<\infty$, i.e. \ref{H1} holds.
Applying Lemma 2.1 of \cite{Ellen}, for any $\tau>0$ there exists a constant $C_\tau>0$ depending only on the domain such that  
\[
\left| \frac{\mathsf{T}_t [f] (x)}{\langle f, \tilde\varphi\rangle \varphi(x)} -1 \right|
=
\left| \int_D \frac{f(y)\varphi_1(y)}{\langle f, \tilde\varphi\rangle}\left[\frac{\mathrm{e}^{\lambda_1 t}p_t^D(x,y)}{\varphi_1(x)\varphi_1(y)} -1 \right] \dd y \right|
\leq 
C_\tau \mathrm{e}^{-(\lambda_2-\lambda_1)t},
\]
for all $t > \tau$, $x \in D$ and $f\in B_1^+(E)$ with $\|f\|=1$. 
Additionally, it is also established that for critical systems \ref{H3} is satisfied, see Theorem 1.3 in aforesaid reference. 
Now we assume 
\begin{equation}\label{eq:pgfN}
\mathcal{E}[z^N] = 1 -(1-z) \mathcal{E}[N] + (1-z)^{1+\alpha} \bar\ell(1-z), \qquad z\in [0,1],
\end{equation}
where $\alpha \in (0,1]$ and  $\bar\ell$ is a slowly varying function at $0$; note this  is the counterpart of \eqref{eq:genh}. It follows from the definition of $\mathsf{A}$ and the previous equality that
\begin{align*}
	\langle \mathsf{A}[y\varphi],\tilde\varphi\rangle 
	&= 
	\beta \int_D  \mathcal{E} \left[(1-y\varphi(x))^N-1+Ny\varphi(x)\right] \varphi(x)\mathrm{d}x
	\\
	&=
	\beta \int_D  (y\varphi(x))^{1+\alpha}\bar\ell(y\varphi(x))\varphi(x)\mathrm{d}x
	\\
	&=
	y^{1+\alpha} \ell(y),
	\end{align*}
	where $\ell(y) = \beta \int_D  \varphi(x)^{2+\alpha}\bar\ell(y\varphi(x))\mathrm{d}x$ is a slowly varying function at 0. Indeed, 
	since  $\varphi(x)>0$ for all $x\in D$ (see page 6 of \cite{Ellen}), and $\bar{\ell}$ is a slowly varying function at 0, we deduce that for any $c>0$,
	\[
	\lim_{y\to 0} \frac{\ell(cy)}{\ell(y)} =
	\lim_{y\to 0} \frac{\int_D  \varphi(x)^{2+\alpha}\frac{\bar\ell(cy\varphi(x))}{\bar\ell(y)}\mathrm{d}x}{\int_D  \varphi(x)^{2+\alpha}\frac{\bar\ell(y\varphi(x))}{\bar\ell(y)}\mathrm{d}x} =  1.
	\]
	Therefore, \ref{H4} holds. 
	    Further, from \eqref{eq:pgfN} we deduce that $\mathcal{E}[N^{1+\delta}]<\infty$ for any $\delta<\alpha$, which implies \ref{H5}. So we can apply our results in this setting to deduce Kolmogorov's estimation and Yaglom's limit for the case of infinite variance completing thus completing the work in \cite{Ellen}.

\bigskip

\noindent{\bf Continuous-state branching processes:}
As a very special example of the superprocess setting, we may assume that particles have no spatial motion (which is akin to $\sP$ corresponding to a particle that is stationary at the origin), and $\phi\equiv0$, that is, there is no non-local branching. This is the setting of a continuous-state branching process (CSBP). Moreover, the local branching mechanism does not depend on $x\in E$ and criticality of \ref{H2} enforces the requirement that $\psi'(0) = 0$, hence 
\[
\psi (\lambda) = c\lambda^2 + \int_{0}^\infty( {\rm e}^{-\lambda y}-1+\lambda y)\nu(\dd y),\;\;\;\lambda\geq 0.
\]
In this setting \ref{H2} is also satisfied with $\varphi = 1$ and $\tilde{\varphi}(\dd x) = \delta_0(\dd x)$. In addition, \ref{H1} is also automatically satisfied. We assume the following integrability condition {(\textit{Grey's condition})}
\[\int_{\theta}^{\infty} \frac{1}{\psi(\lambda)}\mathrm{d}\lambda<\infty\]
for some $\theta>0$, which is a necessary and sufficient condition to have extinction in finite time almost surely (see Grey \cite{grey1974}), in other words, to guarantee \ref{H3}. 
{Further, according to \cite[Theorem 3.14 and Corollary 3.15]{ZL11} and equation \eqref{sufcondH6} in {our} Remark \ref{remH6}, Grey's condition implies \ref{H6}.}
The assumption \ref{H4} can otherwise be worded as requiring that 
\begin{equation*}\label{eq:stablemecha}
\psi(\lambda)= \lambda^{1+\alpha}\ell(\lambda),\text{ as }\lambda\to0.
\end{equation*}
In particular, by Karamata's Tauberian theorem, the above implies that the map $y \mapsto \nu(y,\infty)$ is regularly varying at $\infty$ with index $-\alpha$. Therefore, \ref{H5}, which translates to 
\[\int_{{(1,\infty)}} y^{1+\delta} \nu(\dd y)<\infty,\]
holds for $\delta<\alpha$.

As such our results tell us that, 
there exists a slowly varying function $\tilde\ell$ at $\infty$ such that, 
\[
 \lim_{t\to \infty} \frac{\mathbb{P}_z (\zeta>t)}{t^{-1/\alpha}\tilde{\ell}(t)} = z, 
\]
and
\[
    \lim_{t\to \infty}\mathbb{E}_z\left[\exp\left(-\theta t^{-1/\alpha} \tilde{\ell} (t) Z_t \right) \Big|\Big. \ \zeta>t \right] = 1- \frac{ \theta}{(1+ \theta^{\alpha})^{1/\alpha}},\qquad \theta\geq0.
\]
This conclusion previously appeared in Proposition 1.5 of \cite{labbe2013} and Theorem 3.2 of \cite{RenYangZhao2014}. The case $\alpha=0$ was also carried out in the former reference.

\bigskip

\noindent{\textbf{Multi-type continuous-state branching processes:}} 
    These processes are the multidimensional counterpart of the previous example and were originally introduced by \cite{Watanabe-MCSBP} for the two-dimensional case. They can be described as solutions to stochastic differential equations (see \cite{BLP}) or through their semigroup properties looking at them as superprocess (see \cite{KP}). We follow this second approach that essentially correspond to the setting $E = \{1,\ldots, n\}$, $n\in\mathbb{N}$. Recall that $B(E)$ is the space of bounded measurable functions on $E$, so it is isomorphic to $\mathbb{R}^n$. Likewise, the space $M(E)$ of finite Borel measures on $E$ is isomorphic to $ [0,\infty)^n$.  $\mathsf{P}$ is the identity operator and the local branching is the one given in \eqref{local-branching-mechanism}. For the non-local branching, we restrict ourselves to the setting of \cite{KP}. To be precise, given  $\beta , \tilde\gamma\in B^+(E)$, a probability distribution $\pi_x$ on $E\setminus \{x\}$ for all $x\in E$ (specifically $\pi_x(x)=0$) and a bounded kernel $u\tilde\Gamma(x,\dd u)$ from $E$ to $(0,\infty)$ satisfying 
    \begin{equation}\label{eq:boundGamma}
    \tilde\gamma(x) + \int_0^\infty  u \tilde\Gamma(x,\dd u) \leq 1,
    \end{equation}
    then the non-local branching mechanism is defined by \eqref{non-local-branching-mechanism} with the identification
    \[
    \gamma(x,f) = \tilde\gamma(x) \langle f, \pi_x\rangle - f(x),
    \qquad
    \Gamma(x,\dd \nu) = \int_0^\infty  \delta_{u\pi_x}(\dd \nu) \tilde\Gamma(x,\dd u).
    \]
    In this framework, \ref{H2} would be a natural ergodic assumption where a simple irreducibility condition ensures it. In other words, we recover the classical Perron--Frobenius behaviour for the matrix of the mean semigroup.
    Indeed, consider the matrix $L$ given by
    \[
    (L)_{ij} = b(i) \delta_{i=j} + \beta(i) \left(\tilde\gamma(i) + \int_0^\infty u \tilde\Gamma(i, \dd u) \right)\pi_i(j),
    \]
    where $b\in B(E)$. The matrix $L$ is the infinitesimal generator of the linear semigroup $\mathsf{T}$ (see, for example, Section 4 of \cite{KP}). Let $\varphi$ and $\tilde\varphi$ be respectively the right and left eigenvector of $L$ associated to the 0 leading eigenvalue (critical case).  In such a case, $\mathsf{T}_t[f](x) = \mathrm{e}^{tL}[f](x)\to \langle f, \tilde\varphi\rangle \varphi(x)$ as $t\to\infty$ by Perron--Frobenius theorem, and then \ref{H2} follows taking into account that  
    \[
    \Delta_t = \sup_{f\in B^+(E), \|f\|=1} \max_{i=1,\ldots,n} \left| \frac{\mathrm{e}^{tL} [f] (i)}{\varphi(i) \sum_{j=1}^n f(j) \tilde\varphi(j)} -1\right|.
    \]
    In addition, if
    \[
    \int_1^\infty (y\log y) \nu (i,\dd y) + \int_1^\infty (y\log y) \tilde\Gamma (i,\dd y) < \infty , \qquad \text{for all } i\in E, 
    \]
    then \ref{H3} holds at criticality due to Theorem 2 of \cite{KP} taking into account that in the setting of finite types local and global extinction coincide. 
    We now observe that
     \ref{H4} translates as 
    \begin{align*}
    \langle \mathsf{J} [x\varphi],\tilde\varphi\rangle 
    &= 
    \sum_{i=1}^n \tilde\varphi (i)\left[x^2 c(i)\varphi(i)^2 
    + 
    \int_0^\infty \left(\mathrm{e}^{-x y\varphi(i)}-1+xy\varphi(i)\right)\nu(i,\dd y) \right.\\
   &\quad\quad\quad\quad\quad\quad +
    \left.\beta(i)\int_0^\infty  \left(\mathrm{e}^{-x u\sum_{j=1}^n \varphi(j) \pi_i(j)}-1+xu \sum_{j=1}^n \varphi(j) \pi_i(j)\right)\tilde\Gamma(i,\dd u)\right] \\
    &= x^{1+\alpha} \ell(x), \qquad \text{as } x\to 0.
    \end{align*}
Finally, we show that \ref{H4} implies \ref{H5}, which in this case takes the form
\[
 \max_{i=1,\ldots,n} \left(\int_{1}^\infty y ^{1+\delta} \nu(i,\mathrm{d} y) + \frac{\beta(i)\sum_{j=1}^n \varphi(j) \pi_i(j)}{\varphi(i)} \int_{0}^\infty {(1\vee y^{\delta})} y   \tilde\Gamma (i, \mathrm{d} y)\right) < \infty,
\]
for some $\delta \in (0,\alpha)$. 
Firstly, since $\varphi,\tilde\varphi$ are eigenvectors with strictly positive coordinates, note that \[\max_{i=1,\ldots,n} \varphi(i)^{-1} \beta(i)\sum_{j=1}^n \varphi(j) \pi_i(j) < \infty.\] Further, $ \max_{i=1,\ldots,n} \int_{0}^1  y   \tilde\Gamma (i, \mathrm{d} y) < \infty $ by \eqref{eq:boundGamma}.
Let us see now that
\begin{equation}\label{exampleH5}
	\int_{1}^\infty y ^{1+\delta} \nu(i,\mathrm{d} y) < \infty, \qquad \text{for all } i=1,\ldots,n.
\end{equation}
Using Tonelli's theorem,
we have 
\[
\int_{1}^\infty y ^{1+\delta} \nu(i,\mathrm{d} y) = \nu (i, [1,\infty)) + (1+\delta) \int_{1}^\infty y ^{\delta} \nu (i, [y,\infty)) \mathrm{d} y.
\]
Taking into account that $\int_1^\infty y ^{\delta-1-\alpha} \bar{\ell}(y) \mathrm{d} y < \infty$ if $\delta<\alpha$, with $\bar\ell$ a slowly varying function at $\infty$, it is enough to prove that $ \nu (i, [y,\infty)) = O\left(y^{-1-\alpha} \bar{\ell}(y)\right)$ as $y\to\infty$. 
We use again Tonelli's theorem,
so that
\[
\int_{y}^\infty u \nu(i,\dd u) =  y \nu (i, [y,\infty)) + \int_y^\infty \nu (i, [u,\infty)) \dd u,
\]
and, from the deterministic inequalities $\mathrm{e}^{-z}-1+z \geq \frac12 (1 - \mathrm{e}^{-z/2})z \geq \frac12 (1 - \mathrm{e}^{-1/2})z$ if $z\geq 1$, we deduce
\begin{align*}
	\langle \mathsf{J} [x\varphi],\tilde\varphi\rangle 
	&\geq  
	\tilde\varphi (i)    \int_0^\infty \left(\mathrm{e}^{-x y\varphi(i)}-1+xy\varphi(i)\right)\nu(i,\dd y) 
	\\
	&\geq 
	\frac{1-\mathrm{e}^{-1/2}}{2} \tilde\varphi (i)    \varphi(i) x\int_{\frac{1}{x\varphi(i)}}^\infty y \nu(i,\dd y).
\end{align*}
Since $\varphi,\tilde\varphi$ are eigenvectors with strictly positive coordinates,  
\[\int_{\frac{1}{x\varphi(i)}}^\infty y \nu(i,\dd y) = O \left(x^{\alpha} \ell(x)\right), \qquad \text{as } x\to 0,\]
and 
\[ \nu (i, [y,\infty)) \leq  y^{-1} \int_{y}^\infty u \nu(i,\dd u) = O\left(y^{-1-\alpha} \bar{\ell}(y)\right) , \qquad \text{as } y\to \infty,\]
completing \eqref{exampleH5}. The proof of 
\[
\int_{1}^\infty y ^{1+\delta} \tilde\Gamma(i,\mathrm{d} y) < \infty, \qquad \text{for all } i=1,\ldots,n,
\]
follows the same steps. 
{Finally, the projection from Remark \ref{remH6} in this example takes the form
\begin{align*}
\Phi_1(i,\lambda)
&= 
\psi(i, \lambda) + \beta(i) \left(2-\tilde\gamma(i)-\int_0^\infty u \tilde\Gamma (i,\dd u) \right) \lambda
\\
&\geq
\lambda \min_{i=1,\ldots, n} \left( -b(i) + \beta(i)\right)  + \lambda^2 \min_{i=1,\ldots, n} c(i) + \int_{0}^\infty ( {\rm e}^{-\lambda y}-1+\lambda y)\nu_\ast(\dd y) =: \psi_\ast (\lambda),
\end{align*}
where the inequality follows from \eqref{eq:boundGamma} and from choosing the measure $\nu_\ast$ such that $\nu_\ast(\dd y) \leq \nu(i, \dd y)$ for all $i=1,\ldots,n$. We assume that $\lim_{\lambda \to \infty} \psi_\ast^\prime (\lambda) = \infty$. Therefore,  assumption \ref{H6} is satisfied.
}
	
	\bigskip
	

 \noindent \textbf{Local superprocesses with spatially dependent  stable branching mechanism:} We recover the results from \cite{ren2020limit} by setting the non-local branching mechanism $\phi \equiv 0$, and defining the local branching mechanism as \eqref{local-branching-mechanism} with
\begin{equation*}
	c(x)=0 \quad \text{and} \quad \nu(x, \mathrm{d}y) = \frac{\kappa(x)}{\Gamma(-\gamma(x))y^{1+\gamma(x)}}\mathrm{d}y, \quad x\in E,
\end{equation*}
where $\kappa ,\gamma \in B^+(E)$ such that $\kappa(\cdot)>0$ and $1<\gamma(\cdot)<2$, and $\Gamma$ is the Gamma function defined on the negative half line. Thus
\begin{equation*}
	\psi(x,\lambda)=-b(x) \lambda + \kappa(x)\lambda^{\gamma(x)},\qquad \lambda\geq 0.
\end{equation*}
We assume that $\sP=(\sP_t, t\geq0)$ is absolutely continuous with respect to a Borel finite measure $m$ on $E$ with densities $(p_t(\cdot,\cdot), t>0)$ on $E\times E$, i.e.
\[\sP_t [f](x)= \int_{E}p_t(x,y)f(y) m(\mathrm{d}y), \quad t>0, \  x\in E, \ f\in B^+(E).\]
Moreover, we assume that $m$ and $(p_t(\cdot,\cdot), t>0)$ satisfy Assumption 1 in  \cite{ren2020limit}. Assumptions 2 and 3 in \cite{ren2020limit} imply our hypothesis \ref{H2}.
In this particular setting $\tilde \varphi(\mathrm{d}x) = \varphi^*(x) m(\mathrm{d}x)$, where $\varphi^*$ is the principal eigenfunction of the generator of the dual semigroup of $\sT$. 
Assuming that $\gamma_0 := \sup\{r: m(x: \gamma(x)<r)=0\}>1$ and $\kappa_0:=\sup\{r: m(x: \kappa(x)<r)=0\}>0$, then
\begin{equation*}
\langle \textsf{J}[x\varphi], \tilde{\varphi} \rangle = \langle \kappa (x\varphi)^{\gamma}, \tilde \varphi  \rangle \sim   x^{\gamma_0}\langle \mathbf{1}_{(\gamma=\gamma_0)} \kappa \varphi^{\gamma_0}, \tilde{\varphi}\rangle, \quad \text{as} \ x \to 0,
\end{equation*}
(see Lemma 2.6 or the proof of Proposition 3.3 in \cite{ren2020limit}). Therefore, \ref{H4} is satisfied with $1+\alpha =\gamma_0$ and $\ell$ a function which converges to the constant $\langle \mathbf{1}_{(\gamma=\gamma_0)} \kappa \varphi^{\gamma_0}, \tilde{\varphi}\rangle$ as $x \to 0$. On the other hand, \ref{H3} follows from the specific form of the branching mechanism, together with the assumptions that $\gamma_0>1$ and $\kappa_0>0$ (see Proposition 3.1 in \cite{ren2020limit}). Further, \ref{H5} can be written as
\[\sup_{x\in E} \frac{\kappa(x)}{\Gamma(-\gamma(x)) } \frac{1}{\gamma(x)-1-\delta} <\infty,\]
which holds for some $\delta\in (0,\inf_{x\in E} \gamma(x)-1)$ if $\inf_{x\in E} \gamma(x)>1$. {Finally, since we are assuming $c(x)=0$ for every $x\in E$, then \ref{H6} automatically holds.}


\section{Preliminaries}
In this section, we consider several semigroup evolution equations that will be useful for proving our main results. 

\subsection{Non-local branching Markov processes}
In the setting of the $(\sP, \sG)$-branching Markov process and under \ref{H1}, the evolution equation for the expectation semigroup $(\sT_t, t \ge 0)$ is given by
\begin{equation}\label{eq:BMP-linear}
\sT_t[f](x) = \sP_t[f](x) + \int_0^t \sP_s[\beta (\sm[\sT_{t-s}[f]] - \sT_{t-s}[f])](x) \dd s,
\end{equation}
for $t \ge 0$, $x \in E$ and $f \in B^+(E)$, where
\begin{equation}\label{eq:mparticle}
\sm[f](x)= \mathcal{E}_x[\langle f, \sZ\rangle] .
\end{equation}
See, for example, Lemma 8.1 of \cite{Horton2023}. 

\medskip

Our next evolution equation will relate the non-linear semigroup  to the linear semigroup, which will enable us to use {\ref{H2}}, for example, to study the limiting behaviour of $\sv_t$. For this, we will introduce the following modification to the non-linear semigroup,
\begin{equation}
  \su_t[g](x) := \mathbb E_{\delta_x}\left[1-\prod_{i = 1}^{N_t}g(x_i(t)) \right] = 1-\mathrm{e}^{-\sv_t[-\ln g](x)}, \qquad g\in B^+_1(E).
  \label{1-u has smg property}
\end{equation}
Recalling that we have assumed first moments of the offspring distribution, one can show that
\begin{equation}\label{eq:nonlin-lin}
\su_t[g](x) = \sT_t[1-g](x) - \int_0^t \sT_s[\sA[\su_{t-s}[g]]](x) \dd s, \quad t \ge 0,
\end{equation} 
where, for $g \in B^+_1(E)$ and $x \in E$,  we recall that 
\begin{equation}
\sA[g](x) = \beta(x)\mathcal E_x\left[\prod_{i = 1}^N (1-g(x_i)) - 1 + \sum_{i = 1}^N g(x_i)\right]\geq0.
\label{A}
\end{equation}
We refer the reader to Theorem 8.2 of \cite{Horton2023} for a more general version of \eqref{eq:nonlin-lin}, along with a proof. Further, recall that the expectation semigroup $(\sT_t, t \ge 0)$ is defined for $f\in B^+(E)$. Note that $\sA[g] \in B^+(E)$ for $g\in B^+_1(E)$, since 
\begin{equation*}
	\|\sA[g]\| \le \|\beta\|  \sup_{x\in E}  \mathcal E_x \left[\sum_{i = 1}^N g(x_i)\right] \le  \|\beta\|  \sup_{x\in E}  \mathcal E_x [N],
\end{equation*}
where the right-hand side is finite by assumption \ref{H1} and the fact that  $\beta \in B^+(E)$. {We also have the following useful monotonicity property.}
\begin{lemma}\label{lem:Amono}
		Under \ref{H1},
	for $f,g\in B_1^+(E)$,
	\begin{equation}\label{eq:Amonotone}
		\mathsf{A}[f] \leq \mathsf{A}[g] \qquad \text{if} \quad f\leq g.
	\end{equation}
\end{lemma}
\begin{proof}
	It is enough to see that
	\[
	F(z_1,\ldots,z_N) = \prod_{i=1}^N (1-z_i) -1 + \sum_{i=1}^N z_i, \qquad N\in \mathbb{N}, \quad z_1,\ldots,z_N\in [0,1],
	\]
	is a non-decreasing function in each coordinate. Indeed,
	\[
	\frac{\partial F}{\partial z_j} = 1- \prod_{i=1, i\neq j}^N (1-z_i) \geq 0, \qquad \text{for all }\ j=1,\ldots, N.
	\]
	This completes the proof of the lemma.
\end{proof}

An important part of our analysis will be to consider the behaviour of $\su_t[g](x)$, $t\geq0$, $g\in B^+_1(E)$ and compare it to 
\begin{equation}
a_t : = \langle \su_t[0],\tilde{\varphi}\rangle, \qquad t\geq0,
\label{atdef}
\end{equation}
where we note that 
\begin{equation}\label{def-u-prob}
	\su_t(x): = \su_t[0](x) = \mathbb{P}_{\delta_x}(\zeta>t), \qquad t\geq0.
\end{equation}
Note that \ref{H3} implies that $\lim_{t\to\infty}\su_t(x) = 0$ and dominated convergence  in \eqref{atdef} similarly implies $\lim_{t\to\infty}a_t = 0$. 

The reason for the comparison is that $(a_t,t\geq0)$ solves a simpler evolution equation. From \eqref{eq:nonlin-lin} and \eqref{criticaldef} it is easy to show that 
\begin{equation}\label{eq:at}
a_t=\langle 1, \tilde{\varphi} \rangle -\int_{0}^{t}\langle  \textsf{A}[\su_{t-s}], \tilde{\varphi} \rangle \mathrm{d} s, \qquad t\geq0.
\end{equation}

As an initial indication of how $\su_t[g](x)$ compares to $a_t$ as $t\to\infty$, note from \eqref{eq:nonlin-lin}  that 
\[
\su_{t+s}[g](x) = \sT_t[\su_s[g]] (x)- \int_0^t \sT_\ell[\sA[\su_{t+s-\ell}[g]]](x) \dd \ell \leq \sT_t[\su_s {[g]}](x),
\]
where we use the fact that $\su_{t+s}[g]=\su_t[1-\su_{s}[g]]$ which follows by the semigroup property of $\sv_t[\cdot]$. Therefore,
under the assumptions \ref{H1}-\ref{H3}, since $\su_t[g]\in B_1^+(E)$,
\begin{equation}\label{E:ut-a-bd}
	 \|\varphi^{-1}\su_{t+s}\| \le \| \varphi^{-1}\sT_t[\su_s] \| 
	  \leq {a_s (1+\Delta_t)}.
	 \end{equation}
 As $t\to\infty$ and then $s\to\infty$, we see that 
 \begin{equation}
\lim_{t\to\infty} \|\varphi^{-1}\su_{t}\| = 0 
\qquad 
{\text{and}
\qquad
\lim_{t\to\infty} \|\su_{t}\| = 0.}
\label{uniformutto0}
 \end{equation}

\subsection{Non-local superprocesses}
In the setting of the  $(\sP,\psi,\phi)$-superprocess,
the evolution equation for the expectation semigroup $(\sT_t, t\geq0)$ is well known and satisfies 
\begin{equation}
\sT_t\bra{f}(x)=\sP_t[f](x)
+\int_{0}^{t}\sP_s\bra{\beta (\sm[\sT_{t-s}[f]]-\sT_{t-s}[f])+b\sT_{t-s}[f]} (x)\dd s,
\label{superm21}
\end{equation} 
for $t\geq 0$, $x\in E$ and $f\in B^+(E)$, where 
\begin{equation}\label{supermh}
\sm[f](x) = \gamma(x,f)+\int_{M_{0}(E)}\langle{f},{\nu}\rangle\Gamma(x,\dd \nu).
\end{equation}
See for example  equation (3.24) of \cite{dawson2002nonlocal}.

\medskip

Similarly to the branching Markov process setting, let us re-write {\color{black} an extended version of the non-linear semigroup evolution $({\sV}_t, t\geq0)$, defined in \eqref{non-local-evolution-equation}, i.e. the natural analogue of \eqref{nonlinv},} in terms of the linear semigroup $(\sT_t, t\geq0)$. 
Indeed, from Lemma 4 of \cite{GHK}, we have the following evolution equation, 
\begin{equation}
{\sV}_t[f](x)= \sT_{t}[f](x) - \int_{0}^{t}\sT_{s}\left[\sJ[{\sV}_{t-s}[f]]\right](x)\dd s, \quad f \in B^+(E), x \in E, t\geq 0,
\label{nonlinvJ}
\end{equation}
where, for $h\in B^+(E)$ and $x\in E$, recall that
\begin{align}
\sJ[h](x) 
&=c(x)h(x)^2 +\int_{(0,\infty)}( {\rm e}^{-h(x) y}-1+h(x) y)\nu(x,\dd y)\notag\\
&\hspace{1cm} +\beta(x)\int_{M_{0}(E)}({\rm e}^{-\langle{h},{\nu}\rangle}-1+ \langle h, \nu\rangle)\Gamma(x,\dd \nu).
\label{J} \notag
\end{align}

The analogue of $\su_t(x)$ for the superprocess setting is given by $\sV_t(x) : = \lim_{\theta\to\infty}\sV_t[\theta](x)$. As such 
\begin{equation}\label{surv-prob-sup}
1-{\rm e}^{-\sV_t(x)} = \mathbb{P}_{\delta_x}(\zeta>t), \qquad t\geq0,
\end{equation}
and taking account of \ref{H3}, we also note that 
$\lim_{t\to \infty}\sV_t(x) = 0$ and hence 
\begin{equation}
\lim_{t\to\infty}\frac{\mathbb{P}_{\delta_x}(\zeta>t)}{\sV_t(x)} = 1.
\label{P/V}
\end{equation}

Similarly $(a_t, t\geq0)$ is replaced by 
\begin{equation}\label{atsup}
a_t = \langle \sV_t, \tilde\varphi\rangle, \qquad t\geq0
\end{equation}
and dominated convergence again implies that $\lim_{t\to\infty}a_t = 0$.

Although $V_0= \infty$, we can describe the evolution of $(a_t, t\geq0)$, avoiding this singularity, by choosing $t_0>0$ and writing 
\begin{equation}
a_t = a_{t_0} - \int_{t_0}^t\langle \sJ[{\sV}_{t-s}],\tilde\varphi\rangle\dd s, \qquad t\geq t_0.
\label{atJ}
\end{equation}

Similarly to \eqref{uniformutto0}, under the assumptions \ref{H2} and \ref{H3},  we can use the limiting behaviour of $(a_t, t\geq0)$ to show  
$\lim_{t\to\infty}\|\varphi^{-1}\sV_{t}\|=0.$

{Finally, we note that Lemma \ref{lem:Amono} also admits a counterpart in the superprocess setting.}
\begin{lemma}\label{lem:Jmonot}
	For $f,g\in B^+(E)$,
	\begin{equation}\label{eq:Jmonotone}
		\mathsf{J}[f] \leq \mathsf{J}[g] \qquad  \text{if}  \quad f\leq g.
	\end{equation}
\end{lemma}
\begin{proof}
	The proof easily follows from the fact that the function $x\in[0,\infty) \longmapsto \mathrm{e}^{-x}-1+x$ is increasing.
\end{proof}

\section{Proof of Theorem \ref{theo:kolmogorov}}

We give the proof in detail for the setting of branching Markov processes. We indicate only the minor changes needed for the setting of superprocesses.

\subsection{Non-local branching Markov processes}
The proof is made up of several smaller results. Our first result shows that, as alluded to above, we may uniformly replace $\varphi^{-1}(x)\su_t(x) $ by $a_t$ as $t\to\infty$. Our proof is largely guided by a similar reasoning in Proposition 3.2 of \cite{ren2020limit}.

\begin{lemma}\label{lem:asymu}
Under the assumptions {\ref{H1}-\ref{H3} and \ref{H5}}, it holds that 
\begin{equation}
\sup_{x\in E} \left| \frac{\emph{\su}_t(x)}{a_t \varphi(x)} - 1 \right| \to 0, \quad \text{as}\quad t\to \infty.
\label{ut/at}
\end{equation}
\end{lemma}

\begin{proof}
First we will prove that there is a sequence $(A_t, t \ge 0)$ such that 
\begin{equation}\label{goalAt}
\sup_{x\in E} \left| \frac{\varphi(x)^{-1}\su_t(x)}{A_t} - 1 \right| \to 0, \quad \text{as}\quad t\to \infty.
\end{equation}
In that case it follows by dominated convergence that 
\begin{equation*}
\left| \frac{a_t}{A_t} - 1 \right| \leq  \int_E \left| \frac{\varphi(y)^{-1}\su_t(y)}{A_t}-1  \right|\varphi\tilde\varphi({\rm d} y) \to 0, \quad \text{as}\quad t\to \infty,
\end{equation*}
and hence \eqref{ut/at} follows.

To find the normalisation $A_t$, consider the spine change of measure, 
\begin{equation}
\left.\frac{{\rm d}\mathbb{P}^\varphi_{\delta_x}}{{\rm d}\mathbb{P}_{\delta_x}}\right|_{\mathcal{F}_t}  =  \frac{\langle \varphi, X_t\rangle}{\varphi(x)} : = \varphi(x)^{-1}Y_t, \qquad t\geq0.
\label{COM}
\end{equation}
Note  that 
\begin{equation}
\mathbb{E}^\varphi_{\delta_x}[1/Y_t]  = \varphi(x)^{-1}\mathbb{P}_{\delta_x}(\zeta>t) = \varphi(x)^{-1}\su_t(x)\text{ and }\mathbb{E}^\varphi_{\varphi\tilde\varphi}[1/Y_t]  = a_t.
\label{notethat}
\end{equation}

We recall from Chapter 11 of \cite{Horton2023} that under $\mathbb{P}^\varphi_{\delta_x}$, the process is equivalent in law to the aggregated mass from following spine decomposition:
\begin{itemize}
\item  The initial ancestor is selected and marked `{\it spine}'.

 \item 
 For the marked particle, issue a copy of the process whose motion is determined by the semigroup
\begin{equation}
{\sP}_t^\varphi[f](x) := 
\frac{\sT_t[\varphi f](x)}{\varphi(x)},
 \label{motion-bias}
\end{equation}
for $x\in E$ and $f\in B^+(E)$.


\item   When at $x \in E$, the marked particle undergoes branching at rate
\begin{equation}
\beta^\varphi(x): =
\beta(x)\frac{\sm [\varphi](x)}{\varphi(x)},
\label{gammaphi}
\end{equation}
 at which point, it scatters a random number of particles in type space $E$ according to the random measure given by $(\sZ, {\mathcal P}^\varphi_{x})$ where 
\begin{equation}
\frac{\dd{\mathcal P}^\varphi_{x}}{\dd {\mathcal P}_{x}} = \frac{\langle\varphi, \sZ\rangle}{\sm [\varphi](x)}.
\label{pointprocessCOMgen}
\end{equation}

\item  When the marked particle  scatters  its offspring from an instantaneous type $x\in E$, given $Z = (x_1,\dots, x_N)$, select the continuation of the `spine' marked genealogy using the empirical probabilities ${\varphi(x_i)}/{{\langle \varphi, \sZ\rangle}}$, $i = 1,\dots, N$. The particle marked as a spine repeats the stochastic behaviour from the second bullet point above. 
\item
 Each  unmarked particle $j\neq i$  from the offspring point process $\sZ = \sum_{i=1}^N \delta_{x_i}$ issues an independent copy of $(X, \mathbb{P}_{\delta_{x_j}})$.
\end{itemize}

A consequence of the change of measure \eqref{COM} and its stochastic representation is that, if we denote by $(\xi_t,t\geq0)$ by the position of the particle marked as the spine, then
\begin{equation}
\mathbb{E}_{\delta_x}^\varphi[f(\xi_t)] = \mathbb{E}^\varphi_{\delta_x}\left[\frac{\langle \varphi f, X_t\rangle}{\langle \varphi, X_t\rangle}\right] =\frac{\sT_t[\varphi f](x)}{\varphi(x)}, \qquad f\in B^+(E),\ t\geq0,\ x\in E.
\label{spineunderspine}
\end{equation}
See for example (11.14) and Theorem 11.1 of \cite{Horton2023}.

Next define  for any $t>0$ and $t_0 \in (0,t)$,
\begin{equation}
Y_t^{(t_0,t]} =  \varphi(\xi_t)+\sum_{i = n_{t_0}+1}^{n_t}\sum_{j = 1, j\neq i^*}^{N^i} Y^j_{t-T_i} ,
\label{Y}
\end{equation}
with $T_i$, $i = 1,\dots, n_t$, the times and number of birth events along the spine up to time $t$ and $i^*$ is the index of the spine among siblings. Further, $N^i$ is the number of particles created by individual $i$ for $i=n_{t_0+1}, \dots, n_t$ and  $Y_{t-T_i}^j$ is a copy of $Y_{t-T_i}$ for $j=1,\dots, N^{i}$ with $j\not=i^*$ (see Fig. \ref{fig:spine}).  Note that $Y_t = Y_t^{(0,t]}$ thanks to the spine decomposition. 

\begin{figure}
	\includegraphics[width=7cm]{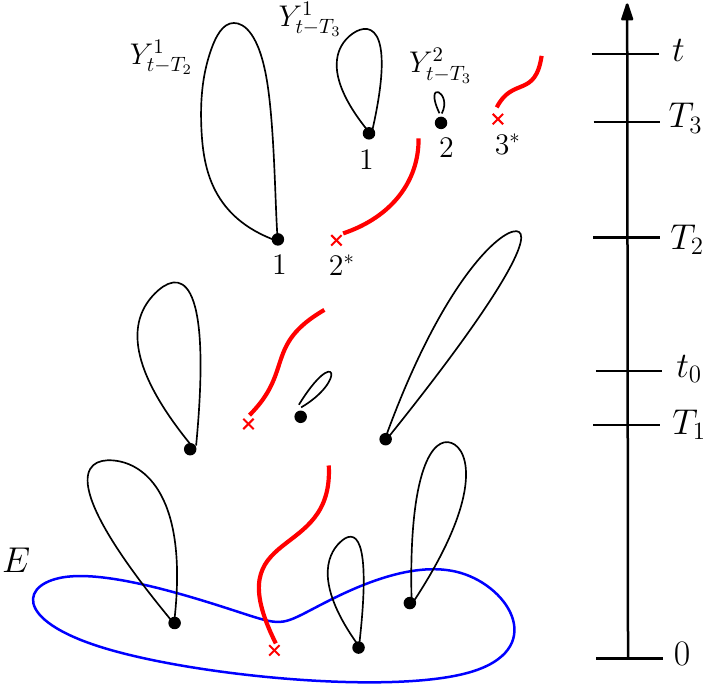}
	\caption{A realization of the components in \eqref{Y}. In this realization $n_{t_0} = 1$, $n_t = 3$.}
	\label{fig:spine}
\end{figure}

Now write 
\[
\epsilon^2_x(t_0, t) = \left|\mathbb{E}^\varphi_{\delta_x}[1/Y_t] - \mathbb{E}^\varphi_{\delta_x}[1/Y_t^{(t_0,t]}]\right| \ \text{ and } \ \epsilon^1_x(t_0, t) = \left|\mathbb{E}^\varphi_{\delta_x}[1/Y_t^{(t_0,t]}] - \mathbb{E}^\varphi_{\varphi\tilde\varphi}[1/Y_t^{(t_0,t]}]\right|.
\]
We claim that, for $t_0=t_0(t)\in(0,t)$,
we have 
\begin{equation}
\left|  \frac{\varphi(x)^{-1}\su_t(x)}{a_{t-t_0} } - 1 \right|=\frac{\left|  \varphi(x)^{-1}\su_t(x) - a_{t-t_0} \right|}{a_{t-t_0} }\leq \frac{\epsilon^1_x(t_0, t)}{a_{t-t_0}} +\frac{\epsilon^2_x(t_0, t)}{a_{t-t_0}}. 
\label{claim}
\end{equation}
To see why, note that if we define $\mathcal{F}^\xi_{s}$  to be the filtration generated by the motion of the spine, denoted here by  $\xi$, up to time $s$. Then, thanks to the Markovian nature of the spine, its immigration process and \eqref{notethat},  
\begin{align}
\mathbb{E}^\varphi_{\delta_x}[1/Y_t^{(t_0,t]}] &= \mathbb{E}^\varphi_{\delta_x}[\mathbb{E}^\varphi_{\delta_x}[1/Y_t^{(t_0,t]}|\mathcal{F}^\xi_{t_0}]] \notag\\
&=  \mathbb{E}^\varphi_{\delta_x}[{\mathbb{E}^\varphi_{\delta_{\xi_{t_0}}}[1/Y_{t-t_0}]}]\notag\\
&= \mathbb{E}^\varphi_{\delta_x}[\varphi(\xi_{t_0})^{-1}\su_{t-t_0}(\xi_{t_0})]
\label{useforepsilon2}
\end{align}
and similarly, by \eqref{spineunderspine}, 
\begin{equation}
\mathbb{E}^\varphi_{\varphi\tilde\varphi}[1/Y_t^{(t_0,t]}] = \mathbb{E}^\varphi_{\varphi\tilde\varphi}[\varphi(\xi_{t_0})^{-1}\su_{t-t_0}(\xi_{t_0})] = \langle \su_{t-t_0},\tilde\varphi\rangle = a_{t-t_0},
\label{othernotethat}
\end{equation}
where we have used the known fact that  $\varphi\tilde\varphi$ is invariant for $\xi$; cf. Chapter 11 of \cite{Horton2023}.
The claim \eqref{claim} now follows from \eqref{notethat} and \eqref{othernotethat} and the triangle inequality.

Next note that  
\begin{align}
\epsilon^1_x(t_0, t)& =\left| \mathbb{E}^\varphi_{\varphi\tilde\varphi}[\varphi(\xi_{t_0})^{-1}\su_{t-t_0}(\xi_{t_0})] -  \mathbb{E}^\varphi_{\delta_x}[\varphi(\xi_{t_0})^{-1}\su_{t-t_0}(\xi_{t_0})]\right|
\notag\\
&= \left| \mathbb{E}^\varphi_{\varphi\tilde\varphi}[f_{t-t_0}(\xi_{t_0})] -  \mathbb{E}^\varphi_{\delta_x}[f_{t-t_0}(\xi_{t_0})]\right| \notag\\
&= \left| \langle f_{t-t_0}, \varphi\tilde\varphi\rangle -  \mathbb{E}^\varphi_{\delta_x}[f_{t-t_0}(\xi_{t_0})]\right| \notag\\
&=\left|\varphi^{-1}(x)\sT_{t_0}[\varphi f_{t-t_0}](x) - \langle\varphi f_{t-t_0}, \tilde\varphi\rangle
\right|\notag\\
&\leq  \Delta_{t_0}  \langle\varphi f_{t-t_0}, \tilde\varphi\rangle\notag\\
&= \Delta_{t_0} a_{t-t_0}, 
\label{e1estimate}
\end{align}
where $f_t (x) = \varphi(x)^{-1}\su_t(x)$, we have used \eqref{spineunderspine} in the fourth equality and for the inequality,  
{we have used  \ref{H2} to deduce that
\[
\left|\frac{\sT_s\bra{g}(x)}{\langle g,\tilde\varphi\rangle\varphi(x)}-1\right| \leq \Delta_{s},
\qquad s\geq 0,\  x\in E, \ g\in B_1^+(E)\setminus \{0\}.
\]
} 
In conclusion, uniformly in $x\in E$,
\begin{equation}
0\leq \frac{\epsilon^1_x(t_0, t)}{a_{t-t_0}} \leq \Delta_{t_0}.
\label{e1deltabound}
\end{equation}

Next define 
\[
Y_t^{(0,t_0]} = Y_t - Y_t^{(t_0, t]} = \sum_{i = 1}^{n_{t_0}}\sum_{j = 1, j\neq i^*}^{N^i} Y^j_{t-T_i}
\] 
and note that 
\begin{align}
\epsilon^2_x(t_0, t)&= 
\left|
\mathbb{E}^\varphi_{\delta_x}[1/Y_t] - \mathbb{E}^\varphi_{\delta_x}[1/Y_t^{(t_0,t]}]\right|
\notag\\
&\le \mathbb{E}^\varphi_{\delta_x}\left[\frac{|Y_t^{(t_0,t]} - Y_t|}{Y_t Y_t^{(t_0,t]}}\right]\notag\\
&=\mathbb{E}^\varphi_{\delta_x}\left[\frac{Y_t^{(0,t_0]}}{Y_t Y_t^{(t_0,t]}}\right]\notag\\
&\leq \mathbb{E}^\varphi_{\delta_x}\left[\mathbf{1}_{(Y_t^{(0,t_0]}>0)}/Y_t^{(t_0,t]}\right]\notag\\
&=\mathbb{E}^\varphi_{\delta_x}\left[
\mathbb{E}^\varphi_{\delta_x} \left[\mathbf{1}_{(Y_t^{(0,t_0]}>0)} | \mathcal{F}^\xi_{t_0}\right] \mathbb{E}^\varphi_{\delta_x}\left[1/Y_t^{(t_0,t]} | \mathcal{F}^\xi_{t_0} \right]
\right].
\label{bigmess}
\end{align}

We now need to deal with the  conditional expectation $\mathbb{E}^\varphi_{\delta_x} [\mathbf{1}_{(Y_t^{(0,t_0]}>0)} | \mathcal{F}^\xi_{t_0}]$. To this end, we note
\begin{align}
\mathbb{E}^\varphi_{\delta_x} [\mathbf{1}_{(Y_t^{(0,t_0]}>0)} | \mathcal{F}^\xi_{t_0}]
&\leq -\log \left(1-\mathbb{E}^\varphi_{\delta_x} [\mathbf{1}_{(Y_t^{(0,t_0]}>0)} | \mathcal{F}^\xi_{t_0}] \right)\notag\\
&=\lim_{\lambda\to\infty} - \log  \mathbb{E}^\varphi_{\delta_x}[{\rm e}^{-\lambda Y^{(0,t_0]}_t} | \mathcal{F}^\xi_{t_0}].
\label{oddLaplace}
\end{align} 
Next, using the spine decomposition for the general set up, recalling the definitions of ${\mathcal P}^\varphi_{x}$ and $\sm$ given in \eqref{eq:mparticle} and \eqref{pointprocessCOMgen} and appealing to the Campbell’s formula, we have
\begin{align}
&\mathbb{E}^\varphi_{\delta_x}[{\rm e}^{-\lambda Y^{(0,t_0]}_t} | \mathcal{F}^\xi_{t_0}] \notag\\
&= 
\exp\left(
- \int_0^{t_0} 
\left(
1- {\mathcal E}^\varphi_{\xi_s}\left[\sum_{i = 1}^N\frac{\varphi(x_i)}{\langle \varphi, \sZ\rangle} \prod_{j \neq i} {\rm e}^{- \sv_{t-s}[\lambda \varphi](x_j)}\right]
\right)
\beta(\xi_s)\frac{\sm[\varphi](\xi_s)}{\varphi(\xi_s)}
{\rm d}s
\right)\notag\\
&=\exp\left(
 -\int_0^{t_0} 
\left(
\sm[\varphi](\xi_s)- {\mathcal E}_{\xi_s}\left[\sum_{i = 1}^N \varphi(x_i) \prod_{j \neq i} {\rm e}^{- \sv_{t-s}[\lambda \varphi](x_j)}\right]
\right)
\frac{\beta(\xi_s)}{\varphi(\xi_s)}
{\rm d}s
\right)\notag\\
&=\exp\left(
 - \int_0^{t_0} 
{\mathcal E}_{\xi_s}\left[\sum_{i = 1}^N \varphi(x_i) \left(1- \prod_{j \neq i} \Big(1- \su_{t-s}[\mathrm{e}^{-\lambda \varphi}](x_j)\Big)\right)\right]
\frac{\beta(\xi_s)}{\varphi(\xi_s)}
{\rm d}s\right)
\label{campbell}
\end{align}
where in the last equality we used \eqref{MBP} and \eqref{1-u has smg property} to get
\[
\mathbb{E}_{\delta_x}[{\rm e}^{-\lambda Y_t}] = \mathbb{E}_{\delta_x}[{\rm e}^{-\lambda \langle\varphi,X_t\rangle}] ={\rm e}^{- \sv_{t}[\lambda\varphi](x)} = 1- \su_t[{\rm e}^{-\lambda \varphi}](x).
\]
We now have by monotone convergence
\begin{align}
&\lim_{\lambda\to\infty} - \log  \mathbb{E}^\varphi_{\delta_x}[{\rm e}^{-\lambda Y^{(0,t_0]}_t} | \mathcal{F}^\xi_{t_0}]\notag\\
&=\lim_{\lambda\to\infty}\int_0^{t_0} 
{\mathcal E}_{\xi_s}\left[\sum_{i = 1}^N \varphi(x_i) \left(1- \prod_{j \neq i} \Big(1- \su_{t-s}[{\rm e}^{-\lambda \varphi}](x_j)\Big)\right)\right]
\frac{\beta(\xi_s)}{\varphi(\xi_s)}
{\rm d}s\notag\\
&=\int_0^{t_0} 
{\mathcal E}_{\xi_s}\left[\sum_{i = 1}^N \varphi(x_i) \left(1- \prod_{j \neq i} \Big(1- \su_{t-s}(x_j)\Big)\right)\right]
\frac{\beta(\xi_s)}{\varphi(\xi_s)}
{\rm d}s.
\label{limlam}
\end{align}
We want this to be bounded by some deterministic $C_t$, which is independent of $x$, such that $\lim_{t\to\infty} C_t = 0$. To this end, we need to appeal to \ref{H5}, use that $\su_{t}$ is non-increasing (recall definition \eqref{def-u-prob}) and the deterministic inequality {$1-(1 -z)^{n-1} \le ((n-1)z)^\delta$} for $z \in [0,1]$, $n\geq 1$, and $\delta\in (0,1)$, so that
\begin{align}
		\mathbb{E}^\varphi_{\delta_x} [\mathbf{1}_{(Y_t^{(0,t_0]}>0)} \mid  \mathcal{F}^\xi_{t_0}] 
		&\leq 
		\int_0^{t_0} 
		{\mathcal E}_{\xi_s}\left[\sum_{i = 1}^N \varphi(x_i) \left(1- \prod_{j \neq i} \Big(1- {\su_{t-t_0}(x_j)}\Big)\right)\right]
		\frac{\beta(\xi_s)}{\varphi(\xi_s)} {\rm d}s
		\notag \\
		&\leq 
		\int_0^{t_0} 
		{\mathcal E}_{\xi_s}\left[\sum_{i = 1}^N \varphi(x_i) \left(1- {\Big(1- \|\su_{t-t_0}\|\Big)^{N-1}}\right)\right]
		\frac{\beta(\xi_s)}{\varphi(\xi_s)} {\rm d}s
		\notag \\
		&\leq
		\int_0^{t_0} {\mathcal E}_{\xi_s}\left[\sum_{i = 1}^N \varphi(x_i) {\Big( (N-1)\| \su_{t-t_0}\|\Big)^{\delta}} \right]
		\frac{\beta(\xi_s)}{\varphi(\xi_s)}
		{\rm d}s 
		\notag \\
		&\leq \|\su_{t-t_0}\|^{\delta} t_0 C,
		\label{condind}
\end{align}
where
\begin{equation}
	C= \sup_{x\in E} \frac{\beta(x)}{\varphi(x)}  {\mathcal E}_{x}\left[(N-1)^{\delta}\sum_{i = 1}^N \varphi(x_i)\right]<\infty.
\end{equation}

Back in \eqref{bigmess}, using \eqref{useforepsilon2}, \eqref{condind}, \eqref{spineunderspine}  and \ref{H2} {(see also \eqref{E:ut-a-bd})},

\begin{align}
\epsilon^2_x(t_0, t)
&\leq 
C t_0\|\su_{t-t_0}\|^{\delta}  \mathbb{E}^\varphi_{\delta_x}[\varphi(\xi_{t_0})^{-1}\su_{t-t_0}(\xi_{t_0})]
\notag \\
& = C t_0\|\su_{t-t_0}\|^{\delta}  \frac{\sT_{t_0}[\su_{t-t_0}](x)}{\varphi(x)}\notag\\
&\leq C t_0\|\su_{t-t_0}\|^{\delta} (1+\Delta_{t_0}) a_{t-t_0}.
\label{separatingtheexpectations}
\end{align}
It follows that, uniformly in $x\in E$,
\begin{equation} \label{e2upperbound}
0\leq \frac{\epsilon^2_x(t_0, t)}{a_{t-t_0}} \leq  C t_0\|\su_{t-t_0}\|^{\delta} (1+\Delta_{t_0}).
\end{equation}
Note in particular that the upper estimates for $\epsilon^1_x(t_0, t)$ and $\epsilon^2_x(t_0, t)$ in \eqref{e1deltabound} and \eqref{e2upperbound} are independent of $x$ and hence, back in \eqref{claim}, 
\begin{equation}\label{plugging}
\sup_{x\in E}\left|  \frac{\varphi(x)^{-1}\su_t(x)}{a_{t-t_0} } - 1 \right|
\leq 
\sup_{x\in E}\frac{\epsilon^1_x(t_0, t)}{a_{t-t_0}} + \sup_{x\in E}\frac{\epsilon^2_x(t_0, t)}{a_{t-t_0}} 
\leq
\Delta_{t_0} + C t_0\|\su_{t-t_0}\|^{\delta} (1+\Delta_{t_0}).
\end{equation}

According to \eqref{uniformutto0}, there exists a map $t \mapsto t_0(t)$ such that,
\[
t_0(t) \to \infty 
\qquad \text{and} \qquad
t_0(t)\|\su_{t-t_0(t)}\|^{\delta} \to 0,
\]
as $t\to\infty$.
For instance, take $ t_0 =  \|\su_{t-\ln(1+t)}\|^{-\delta /2} \wedge \ln(1+t)$.
Plugging this choice of $t_0(t)$ into \eqref{plugging} and taking $t\to\infty$, we complete the proof thanks to \ref{H2} and get the desired assertion \eqref{goalAt}
with $A_t = a_{t-t_0(t)}$.
\end{proof}

The next lemma gives the desired asymptotic of Theorem \ref{theo:kolmogorov} albeit for $(a_t, t\geq0)$. Nonetheless, taking account of the previous lemma, it completes the proof of Theorem \ref{theo:kolmogorov}. Once again, we take inspiration from arguments given in the specific setting of \cite{ren2020limit}; cf. Proposition 3.3 therein.

\begin{lemma}\label{lem:lemat}Suppose that {\ref{H1}-\ref{H5}} hold.
Then
	\[
	a_t \sim t^{-1/\alpha} \tilde{\ell}(t), \quad \text{as} \quad t\to\infty,
	\]
	where $\tilde{\ell}$ is a slowly varying function at $\infty$. In particular, $\tilde{\ell}(t) = 1 / L^\ast (t^{1/\alpha})$ where $L^*$ is the Bruijn conjugate of $s\mapsto \alpha^{-1/\alpha} \ell(1/s)^{-1/\alpha}$.
	
\end{lemma}

\begin{proof}
The proof is based on analysing the asymptotic behaviour of the inverse of the map $t \mapsto a_t$, $t\geq0$, and then applying a classical result from slowly varying functions to get the asymptotic behaviour of $a_t$. 

{Obverse that}
{
	there exists $t_{1/2}>0$ such that $\|a_t^{-1}\varphi^{-1}\su_t - 1\|\leq 1/2$ for all $t\geq t_{1/2}$ by Lemma \ref{lem:asymu}. Therefore, {by Lemma~\ref{lem:Amono}}
\begin{equation}\label{eq:derivat}
	\frac{\mathrm{d} a_t}{\mathrm{d} t} = - \langle \textsf{A}[\su_t], \tilde{\varphi} \rangle \leq - \langle \textsf{A}[ \tfrac{1}{2} a_t \varphi], \tilde{\varphi} \rangle, \qquad t\geq t_{1/2}.
\end{equation}
In particular, $\langle \textsf{A}[ \tfrac{1}{2} a_t \varphi], \tilde{\varphi} \rangle >0$, otherwise we are in contradiction with {\ref{H4}}. This implies that $a_t$ is a strictly decreasing function in $t\in (t_{1/2},\infty)$.}
Further, 
\begin{equation*}
	a_t = \langle \su_t, \tilde{\varphi} \rangle = \langle \varphi^{-1}\su_t \varphi, \tilde{\varphi} \rangle  \le \| \varphi^{-1} \su_t\| \langle \varphi, \tilde{\varphi} \rangle =   \| \varphi^{-1} \su_t\| \to 0, \quad \text{as}\quad t \to \infty,
\end{equation*}
where the convergence follows by \eqref{uniformutto0}. Hence, the mapping $t\mapsto a_t$ has an inverse on 
${(0, a_{t_{1/2}})}$
and we denote its inverse by 
$R: {(0, a_{t_{1/2}})} \to ({t_{1/2}},\infty)$. 
Note that $R$ is also strictly decreasing on 
$u\in {(0, a_{t_{1/2}})}$
and $R(u)\to \infty$ as $u\to 0$. Now, let us introduce the following function
\begin{equation}
\epsilon_t(x)= \frac{\su_t(x)}{a_t \varphi(x)}-1, \quad  t>0, \quad x\in E.
\label{epsilont}
\end{equation}
From this definition and using that $R$ is the inverse of $a_t$, we have
\begin{equation}\label{eq:uepsilon}
\su_t(x)= (\epsilon_{R(a_t)}(x) +1)a_t \varphi(x), \quad t>{t_{1/2}}
\quad x\in E.
\end{equation}
Thus, by \eqref{eq:uepsilon} and \eqref{eq:derivat}, we obtain for ${t_{1/2}}< t\le s$ 
\begin{equation*}
	\begin{split}
	s-t &= \int_{t}^{s} \mathrm{d} r = \int_{s}^{t}  \langle \textsf{A}[\su_r], \tilde{\varphi} \rangle ^{-1}  \mathrm{d} a_r = \int_{s}^{t} \langle \textsf{A}[(\epsilon_{R(a_r)} +1)a_r \varphi],  \tilde{\varphi} \rangle ^{-1}  \mathrm{d} a_r \\ & =   \int_{a_s}^{a_t} \langle \textsf{A}[(\epsilon_{R(u)} +1)u \varphi],  \tilde{\varphi} \rangle ^{-1}  \mathrm{d} u.
	\end{split}
\end{equation*}
Letting $t\to {t_{1/2}}$,  
we get
\begin{equation*}
			s =  \int_{a_s}^{a_{{t_{1/2}}}} \langle \textsf{A}[(\epsilon_{R(u)} +1)u \varphi],  \tilde{\varphi} \rangle ^{-1}  \mathrm{d} u, \quad s\in ({t_{1/2}},\infty).
\end{equation*}
Now, using that $R$ is the inverse of $a_s$, we obtain (by replacing $s=R(a_s)$ and $a_s$ by $r$)
\begin{equation*}
		R(r) =  \int_{r}^{a_{{t_{1/2}}}} \langle \textsf{A}[(\epsilon_{R(u)} +1)u \varphi],  \tilde{\varphi} \rangle ^{-1}  \mathrm{d} u, \quad {r \in (0, a_{t_{1/2}})}.
\end{equation*}
Furthermore, appealing to Lemma \ref{lem:asymu} and taking into account the fact that $R(u)\to \infty$ as $u\to 0$, we deduce that
\[\sup_{x\in E}|\epsilon_{R(u)}(x)| \to 0,\quad \text{as}\quad u\to 0.\]
In turn, combining this with  {\ref{H4}}, we obtain 
\[  \langle \textsf{A}[(\epsilon_{R(u)} +1)u \varphi],  \tilde{\varphi} \rangle  \sim   ((\epsilon_{R(u)} +1) u)^{1+\alpha} \ell((\epsilon_{R(u)} +1)u) \sim u^{1+\alpha} \ell(u), \quad \text{as}\quad u\to 0.\]
Thus 
\begin{equation*}
	R(r) =  \int_{r}^{a_{{t_{1/2}}}} \langle \textsf{A}[(\epsilon_{R(u)} +1)u \varphi],  \tilde{\varphi} \rangle ^{-1}  \mathrm{d} u \sim  \int_{r}^{a_{{t_{1/2}}}} u^{-(1+\alpha)} \ell(u)^{-1}  \mathrm{d} u, \quad \text{as}\quad r\to 0.
\end{equation*}
Since $1/\ell(u) $ is also a slowly varying function at $0$ (see e.g. Proposition 1.3.6 in \cite{bingham1989regular}), appealing to Corollary 2.3 in \cite{ren2020limit}, we get 
\begin{equation*}
 \int_{r}^{a_0} u^{-(1+\alpha)} \ell(u)^{-1} \mathrm{d} u = - \frac{1}{\alpha}  \int_{r}^{a_0} \ell(u)^{-1} \mathrm{d} u^{-\alpha} \sim  \frac{1}{\alpha r^{\alpha} \ell(r)}, \quad \text{as}\quad r\to 0.
\end{equation*}
\begin{equation*}
 \int_{r}^{a_{{t_{1/2}}}} u^{-(1+\alpha)} \ell(u)^{-1} \mathrm{d} u = - \frac{1}{\alpha}  \int_{r}^{a_{{t_{1/2}}}} \ell(u)^{-1} \mathrm{d} u^{-\alpha} \sim  \frac{1}{\alpha r^{\alpha} \ell(r)}, \quad \text{as}\quad r\to 0.
\end{equation*}
 In particular, 
\begin{equation}\label{eq:atlat}
t = R(a_t) \sim \frac{1}{\alpha a_t^\alpha \ell(a_t)},\quad \text{as} \quad  t\to \infty.
\end{equation}

Now, our goal is to apply Proposition 1.5.15 in \cite{bingham1989regular} to get the asymptotic behaviour of the map $r\mapsto a_r$ as $r\to \infty$ and the exact relationship between the slowly varying functions associated to $R(r)$ and $a_r$. To this end, let us define $S(s)=R(1/s)$. Then
\[
S(s) \sim \frac{s^\alpha}{\alpha \ell(1/s)} = s^\alpha L(s)^\alpha
\]
as $s\to \infty$, where $L(s)=(\alpha \ell(1/s))^{-1/\alpha}$ is slowly varying at $\infty$. 
Taking into account that, asymptotically as $t\to\infty$, $S(1/a_t)=t$, by Proposition 1.5.15 in \cite{bingham1989regular},
\begin{equation}\label{eq:asympat}
a_t \sim t^{-1/\alpha} \tilde{\ell}(t), \qquad \text{as } t\to\infty,
\end{equation}
where $\tilde{\ell}(t) = 1 / L^\ast (t^{1/\alpha})$ is a slowly varying function at $\infty$ and $L^\ast$ is the Bruijn conjugate of $L$.
\end{proof}

We are now ready to show Theorem \ref{theo:kolmogorov} which follows directly from Lemmas \ref{lem:asymu} and \ref{lem:lemat}.

\begin{proof}[Proof of Theorem \ref{theo:kolmogorov}]
On the one hand, by Lemma \ref{lem:asymu}, we have
\[\sup_{x\in E} \left| \frac{\su_t(x)}{a_t \varphi(x)} - 1 \right| \to 0, \quad \text{as}\quad t\to \infty.\]
On the other hand, Lemma  \ref{lem:lemat} tells us that
\[
	a_t \sim t^{-1/\alpha} \tilde{\ell}(t), \quad \text{as}\quad  t\to\infty,
\]
Combining these two results, we obtain
\[\sup_{x\in E} \left| \frac{\su_t(x)}{t^{-1/\alpha} \tilde{\ell}(t)\varphi(x)} - 1 \right| \to 0, \quad \text{as}\quad t\to \infty,\]
and therefore,
\[
\frac{\mathbb{P}_\mu (\zeta > t)}{t^{-1/\alpha} \tilde{\ell}(t)} = {\frac{1 - \mathrm e^{ \langle \ln(1-\su_t), \mu \rangle} }{t^{-1/\alpha} \tilde{\ell}(t)} \sim } \frac{\langle \su_t,\mu\rangle }{t^{-1/\alpha} \tilde{\ell}(t)} \to \langle \varphi, \mu \rangle, \quad \text{as}\quad  t \to \infty,
\]
as desired.
\end{proof}

\subsection{Non-local Superprocesses}
The proof of Theorem \ref{theo:kolmogorov} can be deduced from the analogues of the two Lemmas \ref{lem:asymu} and \ref{lem:lemat}. The reader can easily verify that the analogue of Lemma \ref{lem:lemat} is  true. This requires replacing the role of $\su_t$ by $\sV_t$ and making use of the equation \eqref{atJ} and \eqref{J-H5} in place of \eqref{eq:at} and \eqref{A},  the proof is otherwise is almost verbatim the same.

The analogue of Lemma \ref{lem:asymu} on the other hand requires some adjustment in the proof. We state and prove it next for the sake of clarity.
\begin{lemma}\label{lem:asymV}
Under the assumptions {\ref{H2}, \ref{H3}, \ref{H5} and \ref{H6}}, it holds that 
\begin{equation}
\sup_{x\in E} \left| \frac{\emph{\sV}_t(x)}{a_t \varphi(x)} - 1 \right| \to 0, \quad \text{as}\quad t\to \infty.
\end{equation}
\end{lemma}

By analogy with the proof of Lemma \ref{lem:asymu}, a core part of the proof uses the spine decomposition of $(X,\mathbb{P})$. It transpires that \cite{RenSongYang2022} offers a  version of the spine decomposition which is slightly less general. Although Remark 4.8 in that paper gives insight into the more general setting we deal with in this paper, we give a precise statement below, whose proof is dealt with in the \hyperref[appn]{Appendix}. The reader can also familiarise themselves with other less general versions of the spine decomposition for (local) superprocesses in \cite{Alison1992,  Evans1993, LRSS2023, ren2020limit}. 

\begin{theorem}\label{superspine}{Assume \ref{H2} and \ref{H6}}. Under the change of measure \eqref{COM}, for each $x\in E$,  the process $(X, \mathbb{P}^\varphi_{\delta_x})$ is equivalent in law to the aggregated mass from following decomposition:
\begin{itemize}

 \item An independent copy of $(X, \mathbb{P}_{\delta_x})$ is issued. 
 \item A marked particle, called the `spine' is issued  with Markovian evolution, denoted $\xi: = (\xi_t, t\geq0)$, determined by the semigroup
\begin{equation}
\hat{\emph{\sP}}_t^\varphi[f](x) := 
\frac{\emph{\sT}_t[\varphi f](x)}{\varphi(x)}
 \label{motion-bias2}
\end{equation}
for $x\in E$ and $f\in B^+(E)$. The  semigroup $(\hat{\emph{\sP}}_t^\varphi, t\geq0)$ is described in terms of the original Markov process $(\xi, \mathbf{P})$ associated to the underlying motion semigroup $(\emph{\sP}_t, t\geq0)$ of our non-local superprocess as follows.  Consider a process  $\hat{\xi}:= (\hat{\xi}_t, t\geq0)$
with probabilities $\hat{\mathbf{P}} = (\hat{\mathbf{P}}_x, x\in E)$ that experiences the same Markovian increments as  $(\xi, \mathbf{P})$, albeit interlaced with  `special jumps' that arrive at instantaneous rate $ \beta(x)\emph{\sm}[\varphi](x)/\varphi(x)$, $x\in E$ and with jump distribution given by the kernel $y\mapsto\kappa(x,\dd y)$, $x,y\in E$, where 
\[
\int_E f(y)\kappa(x, \dd y) =\frac{1}{\emph{\sm}[\varphi](x)}\Big(\gamma(x,  f\varphi) + \int_{M(E)^\circ}\langle f\varphi, \nu\rangle\Gamma(x, \dd \nu)\Big)=\frac{\emph{\sm}[f\varphi](x)}{\emph{\sm}[\varphi](x)}.
\]

The semigroup $\hat{\emph{\sP}}^\varphi = (\hat{\emph{\sP}}_t^\varphi, t\geq0)$ is that of a conservative Markov process that corresponds to a Doob $h$-transform of $(\hat{\xi}, \hat{\mathbf{P}})$ via 
\[
\hat{\emph{\sP}}_t^\varphi[f](x)= \hat{\mathbf{E}}_x\left[f(\hat{\xi}_t)\frac{\varphi(\hat{\xi}_t)}{\varphi(x)}{\rm e}^{\int_0^t \frac{\beta(\hat{\xi}_s)}{\varphi(\hat{\xi}_s)}(\emph{\sm}[\varphi](\hat{\xi}_s)-\varphi(\hat{\xi}_s)) + b(\hat{\xi}_s)\dd s}\right], \quad x\in E, t\geq0, f\in B^+(E).
\]
Henceforth we suppose that $(\hat{\xi}^\varphi_t, t\geq0)$ is a realisation of the $\hat{\emph\sP}^\varphi$-spine.

\item Given  $(\hat{\xi}^\varphi_t, t\geq0)$, at rate 
\[
2c(\hat{\xi}^\varphi_t)\mathbb{N}_{\hat{\xi}^\varphi_t}(\dd \omega)\dd t
\]
the ${M}(E)$-valued coordinate process $\omega = (\omega_s, s\geq0)$ is immigrated at the space-time point $(\hat{\xi}^\varphi_t, t)$, where 
 $\mathbb{N}_x$ is the Dynkin--Kuznetsov measure on ${M}(E)$-path space associated to $(\emph{\sV}_t[g], t\geq0)$, $g\in B^+(E)$, via the relation
\[
\emph{\sV}_t[g](x) = \mathbb{N}_x(1- {\rm e}^{-\langle g, \omega_t\rangle}), \qquad t\geq0.
\]

\item Given $(\hat{\xi}^\varphi_t, t\geq0)$, at rate 
\[
y\nu({\hat{\xi}^\varphi_t}, \dd y)\mathbb{P}_{y\delta_{\hat{\xi}^\varphi_t}}(\dd \omega)\dd t
\]
the ${M}(E)$-valued coordinate process 
is immigrated at the space-time point $(\hat{\xi}^\varphi_t, t)$.

\item Suppose we denote by  $(T_n, n\in\mathbb{N})$ the special jump times of $(\hat{\xi}^\varphi_t, t\geq0)$. At the $n$-th special jump time, we immigrate additional mass in a way that is correlated to the special jump size such that the joint distribution is given by the extended kernel $\kappa(\hat{\xi}^\varphi_{T_n-}; \dd y,  \dd \nu)$, $y\in E$, $\nu\in M_{0}(E)$ where
\begin{align}
&\int_E\int_{M(E)^\circ} f(y){\rm e}^{-\langle g, \nu\rangle} \kappa(\hat{\xi}^\varphi_{T_n-}; \dd y,  \dd \nu)\notag\\
&=
\frac{1}{\emph{\sm}[\varphi](\hat{\xi}^\varphi_{T_n-})}
\Bigg(\gamma(\hat{\xi}^\varphi_{T_n-}, f\varphi)+
 \int_{M_{0}(E)}{\rm e}^{-\langle g,{\nu}\rangle} \langle f\varphi, \nu\rangle\Gamma(\hat{\xi}^\varphi_{T_n-},\dd \nu)
 \Bigg),
 \label{jump-mass}
\end{align}
 for $f,g\in B^+(E)$.

\end{itemize}
\end{theorem}

\bigskip

\begin{proof}[Proof of Lemma \ref{lem:asymV}]
We start with \eqref{COM}. Recalling \eqref{surv-prob-sup} and \eqref{atsup}, and in a slight difference with \eqref{notethat}, we have that
\begin{equation}
\mathbb{E}^\varphi_{\delta_x}[1/Y_t]  = \varphi(x)^{-1}\mathbb{P}_{\delta_x}(\zeta>t) \sim \varphi(x)^{-1}\sV_t(x)\qquad \text{and}\qquad \mathbb{E}^\varphi_{\varphi\tilde\varphi}[1/Y_t]  \sim a_t.
\label{notethat2}
\end{equation}
as $t\to\infty$.

Based on the description of the spine in the bullet points above, suppose we define $\mathbf{N}(\dd s, \dd \omega)$ as the optional random measure on $[0,\infty)\times {M}(E)$, which describes the immigration along the spine $({\hat{\xi}}^\varphi_t, t\geq0)$. 
More precisely, given $({\hat{\xi}}^\varphi_t, t\geq0)$ the  previsible compensator of $\mathbf{N}$ is given by 
\begin{align}
\tilde{\mathbf{N}}(\dd s,\dd \omega) &= 2c(\hat{\xi}^\varphi_s)\mathbb{N}_{\hat{\xi}^\varphi_s}(\dd \omega)\dd s + \int_{(0,\infty)}y\nu({\hat{\xi}^\varphi_s}, \dd y)\mathbb{P}_{y\delta_{\hat{\xi}^\varphi_s}}(\dd \omega)\dd s\notag\\
&\hspace{3cm}+\beta(\hat{\xi}^\varphi_{s})\frac{ \sm[\varphi](\hat{\xi}^\varphi_{s})}{\varphi(\hat{\xi}^\varphi_{s})}
\int_E  \kappa(\hat{\xi}^\varphi_{s}; \dd y,  \dd \nu)\mathbb{P}_\nu(\dd \omega)\dd s.
\label{bigNtilde}
\end{align}

Analogously to \eqref{Y},   for any $t>0$ and $t_0 \in (0,t)$, define
\begin{equation}\label{Yt0-super}
Y_t^{(t_0,t]} =  \int_{(t_0, t]}\int_{M(E)} \langle\varphi,\omega_{t-s}\rangle \mathbf{N}(\dd s, \dd \omega).
\end{equation}

Again, write 
\[
 \epsilon^1_x(t_0, t) = \left|\mathbb{E}^\varphi_{\delta_x}[1/Y_t^{(t_0,t]}] - \mathbb{E}^\varphi_{\varphi\tilde\varphi}[1/Y_t^{(t_0,t]}]\right| \quad \text{ and } \quad \epsilon^2_x(t_0, t) = \left|\mathbb{E}^\varphi_{\delta_x}[1/Y_t] - \mathbb{E}^\varphi_{\delta_x}[1/Y_t^{(t_0,t]}]\right| .
\]
As before, we have that, for $t_0 =t_0(t)\in(0,t)$, \eqref{claim} holds and hence it suffices to show that  
\begin{equation}
\lim_{t\to\infty} \frac{\epsilon^1_x(t_0, t)}{a_{t-t_0}} +\frac{\epsilon^2_x(t_0, t)}{a_{t-t_0}} = 0. 
\label{claim2}
\end{equation}
 
 To deal with the limit of $\epsilon_x^1(t_0, t)$,  we can appeal to  the same calculations as in  \eqref{useforepsilon2}, \eqref{othernotethat}, \eqref{e1estimate} and \eqref{e1deltabound}, so that
 \begin{equation}\label{e1super}
 	0 \leq \frac{\epsilon_x^1(t_0,t)}{a_{t-t_0}} 
	\leq \frac{\Delta_{t_0} \langle \varphi f_{t-t_0}, \tilde\varphi\rangle}{a_{t-t_0}} 
	\leq \Delta_{t_0},
 \end{equation}
 where $f_t(x) = \varphi(x)^{-1} (1-\mathrm{e}^{-\sV_t(x)})$.
 
 For the limit with $\epsilon_x^2(t_0, t)$, we note that \eqref{bigmess} and \eqref{oddLaplace} are verbatim the same. A difference occurs with the calculation \eqref{campbell}. In the current setting this is replaced by 
 \begin{align}
&\mathbb{E}^\varphi_{\delta_x}[{\rm e}^{-\lambda Y^{(0,t_0]}_t} | \mathcal{F}_{t_0}^{{\hat{\xi}^\varphi}}] \notag\\
&= \exp\Bigg(-\int_0^{t_0}2c(\hat{\xi}^\varphi_s)\mathbb{N}_{\hat{\xi}^\varphi_s}(1-{\rm e}^{-\lambda \langle\varphi, X_{t-s}\rangle} )\dd s\Bigg)\notag\\
&\hspace{1cm}\times\exp\Bigg(- \int_0^{t_0}\int_{(0,\infty)}y\nu({\hat{\xi}^\varphi_s}, \dd y)\Big(1-\mathbb{E}_{y\delta_{\hat{\xi}^\varphi_s}}[{\rm e}^{-\lambda \langle\varphi, X_{t-s}\rangle}]\Big)\dd s \Bigg)\notag\\
&\hspace{2cm}\times  \exp\Bigg(-\int_0^{t_0}\frac{ \beta(\hat{\xi}^\varphi_{s})}{\varphi(\hat{\xi}^\varphi_{s})} 
 \int_{M_{0}(E)}(1-\mathbb{E}_{\nu}[{\rm e}^{-\lambda \langle\varphi, X_{t-s}\rangle}]) \langle \varphi, \nu\rangle\Gamma(\hat{\xi}^\varphi_{s},\dd \nu)
\dd s\Bigg)\notag\\
 &= \exp\Bigg(-\int_0^{t_0}2c(\hat{\xi}^\varphi_s)\sV_{t-s}[\lambda\varphi](\hat{\xi}^\varphi_s)\dd s\Bigg)\notag\\
&\hspace{1cm}\times\exp\Bigg(- \int_0^{t_0}\int_{(0,\infty)}\Big(1-{\rm e}^{- y \sV_{t-s}[\lambda\varphi](\hat{\xi}^\varphi_s)} 
\Big)
y\nu({\hat{\xi}^\varphi_s}, \dd y)\dd s \Bigg)
\notag\\
&\hspace{2cm}\times  \exp\Bigg(-\int_0^{t_0}\frac{ \beta(\hat{\xi}^\varphi_{s})}{\varphi(\hat{\xi}^\varphi_{s})} 
 \int_{M_{0}(E)}(1-{\rm e}^{- \langle\sV_{t-s}[\lambda\varphi],  \nu\rangle}
 ) \langle \varphi, \nu\rangle\Gamma(\hat{\xi}^\varphi_{s},\dd \nu)\dd s\Bigg)
 \end{align}
where we have used \eqref{bigNtilde} and  \eqref{jump-mass} for the first equality. 

Next recall that $\sV_t(x) : = \lim_{\lambda\to\infty}\sV_t[\lambda](x)$ and hence, as in \eqref{limlam}, we have 
\begin{align}
&\mathbb{E}^\varphi_{\delta_x} [\mathbf{1}_{(Y_t^{(0,t_0]}>0)} | \mathcal{F}_{t_0}^{{\hat{\xi}^{\varphi}}}]\notag\\
&\leq
\lim_{\lambda\to\infty} - \log  \mathbb{E}^\varphi_{\delta_x}[{\rm e}^{-\lambda Y^{(0,t_0]}_t} | \mathcal{F}_{t_0}^{{\hat{\xi}^{\varphi}}}]\notag\\
&=\int_0^{t_0}2c(\hat{\xi}^\varphi_s)\sV_{t-s}(\hat{\xi}^\varphi_s)\dd s +\int_0^{t_0}\int_{(0,\infty)}\Big(1-{\rm e}^{- y\sV_{t-s}(\hat{\xi}^\varphi_s)} \Big)
y\nu({\hat{\xi}^\varphi_s}, \dd y)\dd s\notag\\
&\hspace{2cm}+\int_0^{t_0}\frac{ \beta(\hat{\xi}^\varphi_{s})}{\varphi(\hat{\xi}^\varphi_{s})} 
 \int_{M_{0}(E)}(1-{\rm e}^{- \langle\sV_{t-s},  \nu\rangle}
 ) \langle \varphi, \nu\rangle\Gamma(\hat{\xi}^\varphi_{s},\dd \nu) \dd s.
 \label{makeitgotozero}
\end{align}

As in the branching particle setting, we want this to be bounded by some deterministic $C_t$, which is independent of $x$, such that $\lim_{t\to\infty} C_t = 0$.

To this end, let us note first that $V_t$ is non-increasing according to \eqref{surv-prob-sup}. Furthermore,
let us recall the analysis leading to \eqref{uniformutto0}. Using a similar approach, in particular starting from the equation \eqref{nonlinvJ} with $f$ replaced by $\infty$, we can similarly show that 
\begin{equation}
\lim_{t\to\infty}\norm{\varphi^{-1}\sV_t}=0 \qquad \text{and hence} \qquad \lim_{t\to\infty}\norm{\sV_t}=0.
\label{uniformutto0-super}
\end{equation}
Since $c\in B^+(E)$, it follows that we can easily conclude that, for the first term in \eqref{makeitgotozero},
\begin{equation}
\int_0^{t_0}2c(\hat{\xi}^\varphi_s)\sV_{t-s}(\hat{\xi}^\varphi_s)\dd s \leq 2  \|c\| \cdot t_0 \norm{\sV_{t-t_0}}.
\label{C1}
\end{equation}

Recall Remark \ref{rem1} and also that, in \ref{H5}, we assumed in particular that $ y^{1+\delta} \nu(x,\dd y)$ is uniformly integrable for all $x\in E$ and some $\delta \in (0,1)$. This together with the deterministic inequalities $1-\mathrm{e}^{-z} \leq z$  and $1-\mathrm{e}^{-z} \leq z^\delta$ if {$z\geq 0$},
we get
	\begin{align}
		\int_0^{t_0}\int_{(0,\infty)}\Big(1-{\rm e}^{- y\sV_{t-s}(\hat{\xi}^\varphi_s)} \Big)
		y\nu({\hat{\xi}^\varphi_s}, \dd y)\dd s 
		&\leq 
		\int_0^{t_0}\int_{(0,1]} y\sV_{t-s}(\hat{\xi}^\varphi_s)y\nu({\hat{\xi}^\varphi_s}, \dd y)\dd s
		\notag\\
		&\quad + \int_0^{t_0}\int_{{(1,\infty)}} \Big(y\sV_{t-s}(\hat{\xi}^\varphi_s)\Big)^\delta y\nu({\hat{\xi}^\varphi_s}, \dd y)\dd s
		\notag\\
		&\le   C_1\cdot t_0 \norm{V_{t-t_0}} + C_2 \cdot t_0\norm{V_{t-t_0}}^{\delta},
		\label{C2}
\end{align}
where
\[C_1= \sup_{x\in E} \int_{(0,1]} y ^2 \nu(x,\dd y)<\infty, \qquad C_2= \sup_{x\in E} \int_{{(1,\infty)}} y ^{1+\delta} \nu(x,\dd y)<\infty.\]

Finally, note that, \ref{H5} also implies
	\begin{align}
		&\int_0^{t_0}\frac{ \beta(\hat{\xi}^\varphi_{s})}{\varphi(\hat{\xi}^\varphi_{s})}
		\int_{M_{0}(E)}(1-{\rm e}^{- \langle\sV_{t-s},  \nu\rangle}
		) \langle \varphi, \nu\rangle\Gamma(\hat{\xi}^\varphi_{s},\dd \nu)\dd s\notag\\
		&\quad\leq \int_0^{t_0}\frac{ \beta(\hat{\xi}^\varphi_{s})}{\varphi(\hat{\xi}^\varphi_{s})}
		\int_{M_{0}(E)} \mathbf{1}_{(\langle 1,\nu\rangle\leq 1)} \langle \sV_{t-s},  \nu\rangle
		 \langle \varphi, \nu\rangle\Gamma(\hat{\xi}^\varphi_{s},\dd \nu)\dd s\notag\\
		&\qquad+ \int_0^{t_0}\frac{ \beta(\hat{\xi}^\varphi_{s})}{\varphi(\hat{\xi}^\varphi_{s})}
		\int_{M_{0}(E)} \mathbf{1}_{(\langle 1,\nu\rangle > 1)}\langle \sV_{t-s},  \nu\rangle^\delta
		 \langle \varphi, \nu\rangle\Gamma(\hat{\xi}^\varphi_{s},\dd \nu)\dd s\notag\\
		 &\quad\leq \|\sV_{t-t_0}\|\int_0^{t_0}\frac{ \beta(\hat{\xi}^\varphi_{s})}{\varphi(\hat{\xi}^\varphi_{s})}
		\int_{M_{0}(E)} \mathbf{1}_{(\langle 1,\nu\rangle\leq 1)} 
		 \langle \varphi, \nu\rangle\Gamma(\hat{\xi}^\varphi_{s},\dd \nu)\dd s\notag\\
		&\qquad + \|\sV_{t-t_0}\|^\delta \int_0^{t_0} \frac{ \beta(\hat{\xi}^\varphi_{s})}{\varphi(\hat{\xi}^\varphi_{s})}	\int_{M_{0}(E)} \mathbf{1}_{(\langle 1,\nu\rangle > 1)}\langle 1,  \nu\rangle^{\delta}\langle \varphi, \nu\rangle\Gamma(\hat{\xi}^\varphi_{s},\dd \nu)\dd s
		\notag\\
		&\quad\leq K_1  \cdot t_0 \|\sV_{t-t_0}\| + K_2  \cdot t_0 \|\sV_{t-t_0}\|^{\delta}
		\label{useDCT}
\end{align}
where 
\[
K_1 = \sup_{x\in E}\frac{\beta(x)}{\varphi(x)} \int_{M_0(E)}  \mathbf{1}_{(\langle 1,\nu\rangle \leq 1)} \langle\varphi,\nu \rangle \Gamma(x, \mathrm{d}\nu)<\infty,
\]
\[
K_2 = \sup_{x\in E}\frac{\beta(x)}{\varphi(x)} \int_{M_0(E)}  \mathbf{1}_{(\langle 1,\nu\rangle > 1)} \langle 1, \nu \rangle^{\delta} \langle\varphi,\nu \rangle \Gamma(x, \mathrm{d}\nu)<\infty.
\]

Plugging \eqref{C1}, \eqref{C2} and \eqref{useDCT} into the counterpart of \eqref{bigmess}, we have
\begin{align*}
\epsilon^2_x(t_0, t)
&\leq 
\left( (2 \|c\| + C_1+K_1) \cdot t_0 \|\sV_{t-t_0}\| + (C_2 +K_2)  \cdot t_0 \|\sV_{t-t_0}\|^{\delta} \right)   \frac{\sT_{t_0}[\varphi f_{t-t_0}](x)}{\varphi(x)}\notag\\
&\leq \left( (2 \|c\| + C_1+K_1) \cdot t_0 \|\sV_{t-t_0}\| + (C_2 +K_2)  \cdot t_0 \|\sV_{t-t_0}\|^{\delta} \right) (1+\Delta_{t_0}) a_{t-t_0}.
\end{align*}
Combining the above with \eqref{e1super}, we deduce that
\begin{equation*}
 	0 \leq \frac{\epsilon_x^1(t_0,t) + \epsilon_x^2(t_0,t)}{a_{t-t_0}} 
	\leq \Delta_{t_0} + \left( (2 \|c\| + C_1+K_1) t_0 \|\sV_{t-t_0}\| + (C_2 +K_2) t_0 \|\sV_{t-t_0}\|^{\delta} \right) (1+\Delta_{t_0}).
	\end{equation*}

The proof now proceeds as in the particle setting. Since \eqref{uniformutto0-super} holds, it is possible to choose a map $t \mapsto t_0(t)$ such that,
\[
t_0(t) \to \infty
\qquad \text{and} \qquad
t_0(t)\|\sV_{t-t_0(t)}\|^{\delta} \to 0
\]
as $t\to\infty$.
This completes the justification of \eqref{claim2}, which we earlier remarked was sufficient to complete the proof of the lemma.
 \end{proof}

\section{Proof of Theorem \ref{theo:Yaglom}}

Again, we split the proof according to the setting of branching Markov processes  and superprocesses, with greater emphasis on detail for the former. 

\subsection{Non-local branching Markov processes}

We first need three auxiliary lemmas.

\begin{lemma}\label{lem:approxBMP}
Under assumptions \ref{H1}--\ref{H5}, it holds that
    $$
    \lim_{s\to\infty} \sup_{g\in B^+(E)}\frac{\left\| \varphi^{-1}\normalfont \su_{s}[\mathrm{e}^{-g}] -a_{s} [\mathrm{e}^{-g}]\right\|}{a_s} = 0,
    $$
    {where we recall that $\normalfont a_t = \langle \su_t[0],\tilde{\varphi}\rangle$.}
\end{lemma}

\begin{proof}
We follow an approach very similar to Lemma \ref{lem:asymu}. Let us fix $g\in B^+(E)$. Recalling the definition \eqref{1-u has smg property} and the spine change of measure \eqref{COM}, we have
\begin{equation}\label{eq:usuper}
\su_t[{\rm e}^{-g}](x)  = \mathbb{E}_{\delta_x}\left[1- {\rm e}^{-\langle g, X_t\rangle}\right] = \mathbb{E}^\varphi_{\delta_x}\left[\frac{\varphi(x)}{\langle \varphi, X_t\rangle}\left(1- {\rm e}^{-\langle g, X_t\rangle}\right)\right]
\end{equation}
and accordingly we define
\[
a_t[{\rm e}^{-g}]= \langle \su_t[{\rm e}^{-g}], \tilde\varphi\rangle,
\]
which from \eqref{eq:nonlin-lin} and \ref{H2}, verifies
\begin{equation}\label{aaadef}
    a_t[{\rm e}^{-g}] = \langle 1-{\rm e}^{-g}, \tilde\varphi\rangle - \int_0^t \langle \sA[\su_{s}[{\rm e}^{-g}]], \tilde\varphi\rangle \dd s.
\end{equation}
Note also that 
\begin{equation}
	\su_t[{\rm e}^{-g}](x)  = \mathbb{E}_{\delta_x}\left[1- {\rm e}^{-\langle g, X_t\rangle}\right] 
	=  \mathbb{E}_{\delta_x}\left[(1- {\rm e}^{-\langle g, X_t\rangle})\mathbf{1}_{\{\langle 1, X_t\rangle>0\}}\right] \leq \su_t(x)
	\label{u[g]<u}
\end{equation}
which also implies
\begin{equation}
	a_{t}[{\rm e}^{-g}]\leq a_t.
	\label{agggaa}
\end{equation}
For any $f\in B^+(E)$, $t>0$ and $t_0\in(0,t)$, using the spine decomposition and similarly to \eqref{Y}, we write
\[
\langle f, X_{t}^{(t_0,t]}\rangle = f(\xi_t)+\sum_{i = n_{t_0}+1}^{n_t}\sum_{j = 1, j\neq i^*}^{N^i} \langle f, X^j_{t-T_i}\rangle.
\]
Then,
{using Markov property and \eqref{eq:usuper}, we have}
\begin{align}
	\mathbb{E}^\varphi_{\varphi\tilde\varphi}\left[\frac{1}{\langle \varphi, X_{t}^{(t_0,t]}\rangle}\left(1- {\rm e}^{-\langle g, X_{t}^{(t_0,t]}\rangle}\right)\right] & = \mathbb{E}^\varphi_{\varphi\tilde\varphi}\left[\mathbb{E}^\varphi_{\varphi\tilde\varphi}\left[\frac{1}{\langle \varphi, X_{t}^{(t_0,t]}\rangle}\left(1- {\rm e}^{-\langle g, X_{t}^{(t_0,t]}\rangle}\right)\Bigg| \mathcal{F}^\xi_{t_0}\right]\right] \notag\\
	&=\mathbb{E}^\varphi_{\varphi\tilde\varphi}\left[ 
	\mathbb{E}^\varphi_{\delta_{\xi_{t_0}}}\left[\frac{1}{\langle\varphi ,X_{t-t_0}\rangle}(1- {\rm e}^{-\langle g ,X_{t-t_0}\rangle})\right]\right]\notag\\
	&=\mathbb{E}^\varphi_{\varphi\tilde\varphi}\left[
	\varphi(\xi_{t_0})^{-1}{\su_{t-t_0}}[{\rm e}^{-g}](\xi_{t_0}) 
	\right]\notag\\
	&=\langle \varphi^{-1}{\su_{t-t_0}}[{\rm e}^{-g}], \varphi\tilde\varphi\rangle= a_{t-t_0}[{\rm e}^{-g}],
	\label{aggg}
\end{align}
where we have used that $\xi_{t_0}\sim\varphi\tilde\varphi$ under $\mathbb{P}^\varphi_{\varphi\tilde\varphi}$. 
Therefore,
\begin{align*}
		&\left|
		\frac{\su_t[{\rm e}^{-g}](x)}{\varphi(x)}- a_{t-t_0}[{\rm e}^{-g}]\right|
        \\
        &\quad =
		\left|\mathbb{E}^\varphi_{\delta_x}\left[\frac{1}{\langle \varphi, X_{t}\rangle}\left(1- {\rm e}^{-\langle g, X_{t}\rangle}\right)\right]
		-\mathbb{E}^\varphi_{\varphi\tilde\varphi}\left[\frac{1}{\langle \varphi, X_{t}^{(t_0,t]}\rangle}\left(1- {\rm e}^{-\langle g, X_{t}^{(t_0,t]}\rangle}\right)\right]\right| \notag\\
		&\quad\leq \left|\mathbb{E}^\varphi_{\delta_x}\left[\frac{1}{\langle \varphi, X_{t}\rangle}\left(1- {\rm e}^{-\langle g, X_{t}\rangle}\right)\right]
		-\mathbb{E}^\varphi_{\delta_x}\left[\frac{1}{\langle \varphi, X_{t}^{(t_0,t]}\rangle}\left(1- {\rm e}^{-\langle g, X_{t}^{(t_0,t]}\rangle}\right)\right]\right|\notag\\
		&\qquad+\left|\mathbb{E}^\varphi_{\delta_x}\left[\frac{1}{\langle \varphi, X_{t}^{(t_0,t]}\rangle}\left(1- {\rm e}^{-\langle g, X_{t}^{(t_0,t]}\rangle}\right)\right]
		-\mathbb{E}^\varphi_{\varphi\tilde\varphi}\left[\frac{1}{\langle \varphi, X_{t}^{(t_0,t]}\rangle}\left(1- {\rm e}^{-\langle g, X_{t}^{(t_0,t]}\rangle}\right)\right]\right|\notag\\
		&\quad =:\epsilon^2_x[g](t_0, t) + \epsilon^1_x[g](t_0, t).
	\end{align*}
	Next, with the help of \eqref{aggg} and \eqref{agggaa}, we deduce that 
	\begin{align*}
		&\epsilon^1_x[g](t_0, t)
        \\
        & =\left|\mathbb{E}^\varphi_{\delta_x}\left[\mathbb{E}^\varphi_{\delta_{\xi_{t_0}}}\left[\frac{1}{\langle \varphi ,X_{t-t_0}\rangle}\left(1- {\rm e}^{-\langle g ,X_{t-t_0}\rangle}\right)\right]\right]
		-\mathbb{E}^\varphi_{\varphi\tilde\varphi}\left[\mathbb{E}^\varphi_{\delta_{\xi_{t_0}}}\left[\frac{1}{\langle \varphi ,X_{t-t_0}\rangle}\left(1- {\rm e}^{-\langle g ,X_{t-t_0}\rangle}\right)\right]\right]\right|\notag\\
		&=\left|\frac{\sT_{t_0}[\varphi \varphi^{-1}\su_{t-t_0}[{\rm e}^{-g}]]}{\varphi(x)}- \langle \su_{t-t_0}[{\rm e}^{-g}], {\tilde\varphi}\rangle\right|
        \leq \Delta_{t_0}\langle\su_{t-t_0}[{\rm e}^{-g}], {\tilde\varphi}\rangle
		=\Delta_{t_0} a_{t-t_0}[{\rm e}^{-g}] \leq\Delta_{t_0} a_{t-t_0},
	\end{align*}
and, as a consequence,
	\begin{equation}\label{e1gbound}
	    \sup_{x\in E, g\in B^{+}(E)}\frac{\epsilon^1_x[g](t_0, t)}{a_{t-t_0}} \leq \Delta_{t_0}.
	\end{equation}
    On the other hand, writing
	\[
	\langle f, X^{(0,t_0]}_t \rangle = \langle f,  X_t\rangle - \langle f,  X^{(t_0,t]}_t\rangle,
	\]
	we see that
	\begin{align}
		&\epsilon^2_x[g](t_0, t)\\
        & = \left|\mathbb{E}^\varphi_{\delta_x}\left[\frac{1}{\langle \varphi, X_{t}\rangle}\left(1- {\rm e}^{-\langle g, X_{t}\rangle}\right)-\frac{1}{\langle \varphi, X_{t}^{(t_0,t]}\rangle}\left(1- {\rm e}^{-\langle g, X_{t}^{(t_0,t]}\rangle}\right)\right]\right|
        \notag\\
		&=
        \left|\mathbb{E}^\varphi_{\delta_x}\left[
		\frac{\langle \varphi , X^{(t_0, t]}_t \rangle\left(1- {\rm e}^{-\langle g , X^{(0, t_0]}_t \rangle-\langle g , X^{(t_0, t]}_t \rangle}\right) - \left(\langle \varphi , X^{(0, t_0]}_t \rangle+\langle \varphi , X^{(t_0, t]}_t \rangle\right)\left(1- {\rm e}^{-\langle g , X^{(t_0, t]}_t \rangle}\right)}{\langle \varphi, X_{t}\rangle\langle \varphi, X_{t}^{(t_0,t]}\rangle}
		\right]\right|
        \notag\\
		&= 
        \left|\mathbb{E}^\varphi_{\delta_x}\left[
		\frac{\langle \varphi , X^{(t_0, t]}_t \rangle{\rm e}^{-\langle g , X^{(t_0, t]}_t \rangle}\left(1- {\rm e}^{-\langle g , X^{(0, t_0]}_t \rangle}\right) - \langle \varphi , X^{(0, t_0]}_t \rangle\left(1- {\rm e}^{-\langle g , X^{(t_0, t]}_t \rangle}\right)}{\langle \varphi, X_{t}\rangle\langle \varphi , X^{(t_0, t]}_t \rangle}
		\right]\right|
        \notag\\
		&\leq \mathbb{E}^\varphi_{\delta_x}\left[
		\frac{{\rm e}^{-\langle g , X^{(t_0, t]}_t \rangle}}{\langle \varphi , X^{(t_0, t]}_t \rangle} (1- {\rm e}^{-\langle g , X^{(0, t_0]}_t \rangle})
		\right]
		+
		\mathbb{E}^\varphi_{\delta_x}\left[
		\mathbf{1}_{\{\langle \varphi , X^{(0, t_0]}_t \rangle>0\}}\frac{1}{\langle \varphi , X^{(t_0, t]}_t \rangle}(1- {\rm e}^{-\langle g , X^{(t_0, t]}_t \rangle})
		\right]
        \notag\\
		&\leq 
		\mathbb{E}^\varphi_{\delta_x}\left[
		\mathbb{E}^\varphi_{\delta_{\xi_{t_0}}}
		\left[
		\frac{1-(1-{\rm e}^{-\langle g , X_{t-t_0} \rangle})}{\langle \varphi , X_{t-t_0} \rangle } 
		\right]
		\mathbb{E}^\varphi\left[
		(1- {\rm e}^{-\langle g , X^{(0, t_0]}_t \rangle})|\mathcal{F}^\xi_{t_0}
		\right]
		\right]\notag
		\\
		&\qquad\quad +
		\mathbb{E}^\varphi_{\delta_x}\left[
		\mathbb{E}^\varphi_{\delta_{\xi_{t_0}}}\left[\frac{1}{\langle \varphi , X_{t-t_0} \rangle}(1- {\rm e}^{-\langle g , X_{t-t_0} \rangle})\right]
		\mathbb{E}^\varphi\left[
		\mathbf{1}_{\{\langle \varphi , X^{(0, t_0]}_t \rangle>0\}}|\mathcal{F}^\xi_{t_0}
		\right]
		\right].
		\label{bigRHS}
	\end{align}
	To complete the calculation in \eqref{bigRHS}, we note that the calculations that lead to \eqref{condind} equally give us the same upper bound for $\mathbb{E}^\varphi_{\delta_x}[
	\mathbf{1}_{\{\langle g, X^{(0,t_0]}_t\rangle >0\}}|\mathcal{F}^\xi_{t_0}
	]$, which leads us to 
	\begin{align*}
		\mathbb{E}^\varphi_{\delta_x}\left[
		(1- {\rm e}^{-\langle g , X^{(0, t_0]}_t \rangle})|\mathcal{F}^\xi_{t_0}
		\right] & = \mathbb{E}^\varphi_{\delta_x}\left[
		(1- {\rm e}^{-\langle g , X^{(0, t_0]}_t \rangle}) \mathbf{1}_{\{\langle g , X^{(0, t_0]}_t \rangle>0\}}\mid\mathcal{F}^\xi_{t_0}
		\right]\notag\\
		&\leq 
		\mathbb{E}^\varphi_{\delta_x}\left[
		\mathbf{1}_{\{\langle g , X^{(0, t_0]}_t \rangle>0\}}\mid \mathcal{F}^\xi_{t_0}
		\right]\notag\\
		&\leq Ct_0\|\su_{t-t_0}\|^\delta.
	\end{align*}
    Thus, from \eqref{bigRHS}, we get
    \begin{align}
		\epsilon^2_x[g](t_0, t)
		&\leq 
		\mathbb{E}^\varphi_{\delta_x}\left[
		\varphi(\xi_{t_0})^{-1}\Big(\su_{t-t_0}(\xi_{t_0}) - \su_{t-t_0}[{\rm e}^{-g}](\xi_{t_0})\Big)\cdot Ct_0\|\su_{t-t_0}\|^\delta
		\right]
		\notag\\
		&\qquad+
		\mathbb{E}^\varphi_{\delta_x}\left[\varphi(\xi_{t_0})^{-1}\su_{t-t_0}[{\rm e}^{-g}](\xi_{t_0}) \cdot Ct_0\|\su_{t-t_0}\|^\delta
		\right]\notag\\
		&\leq Ct_0\|\su_{t-t_0}\|^\delta \frac{\sT_{t_0}[\su_{t-t_0} ](x)}{\varphi(x)}\leq Ct_0\|\su_{t-t_0}\|^\delta
		(1+\Delta_{t_0})a_{t-t_0},
		\notag
	\end{align}
{where in the second inequality we have used \eqref{spineunderspine} and in the last inequality \ref{H2}.}
	In conclusion,
	\begin{equation}\label{e2gbound}
	    \sup_{x\in E, g\in B^+(E)}\frac{\epsilon^2_x[g](t_0, t)}{a_{t-t_0}}\leq Ct_0\|\su_{t-t_0}\|^\delta
	(1+\Delta_{t_0}).
	\end{equation}
    Combining \eqref{e1gbound} and \eqref{e2gbound}, we can choose $t_0=t_0(t)$ with 
	\[
	t_0 \to \infty 
	\qquad \text{and} \qquad
	t-t_0\to\infty
	\qquad \text{and} \qquad
	t_0\|\su_{t-t_0}\|^{\delta} \to 0, \qquad \text{as } t\to\infty,
	\]
    so that
    \begin{equation}
		\lim_{t\to\infty}\sup_{g\in B^+(E)}\frac{1}{a_{t-t_0}}\left\|\frac{\su_t[{\rm e}^{-g}]}{\varphi}- a_{t-t_0}[{\rm e}^{-g}]\right\| = 0.
		\label{bootstrap1}
	\end{equation}
    For instance, take again $t_0 =  \|\su_{t-\ln(1+t)}\|^{-\delta /2} \wedge \ln(1+t)$.
    Finally, taking into account that
    \begin{align*}
		\frac{\|\varphi^{-1}\su_t[{\rm e}^{-g}]-a_t[{\rm e}^{-g}]\|}{a_t} \leq \frac{a_{t-t_0}}{a_t} \frac{\|\varphi^{-1}\su_t[{\rm e}^{-g}]-a_{t-t_0}[{\rm e}^{-g}]\| + \left(a_{t-t_0}[{\rm e}^{-g}]-a_t[{\rm e}^{-g}]\right)}{a_{t-t_0}}
	\end{align*}
	and {appealing to \eqref{eq:usuper}, \eqref{agggaa} and Lemma \ref{lem:Amono}, we get}
	\begin{align*}
		0 \leq \frac{a_{t-t_0}[{\rm e}^{-g}]-a_t[{\rm e}^{-g}]}{a_{t-t_0}} = \frac{1}{a_{t-t_0}} \int_{t-t_0}^t \langle \sA [ \su_s [{\rm e}^{-g}]] , \tilde\varphi \rangle \dd s \leq \frac{1}{a_{t-t_0}} \int_{t-t_0}^t \langle \sA [ \su_s ] , \tilde\varphi \rangle \dd s = \frac{a_{t-t_0}-a_t}{a_{t-t_0}}.
	\end{align*}
	We  conclude the proof of the result because 
	\[\frac{a_{t-t_0}}{a_t}\to 1  \qquad \text{and} \qquad \frac{a_{t-t_0}-a_t}{a_{t-t_0}} \to 0,\]
	as $t\to\infty$, according to Lemma \ref{lem:lemat}.
\end{proof}

\begin{lemma}\label{lem-ratio-BMP}
Assume that hypotheses \ref{H1}--\ref{H5} hold. For all $f\in B^+(E)$ with $\langle f, \tilde\varphi\rangle >0$, there exists $\kappa_0=\kappa_0(f)>0$ small enough such that
\begin{equation*}
    \lim_{t\to\infty} \left\| \frac{\normalfont \su_{\kappa t}[\mathrm{e}^{-a_t f}]}{a_{\kappa t} [\mathrm{e}^{-a_t f}]\varphi} -1\right\| = 0,
\end{equation*}
for all $\kappa\in(0,\kappa_0).$
\end{lemma}
\begin{proof}
{Let us fix $f\in B^+(E)$ with $\langle f, \tilde\varphi\rangle >0$. We begin by observing that since}
\begin{equation}\label{decomp}
    \left\| \frac{\su_{\kappa t}[\mathrm{e}^{-a_t f}]}{a_{\kappa t} [\mathrm{e}^{-a_t f}]\varphi} -1\right\| 
    =
    \frac{a_{\kappa t}}{a_t}\frac{a_t}{a_{\kappa t} [\mathrm{e}^{-a_t f}]} \frac{\left\| \frac{\su_{\kappa t}[\mathrm{e}^{-a_t f}]}{\varphi} -a_{\kappa t} [\mathrm{e}^{-a_t f}]\right\|}{a_{\kappa t}},
\end{equation}
it is enough to prove that
\begin{equation}\label{bound-aaaaa}
    0 < \lim_{t\to\infty} \frac{a_{\kappa t} [\mathrm{e}^{-a_t f}]}{a_t} \leq \langle f, \tilde\varphi\rangle,
\end{equation}
because from Lemmas \ref{lem:lemat} and \ref{lem:approxBMP} we have
$$
\lim_{t\to\infty} \frac{a_{\kappa t}}{a_t} = \kappa^{-1/\alpha} \qquad \text{and} \qquad \lim_{t\to\infty} \frac{\left\| \frac{\su_{\kappa t}[\mathrm{e}^{-a_t f}]}{\varphi} -a_{\kappa t} [\mathrm{e}^{-a_t f}]\right\|}{a_{\kappa t}}=0.
$$
Recalling \eqref{A} and \eqref{aaadef}, the right hand side of \eqref{bound-aaaaa} clearly follows from
$$
a_{\kappa t}[\mathrm{e}^{-a_t f}] \leq \langle 1-\mathrm{e}^{-a_t f},\tilde\varphi\rangle \leq \langle a_t f,\tilde\varphi\rangle.
$$
Further, from \eqref{eq:nonlin-lin} we also have
$$
\su_{r}[\mathrm{e}^{-a_t f}] \leq \mathsf{T}_r [1-\mathrm{e}^{-a_t f}] \leq \mathsf{T}_r [a_t f].
$$
In addition, it is easy to check that
$
\mathsf{A}[g] \leq \beta \mathsf{m}[g], \ g\in B_1^+(E),
$
and, from the many-to-one formula, see \cite[Lemma 8.2]{Horton2023}, we deduce that
$$
\mathsf{T}_r[1] \leq \mathrm{e}^{r \|(\mathsf{m}[1]-1)\beta\|}, \quad  r\geq 0.
$$
With all the previous inequalities {and using Lemma \ref{lem:Amono},} we are able to derive the following lower bound, {for $\tau \in (0, \kappa t)$,}
\begin{align}
    a_{\kappa t}[\mathrm{e}^{-a_t f}] 
    &= 
    a_t \langle f , \tilde\varphi\rangle - \langle \mathrm{e}^{-a_t f} -1 + a_t f , \tilde\varphi\rangle - \int_{0}^{\tau} \langle \mathsf{A} [\mathsf{u}_r[\mathrm{e}^{-a_t f}]], \tilde\varphi\rangle \dd r - \int_{\tau}^{\kappa t} \langle \mathsf{A} [\mathsf{u}_r[\mathrm{e}^{-a_t f}]], \tilde\varphi\rangle \dd r
    \notag
    \\
    &\geq
    a_t \langle f , \tilde\varphi\rangle - \frac{1}{2}a_t^2 \langle f^2 , \tilde\varphi\rangle - \int_{0}^{\tau} \langle \beta \mathsf{m} [\mathsf{u}_r[\mathrm{e}^{-a_t f}]], \tilde\varphi\rangle \dd r - \int_{\tau}^{\kappa t} \langle \mathsf{A} [\mathsf{T}_r [a_t f]], \tilde\varphi\rangle \dd r 
    \notag\\
    &\geq
    a_t \langle f , \tilde\varphi\rangle - \frac{1}{2}a_t^2 \langle f^2 , \tilde\varphi\rangle - \int_{0}^{\tau} \langle \beta \mathsf{m} [\mathsf{T}_r [a_t f]], \tilde\varphi\rangle \dd r - \int_{\tau}^{\kappa t} \langle \mathsf{A} [(1 + \Delta_r) \langle f,\tilde\varphi\rangle a_t \varphi], \tilde\varphi\rangle \dd r 
    \notag\\
    &\geq
    a_t \langle f , \tilde\varphi\rangle - \frac{1}{2}a_t^2 \langle f^2 , \tilde\varphi\rangle - a_t  \int_{0}^{\tau} \langle \beta \mathsf{m} [{\|f\|}\mathsf{T}_r [1]], \tilde\varphi\rangle \dd r \notag\\
    &\quad- \kappa t \Big \langle \mathsf{A} \big[(1 + \sup_{r\geq \tau}\Delta_r) \langle f,\tilde\varphi\rangle a_t \varphi\big], \tilde\varphi\Big\rangle 
    \notag\\
    &\geq
    a_t \langle f , \tilde\varphi\rangle - \frac{1}{2}a_t^2 \langle f^2 , \tilde\varphi\rangle - a_t \langle \beta \mathsf{m} [\|f\|], \tilde\varphi\rangle \int_{0}^{\tau}  \mathrm{e}^{r\|(\mathsf{m}[1]-1)\beta\|} \dd r \label{aktttt} \\
    &\quad - \kappa_0 t \Big \langle \mathsf{A} \big[(1 + \sup_{r\geq \tau}\Delta_r) \langle f,\tilde\varphi\rangle a_t \varphi\big], \tilde\varphi\Big \rangle, \notag
\end{align}
{where in the second inequality we have used that $\mathrm{e}^{-x}-1+x\le x^2/2$ for $x\ge 0$. Thus, by using that $a_t \to 0$ as $t\to \infty$, \ref{H2} and \eqref{eq:atlat}, we have}
\begin{align*}
    &\lim_{t\to\infty} \frac{a_{\kappa t} [\mathrm{e}^{-a_t f}]}{a_t} \\
    &\quad\geq 
    \langle f , \tilde\varphi\rangle -  \langle \beta \mathsf{m} [\|f\|], \tilde\varphi\rangle \int_{0}^{\tau}  \mathrm{e}^{r\|(\mathsf{m}[1]-1)\beta\|} \dd r - \lim_{t\to\infty}\frac{\kappa_0 t}{a_t} \Big \langle \mathsf{A} \big [(1 + \sup_{r\geq \tau}\Delta_r) \langle f,\tilde\varphi\rangle a_t \varphi \big], \tilde\varphi \Big \rangle 
    \\
    &\quad\geq 
    \langle f , \tilde\varphi\rangle -  \langle \beta \mathsf{m} [\|f\|], \tilde\varphi\rangle \int_{0}^{\tau}  \mathrm{e}^{r\|(\mathsf{m}[1]-1)\beta\|} \dd r -  \frac{\kappa_0(1 + \sup_{r\geq \tau}\Delta_r)^{1+\alpha} \langle f,\tilde\varphi\rangle^{1+\alpha}}{\alpha}.
\end{align*}
Choosing $\tau=\tau(f)>0$ small enough such that
$$\langle \beta \mathsf{m} [\|f\|], \tilde\varphi\rangle \int_{0}^{\tau}  \mathrm{e}^{r\|(\mathsf{m}[1]-1)\beta\|} \dd r < \frac{1}{3}\langle f, \tilde\varphi\rangle,$$
and then taking $\kappa_0=\kappa_0(f)>0$ small enough so that 
$$\frac{\kappa_0(1 + \sup_{r\geq \tau}\Delta_r)^{1+\alpha} \langle f,\tilde\varphi\rangle^{1+\alpha}}{\alpha}<\frac{1}{3}\langle f, \tilde\varphi\rangle,$$
we conclude
$$\lim_{t\to\infty} \frac{a_{\kappa t} [\mathrm{e}^{-a_t f}]}{a_t} \geq \frac{1}{3}\langle f, \tilde\varphi\rangle > 0,$$
as claimed in \eqref{bound-aaaaa}.
\end{proof}

\begin{lemma}\label{lem:ut0} Assume \ref{H1}--\ref{H5}. For all $f\in B^+(E)$ 
and $\varepsilon\in(0,1)$, there exists $\kappa=\kappa(f,\varepsilon)\in(0,\varepsilon)$ small enough such that
\begin{equation}\label{upperandlower}
    (1-\varepsilon)^2\langle f, \tilde\varphi\rangle a_t\varphi \leq \normalfont \su_{\kappa t}[\mathrm{e}^{-a_t f}] \leq (1+\varepsilon)\langle f, \tilde\varphi\rangle a_t\varphi
\end{equation}
for all $t$ sufficiently large.
\end{lemma}
\begin{proof}
Fix $f\in B^+(E)$ and $\varepsilon\in(0,1)$. For the  upper bound, recall from \ref{H2} that
	\begin{equation}
	\langle f,\tilde{\varphi}\rangle (1- \Delta_{s})\varphi
	\leq 
	\sT_{s} [f] 
	\leq 
	\langle g,\tilde{\varphi}\rangle (1+  \Delta_{s})\varphi, 
	\qquad
	f\in B^+(E), \quad s > 0,
	\label{useupperandlower}
	\end{equation}
	and, hence,  from \eqref{eq:nonlin-lin} and the upper bound in \eqref{useupperandlower},
	\[
	0 
    \leq
    \su_{\kappa t} [{\rm e}^{-a_t f}] 
	\leq 
	\sT_{\kappa t} [1- {\rm e}^{- a_t f}] 
	\leq
	\sT_{\kappa t} [a_t f] 
	\leq 
	 (1+ \Delta_{\kappa t})\langle f ,\tilde{\varphi}\rangle a_t\varphi.
	\]
	Since $\Delta_s \to 0$ as $s\to\infty$, we can take {$t$} large enough such that 
	${\Delta_{{\kappa t}} \leq \varepsilon}$, and the upper bound in \eqref{upperandlower} follows.

    If $\langle f, \tilde\varphi \rangle =0$, \eqref{upperandlower} trivially holds. Let us prove then the lower bound in \eqref{upperandlower} for the case $\langle f, \tilde\varphi \rangle >0$.
From Lemma \ref{lem-ratio-BMP}, there exists $t_0=t_0(f,\varepsilon)>0$ large enough such that
$$\left\| \frac{\su_{\kappa t}[\mathrm{e}^{-a_t f}]}{a_{\kappa t} [\mathrm{e}^{-a_t f}]\varphi} -1\right\| \leq \varepsilon, \qquad t\geq t_0.$$
In particular,
$$
\su_{\kappa t}[\mathrm{e}^{-a_t f}] \geq (1-\varepsilon)a_{\kappa t}[\mathrm{e}^{-a_t f}]\varphi , \qquad t\geq t_0,
$$
so let us show that $a_{\kappa t}[\mathrm{e}^{-a_t f}] \geq (1-\varepsilon) \langle f, \tilde\varphi\rangle a_t$ to conclude the proof. Indeed, appealing to \eqref{aktttt},
\begin{align}\label{lower0}
    a_{\kappa t}[\mathrm{e}^{-a_t f}] 
    &\geq 
    a_t \langle f , \tilde\varphi\rangle - \frac{1}{2}a_t^2 \langle f^2 , \tilde\varphi\rangle - a_t \langle \beta \mathsf{m} [\|f\|], \tilde\varphi\rangle \int_{0}^{\tau}  \mathrm{e}^{r\|(\mathsf{m}[1]-1)\beta\|} \dd r 
    \\
    &\quad - \kappa t \Big \langle \mathsf{A} \big[(1 + \sup_{r\geq \tau}\Delta_r) \langle f,\tilde\varphi\rangle a_t \varphi\big], \tilde\varphi \Big\rangle. \notag
\end{align}
We take $t$ large enough such that
	\[
	a_t \leq \frac{\frac{\varepsilon}{3} \langle f, \tilde\varphi\rangle}{\frac{1}{2}\langle f^2, \tilde\varphi\rangle},
	\]
	or equivalently,
		\begin{equation}
	\frac{1}{2}a_t^2 \langle f^2 , \tilde\varphi\rangle 
	\leq 
	\frac{\varepsilon}{3} \langle f,\tilde{\varphi}\rangle a_t.
\label{lower1}
	\end{equation}
{Note that we are using that $a_t\to 0$ as $t\to \infty$, which is valid due to assumption \ref{H3}.}
Now we choose $\tau=\tau(f,\varepsilon)>0$ small enough so that
\begin{align}\label{lower2}
    \langle \beta \mathsf{m} [\|f\|], \tilde\varphi\rangle \int_{0}^{\tau}  \mathrm{e}^{r\|(\mathsf{m}[1]-1)\beta\|} \dd r \leq \frac{\varepsilon}{3} \langle f , \tilde\varphi\rangle.
\end{align}
Finally, applying \ref{H4} and \eqref{eq:atlat} from Lemma \ref{lem:lemat}, which requires \ref{H1}--\ref{H5}, we conclude that, as $t\to\infty$,
\begin{align*}
	{\kappa t} \Big \langle {\mathsf{A}} \big[(1 + \sup_{r\geq \tau}\Delta_r) \langle f,\tilde\varphi\rangle a_t \varphi\big], \tilde\varphi \Big\rangle  
	&\sim
    {\kappa t} (1 + \sup_{r\geq \tau}\Delta_r)^{1+\alpha} \langle f,\tilde\varphi\rangle^{1+\alpha}a_t^{1+\alpha} \ell(a_t)
	\notag\\
	&\sim
	a_t \langle f,\tilde\varphi\rangle  \frac{\kappa (1 + \sup_{r\geq \tau}\Delta_r)^{1+\alpha} \langle f,\tilde\varphi\rangle^{\alpha}}{\alpha}. 
	\end{align*}
Therefore, we can take $\kappa=\kappa(f,\varepsilon)\in(0,\varepsilon)$ small enough so that
\begin{equation}\label{lower3}
\frac{\kappa (1 + \sup_{r\geq \tau}\Delta_r)^{1+\alpha} \langle f,\tilde\varphi\rangle^{\alpha}}{\alpha} < \frac{\varepsilon}{3}.    
\end{equation}
Then, for all $t$ sufficiently large, using \eqref{lower1}, \eqref{lower2} and \eqref{lower3} in \eqref{lower0}, we have
$
a_{\kappa t}[\mathrm{e}^{-a_t f}] \geq (1-\varepsilon) a_t \langle f , \tilde\varphi\rangle,
$
as desired.
\end{proof}

We are now ready for the proof of our second main result for non-local branching Markov processes.
\begin{proof}[Proof of Theorem~\ref{theo:Yaglom}]
	It is enough to show {the result with  $\mu =  \delta_x$ for some $x\in E$, that is} 
\begin{equation}
	\lim\limits_{t\to \infty} \frac{\mathbb{E}_{\delta_x}\left[1-\exp\left(-a_t \langle f,  X_t\rangle \right) \right]}{\mathbb{P}_{\delta_x}(\zeta >t)} 
	=  
	\lim_{t\to\infty} \frac{\su_t[{\rm e}^{- a_t f}](x)}{\su_t(x)}
	=
	\frac{\langle f, \tilde{\varphi} \rangle}{(1+ \langle f, \tilde{\varphi} \rangle ^{\alpha})^{1/\alpha}}.
	\label{itsenoughtoprove}
\end{equation}

	{Indeed, let $\mu \in N(E)$ then $\mu= \sum_{i=1}^{n}\delta_{x_i}$ for some $n\in \mathbb{N}$ and $x_1, \dots, x_n \in E$. Thus
\begin{align*}
\frac{\mathbb{E}_\mu [1-\mathrm e^{- a_t \langle f, X_t\rangle}]}{\mathbb{P}_\mu(\zeta >t)}
&= 
\frac{ 1 - \mathrm e^{ \langle \ln(1-\su_t[{\rm e}^{- a_t f}]), \mu \rangle}}{1 - \mathrm e^{ \langle \ln(1-\su_t), \mu \rangle}}
= 
\frac{1-\prod_{i=1}^n (1-\su_t[{\rm e}^{- a_t f}](x_i))} {1-\prod_{i=1}^n (1-\su_t(x_i))}
\\
&\sim
\frac{\sum_{i=1}^n \su_t[{\rm e}^{- a_t f}](x_i)} {\sum_{i=1}^n \su_t(x_i)} 
= 
\frac{\sum_{i=1}^n  \su_t(x_i) \frac{\su_t[{\rm e}^{- a_t f}] (x_i)}{\su_t(x_i)}} {\sum_{i=1}^n \su_t(x_i)} 
\to 
\frac{\langle  f, \tilde{\varphi} \rangle}{(1+ \langle  f, \tilde{\varphi} \rangle ^{\alpha})^{1/\alpha}},
\end{align*}
as $t \to \infty,$ where in the last convergence we have used that \eqref{itsenoughtoprove} holds for any $x_i \in E$.} {Note also that the asymptotic equivalence $\sim$ is coming from the Taylor theorem.}

{Now for $f\in B^+(E)$ and $\varepsilon\in(0,1)$, let us define the integer sequences
\begin{equation*}
    k(t) = \max \left\{ m \in \mathbb{N} : a_m \geq \frac{1+\varepsilon}{1-\varepsilon} \langle f, \tilde\varphi\rangle a_t \right\}, \; n(t) = \min \left\{ m \in \mathbb{N} : a_m \leq \frac{(1-\varepsilon)^2}{1+\varepsilon} \langle f, \tilde\varphi\rangle a_t \right\}.
\end{equation*}
By definition, we have
\begin{equation}\label{ineq-ak}
    a_{k(t)+1} \leq \frac{1+\varepsilon}{1-\varepsilon} \langle f, \tilde\varphi\rangle a_t \leq a_{k(t)}\quad \text{and}\quad   a_{n(t)} \leq \frac{(1-\varepsilon)^2}{1+\varepsilon} \langle f, \tilde\varphi\rangle a_t \leq a_{n(t)-1}.
\end{equation}
Further, $n(t)\to\infty$ and $k(t)\to\infty$ as $t\to\infty$. Moreover, from \eqref{ineq-ak}, we get
\begin{equation}
	\lim_{t\to\infty} \frac{a_{k(t)}}{a_{t}}= \langle f, \tilde{\varphi} \rangle \frac{1+\varepsilon}{1-\varepsilon}\quad \text{and}\quad 	\lim_{t\to\infty} \frac{a_{n(t)}}{a_{t}}= \langle f, \tilde{\varphi} \rangle \frac{(1-\varepsilon)^2}{1+\varepsilon}.
	\label{akt/at}
\end{equation}
Combining Lemmas \ref{lem:asymu} and \ref{lem:ut0} with inequalities in \eqref{ineq-ak}, we deduce for $t$ large enough that}
\begin{equation}\label{upperu}
    \su_{\kappa t}[\mathrm{e}^{-a_t f}] \leq (1+\varepsilon)\langle f, \tilde\varphi\rangle a_t\varphi\leq (1-\varepsilon)a_{k(t)}\varphi \leq \su_{k(t)}
\end{equation}
and
\begin{equation}\label{loweru}
    \su_{\kappa t}[\mathrm{e}^{-a_t f}] \geq (1-\varepsilon)^2\langle f, \tilde\varphi\rangle a_t\varphi \geq (1+\varepsilon)a_{n(t)}\varphi \geq \su_{n(t)}.
\end{equation}
Putting \eqref{upperu} and \eqref{loweru} together, we have shown that 
		$
	\su_{n(t)} \leq \su_{\kappa t}[{\rm e}^{-a_t f}] \leq \su_{k(t)}.$
Using that $\su_t[h] \geq \su_t[g]$ if $h \leq g$ {with $h,g\in B_1^+(E)$}, we deduce
	\[
	1- \su_{(1-\kappa)t}[1-\su_{n(t)}] 
	 \geq 
	1- \su_{(1-\kappa)t}[1-\su_{\kappa t}[{\rm e}^{-a_t f}] ] 
	 \geq
	1- \su_{(1-\kappa)t}[1-\su_{k(t)}].
	\]
	When divided through by $\su_t$, taking into account that $\su_{t+s}[g]=\su_t[1-\su_{s}[g]]$, we get
	\[
	\frac{\su_{(1-\kappa)t+n(t)}}{\su_t} 
	 \leq 
	 \frac{\su_{t}[{\rm e}^{- a_t f}]}{\su_t}
	 \leq 
	 \frac{\su_{(1-\kappa)t+k(t)}}{\su_t}.
	\]
To complete the proof, our aim is to show that 
	\begin{equation}\label{eq:akratio}
		\lim_{t\to\infty} \frac{\su_{(1-\kappa)t+k(t)}}{\su_t} 
		= 
		\lim_{t\to\infty} \frac{a_{(1-\kappa)t+k(t)}}{a_t} 
		=
		(1-\kappa+[\langle  f , \tilde{\varphi}\rangle (1+\varepsilon)/(1-\varepsilon)]^{-\alpha})^{-1/\alpha}
	\end{equation}
	\begin{equation}\label{eq:anratio}
		\lim_{t\to\infty} \frac{\su_{(1-\kappa)t+n(t)}}{\su_t} 
		= 
		\lim_{t\to\infty} \frac{a_{(1-\kappa)t+n(t)}}{a_t} 
		=
		(1-\kappa+[\langle f , \tilde{\varphi}\rangle (1-\varepsilon)^2/(1+\varepsilon)]^{-\alpha})^{-1/\alpha},
	\end{equation}
	in which case, since we may take $\varepsilon$ arbitrarily small {and $\kappa\in(0,\varepsilon)$}, we get
	\begin{equation}
		\lim_{t\to\infty} \frac{\su_{t}[{\rm e}^{- a_t f}]}{\su_t} 
		=
		(1+\langle f , \tilde{\varphi}\rangle ^{-\alpha})^{-1/\alpha},
	\end{equation}
which suffices from \eqref{itsenoughtoprove} to prove the theorem.

We thus conclude the proof by showing \eqref{eq:akratio}, noting that the proof of \eqref{eq:anratio} is essentially the same. 
	To this end, we {recall Lemma \ref{lem:lemat} and \eqref{akt/at}}.
	From Lemma 1 in \cite{slack1968branching}, we know that if $(x_n, n\ge 1)$, $(y_n, n\ge 1)$ are both positive sequences and tend to $\infty$ as $n\to \infty$,
	and for $n$ large enough there are constants $K_1,K_2$ such that
	\[
	0 < K_1 < \frac{x_n}{y_n} < K_2 < \infty,
	\]
	then for any slowly varying function $\tilde{L}$ at $\infty$, $\tilde{L}(x_n)/\tilde{L}(y_n)\to 1$ as $n\to \infty$.
	Thus, taking into account \eqref{eq:atlat} and \eqref{akt/at}, we note that 
	\[
	\lim_{t\to\infty} \frac{t}{k(t)}
	= 
	\lim_{t\to\infty} \frac{a_{k(t)}^\alpha \ell(a_{k(t)})}{a_{t}^\alpha \ell(a_t)}
	= 
	\left(\langle f, \tilde{\varphi} \rangle \frac{1+\varepsilon}{1-\varepsilon}\right)^\alpha.
	\]
	Finally, using Lemma \ref{lem:lemat} and again  Lemma 1 in \cite{slack1968branching}, we get
    \begin{align*}
        \lim_{t\to\infty} \frac{a_{(1-\kappa)t+k(t)}}{a_t} 
    &=
	\lim_{t\to\infty} \frac{[(1-\kappa)t+k(t)]^{-1/\alpha} \tilde{\ell} ((1-\kappa)t+k(t))}{t^{-1/\alpha} \tilde{\ell}(t)} 
	\\
    &=
	(1-\kappa+[\langle  f , \tilde{\varphi}\rangle (1+\varepsilon)/(1-\varepsilon)]^{-\alpha})^{-1/\alpha},
    \end{align*}
as required.
\end{proof}

\subsection{Non-local Superprocesses}

{This section is very similar to the previous one, replacing the role of $\su_t[\mathrm{e}^{-a_t f}]$ by $\sV_t[a_t f]$. Recalling the definition \eqref{non-local-transition-semigroup}, \eqref{surv-prob-sup} and \eqref{atsup}, we introduce
$$a_t[f]:=\langle \sV_t[f],\tilde\varphi\rangle, \qquad t \geq 0, \quad f\in B^+(E),$$
which according to \eqref{nonlinvJ} and \ref{H2} satisfies
\begin{equation}\label{atf}
	a_t[f] = \langle f, \tilde\varphi\rangle - \int_0^t \langle \mathsf{J}[\sV_s[f]],\tilde\varphi\rangle \dd s.
\end{equation}

Before proving Theorem \ref{theo:Yaglom}, we include the counterpart results of Lemmas \ref{lem:approxBMP}, \ref{lem-ratio-BMP} and \ref{lem:ut0}.}

\begin{lemma}\label{lem:approx}

Under assumptions \ref{H2}--\ref{H6}, it holds that
    $$
    \lim_{s\to\infty} \sup_{g\in B^+(E)}\frac{\left\| \varphi^{-1}\normalfont \sV_{s}[g] -a_{s} [g]\right\|}{a_s} = 0.
    $$
\end{lemma}
\begin{proof}

    We follow the steps of Lemma \ref{lem:approxBMP}, but now in the superprocesses setting. 
    Firstly, let us define 
    $$\sU_t[g](x) := 1 - \mathrm{e}^{-\sV_t[g](x)}, \qquad g\in B^+(E).$$
    By applying the superprocess spine change of measure from Theorem \ref{superspine}, we can write 
$$\sU_t[g](x) = \mathbb{E}_{\delta_x}\left[1 - {\rm e}^{-\langle g, X_t\rangle}\right] = \varphi(x) \mathbb{E}_{\delta_x}^\varphi \left[ \frac{1}{\langle \varphi, X_t\rangle} \left(1 - e^{-\langle g, X_t\rangle}\right) \right].$$
We also define the auxiliary integral $\tilde{a}_t[g] = \langle \sU_t[g], \tilde\varphi\rangle$. 

For a given $t_0 \in (0, t)$, we split the superprocess mass into two components, 
$$X_t = X_t^{(0, t_0]} + X_t^{(t_0, t]},$$ 
representing mass immigrated before and after $t_0$, so that for a test function $f\in B^+(E)$ we have
$$\langle f, X_t^{(t_0,t]}\rangle =  \int_{(t_0, t]}\int_{M(E)} \langle f,\omega_{t-s}\rangle \mathbf{N}(\dd s, \dd \omega),$$
in the same way as in \eqref{Yt0-super}. Taking into account that
\begin{align}
	\mathbb{E}^\varphi_{\varphi\tilde\varphi}\left[\frac{1}{\langle \varphi, X_{t}^{(t_0,t]}\rangle}\left(1- {\rm e}^{-\langle g, X_{t}^{(t_0,t]}\rangle}\right)\right] & = \mathbb{E}^\varphi_{\varphi\tilde\varphi}\left[\mathbb{E}^\varphi_{\varphi\tilde\varphi}\left[\frac{1}{\langle \varphi, X_{t}^{(t_0,t]}\rangle}\left(1- {\rm e}^{-\langle g, X_{t}^{(t_0,t]}\rangle}\right)\Bigg| \mathcal{F}^{\hat{\xi}^\varphi}_{t_0}\right]\right] \notag\\
	&=\mathbb{E}^\varphi_{\varphi\tilde\varphi}\left[ 
	\mathbb{E}^\varphi_{\delta_{{\hat{\xi}^\varphi}_{t_0}}}\left[\frac{1}{\langle\varphi ,X_{t-t_0}\rangle}(1- {\rm e}^{-\langle g ,X_{t-t_0}\rangle})\right]\right]\notag\\
	&=\mathbb{E}^\varphi_{\varphi\tilde\varphi}\left[
	\varphi({\hat{\xi}^\varphi}_{t_0})^{-1}\sU_{t-t_0}[g]({\hat{\xi}^\varphi}_{t_0}) 
	\right]\notag\\
	&=\langle \varphi^{-1}\sU_{t-t_0}[g], {\varphi\tilde\varphi}\rangle\notag\\
	&= \tilde{a}_{t-t_0}[g],
	\label{tildeaggg}
\end{align}
we can bound the following deviation using as again two error terms 
\begin{align*}
    &\left| \frac{\sU_t[g](x)}{\varphi(x)} - \tilde{a}_{t-t_0}[g] \right| 
    \\
    &\quad\leq 
    \left|\mathbb{E}^\varphi_{\delta_x}\left[\frac{1}{\langle \varphi, X_{t}^{(t_0,t]}\rangle}\left(1- {\rm e}^{-\langle g, X_{t}^{(t_0,t]}\rangle}\right)\right]
		-\mathbb{E}^\varphi_{\varphi\tilde\varphi}\left[\frac{1}{\langle \varphi, X_{t}^{(t_0,t]}\rangle}\left(1- {\rm e}^{-\langle g, X_{t}^{(t_0,t]}\rangle}\right)\right]\right|
        \\
		&\qquad+
        \left|\mathbb{E}^\varphi_{\delta_x}\left[\frac{1}{\langle \varphi, X_{t}\rangle}\left(1- {\rm e}^{-\langle g, X_{t}\rangle}\right)\right]
		-\mathbb{E}^\varphi_{\delta_x}\left[\frac{1}{\langle \varphi, X_{t}^{(t_0,t]}\rangle}\left(1- {\rm e}^{-\langle g, X_{t}^{(t_0,t]}\rangle}\right)\right]\right|
        \\
    &\quad =: \epsilon_x^1[g](t_0, t) + \epsilon_x^2[g](t_0, t),
\end{align*}
where, {roughly speaking,} $\epsilon_x^1[g](t_0, t)$ measures the convergence to the {measure $\varphi\tilde\varphi$} for the post-$t_0$ mass and $\epsilon_x^2[g](t_0, t)$ measures the effect of the mass immigrated prior to $t_0$ surviving to time $t$.
By conditioning on $\mathcal{F}_{t_0}^{\hat{\xi}^\varphi}$ and applying the Markov property, as in \eqref{tildeaggg}, we get
\begin{align*}
    \epsilon_x^1[g](t_0, t) &= \left|\mathbb{E}^\varphi_{\delta_x}\left[
	\varphi({\hat{\xi}^\varphi}_{t_0})^{-1}\sU_{t-t_0}[g]({\hat{\xi}^\varphi}_{t_0}) 
	\right]-\langle \sU_{t-t_0}[g],\tilde\varphi\rangle\right|\\
    &= \left|\varphi(x)^{-1}\sT_{t_0}[
	\sU_{t-t_0}[g]](x)-\langle \sU_{t-t_0}[g],\tilde\varphi\rangle\right|,
\end{align*}
where we recall that the expectation of any test function $f$ evaluated at the spine's position at time $s$ satisfies 
$$\mathbb{E}^\varphi_{\delta_x}\left[ f({\hat{\xi}^\varphi}_s) \right] = \frac{1}{\varphi(x)} \mathsf{T}_s[\varphi f](x).$$
Using the Perron-Frobenius type assumption \ref{H2}, this yields
\begin{equation}\label{eq:eps1super}
\sup_{x\in E, g\in B^+(E)} \frac{\epsilon_x^1[g](t_0, t)}{a_{t-t_0}} \le \frac{\Delta_{t_0}\tilde{a}_{t-t_0}[g]}{a_{t-t_0}}\leq \frac{\Delta_{t_0}a_{t-t_0}[g]}{a_{t-t_0}}\leq \Delta_{t-t_0},
\end{equation}
where we {recall that $\tilde{a}_{t-t_0}[g] = \langle \sU_{t-t_0}[g], \tilde\varphi\rangle$ and} used the inequalities
$$
 {\sU_{t-t_0}[g](x)=1-\mathrm{e}^{-\sV_{t-t_0}[g](x)} \leq \sV_{t-t_0}[g](x) \leq \sV_{t-t_0}(x).}
$$
Regarding the second error term, we reproduce the steps of \eqref{bigRHS} but conditioning on $\mathcal{F}_{t_0}^{\hat{\xi}^\varphi}$ instead of $\mathcal{F}_{t_0}^\xi$. In addition, taking into account that
\begin{align*}
		\mathbb{E}^\varphi_{\delta_x}\left[
		(1- {\rm e}^{-\langle g , X^{(0, t_0]}_t \rangle})\mid \mathcal{F}^{\hat{\xi}^\varphi}_{t_0}
		\right] \leq 
		\mathbb{E}^\varphi_{\delta_x}\left[
		\mathbf{1}_{\{\langle g , X^{(0, t_0]}_t \rangle>0\}}\mid \mathcal{F}^{\hat{\xi}^\varphi}_{t_0}\right],
	\end{align*}
and with the arguments from \eqref{makeitgotozero}, \eqref{C1}, \eqref{C2} and \eqref{useDCT}, we see that
\begin{align*}
    &\mathbb{E}^\varphi_{\delta_x}\left[
		\mathbf{1}_{\{\langle g , X^{(0, t_0]}_t \rangle>0\}}\mid \mathcal{F}^{\hat{\xi}^\varphi}_{t_0}\right] 
        \\
        &\qquad \leq \left( (2 \|c\| + C_1+K_1) \cdot t_0 \|\sV_{t-t_0}\| + (C_2 +K_2)  \cdot t_0 \|\sV_{t-t_0}\|^{\delta} \right) =: B(t_0,t).
\end{align*}

Therefore, similarly to \eqref{bigRHS},
\begin{align}
		\epsilon^2_x[g](t_0, t)
		&\leq 
		\mathbb{E}^\varphi_{\delta_x}\left[
		\mathbb{E}^\varphi_{\delta_{\hat{\xi}_{t_0}}}
		\left[
		\frac{1-(1-{\rm e}^{-\langle g , X_{t-t_0} \rangle})}{\langle \varphi , X_{t-t_0} \rangle } 
		\right]
		\mathbb{E}^\varphi\left[
		\mathbf{1}_{\{\langle g , X^{(0, t_0]}_t \rangle>0\}}\mid \mathcal{F}^{\hat{\xi}^\varphi}_{t_0}\right]
		\right]
		\notag\\
		&\qquad+
		\mathbb{E}^\varphi_{\delta_x}\left[
		\mathbb{E}^\varphi_{\delta_{\hat{\xi}_{t_0}}}\left[\frac{1}{\langle \varphi , X_{t-t_0} \rangle}(1- {\rm e}^{-\langle g , X_{t-t_0} \rangle})\right]
		\mathbb{E}^\varphi\left[
		\mathbf{1}_{\{\langle \varphi , X^{(0, t_0]}_t \rangle>0\}}\mid \mathcal{F}^{\hat{\xi}^\varphi}_{t_0}\right]
		\right]
        \notag
        \\
        &\leq
        \mathbb{E}^\varphi_{\delta_x}\left[
		\mathbb{E}^\varphi_{\delta_{\hat{\xi}_{t_0}}}\left[\frac{1}{\langle \varphi , X_{t-t_0} \rangle}\right]
		B(t_0,t)
		\right]
        =
        B(t_0,t) \mathbb{E}^\varphi_{\delta_x}\left[
		\frac{1-\mathrm{e}^{-\sV_{t-t_0}(\hat{\xi}_{t_0}^\varphi)}}{\varphi(\hat{\xi}_{t_0}^\varphi)}
		\right]\notag
        \\
        &=
        B(t_0,t) 
		\frac{\sT_{t_0}[1-\mathrm{e}^{-\sV_{t-t_0}}](x)}{\varphi(x)}. \notag
	\end{align}
Taking the supremum over $E$ and $B^+(E)$ and using \ref{H2}, we get
\begin{equation}\label{eq:eps2super}
\sup_{x\in E, g\in B^+(E)} \frac{\epsilon_x^2[g](t_0, t)}{a_{t-t_0}} \le B(t_0,t)(1 + \Delta_{t_0}),
\end{equation}
{where we also used that $\langle 1-\mathrm{e}^{-\sV_{t-t_0}}, \tilde \varphi \rangle \le \langle \sV_{t-t_0},  \tilde \varphi \rangle = a_{t-t_0}$.}
We can strategically choose a mapping $t\mapsto t_0(t)$ so that $\Delta_{t-t_0} \to 0$ and $B(t_0,t)(1 + \Delta_{t_0}) \to 0$ as $t\to \infty$.
This guarantees, {together with \eqref{eq:eps1super} and \eqref{eq:eps2super}} that
$$\lim_{t\to\infty} \sup_{g\in B^+(E)} \frac{1}{a_{t-t_0}} \left\| \frac{\sU_t[g]}{\varphi} - \tilde{a}_{t-t_0}[g] \right\| = 0.$$
Proceeding exactly as in the proof of Lemma \ref{lem:approxBMP}, see the calculations after \eqref{bootstrap1}, we deduce that 
$$\lim_{t\to\infty} \sup_{g\in B^+(E)} \frac{1}{a_{t}} \left\| \frac{\sU_t[g]}{\varphi} - \tilde{a}_{t}[g] \right\| = 0.$$
Finally, we translate the limits from the exponential formulation $\sU_t[g]$ to the linear formulation $\sV_t[g]$. For any $z \ge 0$, Taylor's theorem dictates that $0 \le z - (1 - e^{-z}) \le \frac{1}{2} z^2$.
Applying this, we have
$$0 \le \sV_t[g] - \sU_t[g] \le \frac{1}{2} \sV_t[g]^2 \le \frac{1}{2} \sV_t^2\quad \text{and}\quad 0 \le a_t[g] - \tilde{a}_t[g] \le \frac{1}{2} \langle \sV_t^2, \tilde\varphi\rangle \le \frac{1}{2} \|\sV_t\| a_t.$$
By applying the triangle inequality, we obtain
$$\left\| \frac{\sV_t[g]}{\varphi} - a_t[g] \right\| \le \left\| \frac{\sV_t[g] - \sU_t[g]}{\varphi} \right\| + \left\| \frac{\sU_t[g]}{\varphi} - \tilde{a}_t[g] \right\| + | \tilde{a}_t[g] - a_t[g] |.$$
Dividing by $a_t$ and taking the supremum over $g$, all three terms independently vanish as $t \to \infty$, successfully yielding our final result.
\end{proof}

\begin{lemma}\label{lem-ratio}
Assume that hypotheses \ref{H2}--\ref{H6} hold. For all $f\in B^+(E)$ with $\langle f, \tilde\varphi\rangle >0$, there exists $\kappa_0=\kappa_0(f)>0$ small enough such that
\begin{equation*}
    \lim_{t\to\infty} \left\| \frac{\normalfont \sV_{\kappa t}[a_t f]}{a_{\kappa t} [a_t f]\varphi} -1\right\| = 0,
\end{equation*}
for all $\kappa\in(0,\kappa_0).$
\end{lemma}
\begin{proof}
Fix $f\in B^+(E)$ with $\langle f, \tilde\varphi\rangle >0$. Using the same decomposition as in \eqref{decomp}, it is enough to prove that
$$\lim_{t\to\infty}\frac{a_{\kappa t}[a_t f]}{a_t}>0.$$
We look for a lower bound of $a_{\kappa t}[a_t f]$. In order to do that, from the many-to-one formula, see \cite[Lemma 3]{GHK}, we deduce that
$
\mathsf{T}_r[1] \leq \mathrm{e}^{r \|(\mathsf{m}[1]-1)\beta+b\|}, \ r\geq 0,
$
and we also use that
$
\mathsf{J}[f] \leq  cf^2 +\mathsf{M}[f], \ f\in B^+(E),
$
where $$\mathsf{M}[f](x):=f(x) \int_{(0,\infty)} y\nu(x, \dd y) + \beta(x) \int_{M_0(E)} \langle f, \nu\rangle \Gamma (x, \dd \nu).$$ Then, together with the help of Lemma \ref{lem:Jmonot}, {\eqref{atf} and \eqref{nonlinvJ}, we deduce}
\begin{align}
    a_{\kappa t}[a_t f] 
    &= 
    a_t \langle f , \tilde\varphi\rangle -\int_{0}^{\tau} \langle \mathsf{J} [\mathsf{V}_r[a_t f]], \tilde\varphi\rangle \dd r - \int_{\tau}^{\kappa t} \langle \mathsf{J} [\mathsf{V}_r[a_t f]], \tilde\varphi\rangle \dd r
    \notag
    \\
    &\geq
    a_t \langle f , \tilde\varphi\rangle - \int_{0}^{\tau} \langle c \mathsf{V}_r[a_t f]^2, \tilde\varphi\rangle \dd r - \int_{0}^{\tau} \langle \mathsf{M} [\mathsf{V}_r[a_t f]], \tilde\varphi\rangle \dd r - \int_{\tau}^{\kappa t} \langle \mathsf{J} [\mathsf{T}_r [a_t f]], \tilde\varphi\rangle \dd r 
    \notag\\
    &\geq
    a_t \langle f , \tilde\varphi\rangle - \|c\| \int_{0}^{\tau} \langle \mathsf{T}_r[a_t f]^2, \tilde\varphi\rangle \dd r - \int_{0}^{\tau} \langle \mathsf{M} [\mathsf{T}_r [a_t f]], \tilde\varphi\rangle \dd r \notag\\
    &\quad - \int_{\tau}^{\kappa t} \langle \mathsf{J} [(1 + \Delta_r) \langle f,\tilde\varphi\rangle a_t \varphi], \tilde\varphi\rangle \dd r 
    \notag\\
    &\geq
    a_t \langle f , \tilde\varphi\rangle - \left(a_t^2 \langle f , \tilde\varphi\rangle \|c\|\|f\| + a_t \langle \mathsf{M} [\|f\|], \tilde\varphi\rangle \right)\int_{0}^{\tau}  \mathrm{e}^{r\|(\mathsf{m}[1]-1)\beta+b\|} \dd r \label{akttttsup}\\
    &\quad - \kappa_0 t  \Big \langle \mathsf{J} \big[(1 + \sup_{r\geq \tau}\Delta_r) \langle f,\tilde\varphi\rangle a_t \varphi \big], \tilde\varphi \Big \rangle. \notag
\end{align}
Hence
$$
\lim_{t\to\infty}\frac{a_{\kappa t}[a_t f]}{a_t} \geq \langle f , \tilde\varphi\rangle -  \langle \mathsf{M} [\|f\|], \tilde\varphi\rangle \int_{0}^{\tau}  \mathrm{e}^{r\|(\mathsf{m}[1]-1)\beta+b\|} \dd r - \frac{\kappa_0(1 + \sup_{r\geq \tau}\Delta_r)^{1+\alpha} \langle f,\tilde\varphi\rangle^{1+\alpha}}{\alpha}.
$$
Taking first $\tau=\tau(f)>0$ and then $\kappa_0=\kappa_0(f)>0$ both small enough such that 
$$\langle \mathsf{M} [\|f\|], \tilde\varphi\rangle \int_{0}^{\tau}  \mathrm{e}^{r\|(\mathsf{m}[1]-1)\beta+b\|} \dd r < \frac{1}{3}\langle f, \tilde\varphi\rangle$$
and
$$\frac{\kappa_0(1 + \sup_{r\geq \tau}\Delta_r)^{1+\alpha} \langle f,\tilde\varphi\rangle^{1+\alpha}}{\alpha}<\frac{1}{3}\langle f, \tilde\varphi\rangle,$$
we get the desired result.
\end{proof}

\begin{lemma}\label{lemineq}

Assume \ref{H2}--\ref{H6}. For all $f\in B^+(E)$ and $\varepsilon\in (0,1 )$, there exists $\kappa=\kappa(f,\varepsilon)\in(0,\varepsilon)$ small enough such that
$$
(1-\varepsilon)^2\langle f, \tilde \varphi\rangle a_t\varphi \leq \normalfont \sV_{\kappa t}[a_t f] \leq (1+\varepsilon)\langle f, \tilde\varphi\rangle a_t\varphi
$$
for all $t$ sufficiently large.
\end{lemma}
\begin{proof}

The upper and lower bounds for the linear semigroup in \eqref{useupperandlower} still hold in the superprocess setting due to \ref{H2}. Therefore, by \eqref{nonlinvJ}, 
	\[
	\sV_{\kappa t} [a_t f] 
	\leq 
	\sT_{\kappa t} [a_t f] 
	\leq 
	\langle f ,\tilde{\varphi}\rangle (1 + \Delta_{\kappa t})a_t\varphi.
	\]
	and we can choose $t$ large enough such that 
	$\Delta_{\kappa t} \leq \varepsilon $.
	
	For the lower bound in the case $\langle  f, \tilde\varphi\rangle>0$ (otherwise, the {inequality trivially holds}), we proceed as in the particle process setting. We use Lemma \ref{lem-ratio} so that there is $t_0=t_0(f,\varepsilon)>0$ large enough verifying
    $$
\sV_{\kappa t}[a_t f] \geq (1-\varepsilon)a_{\kappa t}[a_t f]\varphi , \qquad t\geq t_0.
$$
Then, from \eqref{akttttsup},
\begin{align*}
    a_{\kappa t}[a_t f] 
    &\geq
    a_t \langle f , \tilde\varphi\rangle - \left(a_t^2 \langle f , \tilde\varphi\rangle \|c\|\|f\| + a_t \langle \mathsf{M} [\|f\|], \tilde\varphi\rangle \right) \int_{0}^{\tau} \mathrm{e}^{r\|(\mathsf{m}[1]-1)\beta+b\|} \dd r \\
    &\quad - {\kappa t} \langle \mathsf{J} [(1 + \sup_{r\geq \tau}\Delta_r) \langle f,\tilde\varphi\rangle a_t \varphi], \tilde\varphi\rangle.
\end{align*}
We take $\tau=\tau(f,\varepsilon)>0$ small enough so that
$$
\langle \mathsf{M} [\|f\|], \tilde\varphi\rangle  \int_{0}^{\tau} \mathrm{e}^{r\|(\mathsf{m}[1]-1)\beta+b\|} \dd r < \frac{\varepsilon}{3} \langle f , \tilde\varphi\rangle.
$$
and $\kappa=\kappa(f,\varepsilon)\in(0,\varepsilon)$ small enough such that
$$
\frac{\kappa (1 + \sup_{r\geq \tau}\Delta_r)^{1+\alpha} \langle f,\tilde\varphi\rangle^{\alpha}}{\alpha} < \frac{\varepsilon}{3},
$$
because
\begin{align*}
    {\kappa t} \langle \mathsf{J} [(1 + \sup_{r\geq \tau}\Delta_r) \langle f,\tilde\varphi\rangle a_t \varphi], \tilde\varphi\rangle 
    \sim
    a_t \langle f,\tilde\varphi\rangle  \frac{\kappa (1 + \sup_{r\geq \tau}\Delta_r)^{1+\alpha} \langle f,\tilde\varphi\rangle^{\alpha}}{\alpha}, \qquad \text{as } t\to\infty.
\end{align*}
In addition, for all $t$ sufficiently large, we have
$$a_t \|c\|\|f\|  \int_{0}^{\tau} \mathrm{e}^{r\|(\mathsf{m}[1]-1)\beta+b\|} \dd r< \frac{\varepsilon}{3}.$$
Then we showed that
$a_{\kappa t}[a_t f] \geq (1-\varepsilon) a_t \langle f , \tilde\varphi\rangle,$
concluding the proof.
\end{proof}

To conclude we prove the Yaglom limit for non-local superprocesses.

\begin{proof}[Proof of Theorem~\ref{theo:Yaglom}]

Taking into account the decomposition
$$\frac{\mathbb{E}_{\mu}\left[1-\mathrm{e}^{-a_t \langle f,  X_t\rangle} \right]}{\mathbb{P}_{\mu}(\zeta >t)}=\frac{1-\mathrm{e}^{-\langle \sV_t[a_t f],\mu \rangle }}{\langle \sV_t[a_t f],\mu \rangle}\frac{\langle \sV_t[a_t f],\mu \rangle}{a_t \langle \varphi,\mu\rangle}\frac{a_t \langle \varphi,\mu\rangle }{\mathbb{P}_{\mu}(\zeta >t)},$$
it is enough to show that
\begin{equation}\label{itsenoughtoprovesup}
    \lim_{t\to\infty} \frac{a_t[a_tf]}{a_t} = (1+\langle  f , \tilde{\varphi}\rangle ^{-\alpha})^{-1/\alpha}.
\end{equation}
Indeed, clearly $(1-\mathrm{e}^{-\langle \sV_t[a_t f],\mu \rangle})/\langle \sV_t[a_t f],\mu \rangle\to 1$ as $t\to\infty$ and, by Theorem \ref{theo:kolmogorov},
$$\lim_{t\to\infty} \frac{a_t \langle \varphi,\mu\rangle }{\mathbb{P}_{\mu}(\zeta >t)} = 1.$$
Further,
$$\frac{\langle \sV_t[a_t f],\mu \rangle}{a_t \langle \varphi,\mu\rangle} = \frac{1}{\langle \varphi,\mu\rangle}\Big\langle \varphi \frac{\varphi^{-1}\sV_t[a_t f]-a_t[a_tf]}{a_t},\mu \Big\rangle+\frac{a_t[a_tf]}{a_t},$$
so we just need to show \eqref{itsenoughtoprovesup} because, from Lemma \ref{lem:approx},
$$\lim_{t\to \infty} \left\|\varphi \frac{\frac{\sV_t[a_t f]}{\varphi}-a_t[a_tf]}{a_t}\right\| =0.$$
From this point, the proof is verbatim the particle setting. For $f\in B^+(E)$ and $\varepsilon>0$, we choose
$$n(t) = \min \left\{ m \in \mathbb{N} : a_m \leq \frac{(1-\varepsilon)^2}{1+\varepsilon} \langle f, \tilde\varphi\rangle a_t \right\}$$
and
$$k(t) = \max \left\{ m \in \mathbb{N} : a_m \geq \frac{1+\varepsilon}{1-\varepsilon} \langle f, \tilde\varphi\rangle a_t \right\},$$
so that
$$
a_{k(t)+1} \leq \frac{1+\varepsilon}{1-\varepsilon} \langle f, \tilde\varphi\rangle a_t \leq a_{k(t)} \qquad \text{and} \qquad a_{n(t)} \leq \frac{(1-\varepsilon)^2}{1+\varepsilon} \langle f, \tilde\varphi\rangle a_t \leq a_{n(t)-1}.
$$
	Combining above with Lemmas \ref{lem:asymV} and \ref{lemineq}, we deduce for $t$ large enough that
$$ \sV_{n(t)}  \leq (1+\varepsilon)a_{n(t)}\varphi \leq  (1-\varepsilon)^2\langle f, \tilde\varphi\rangle a_t\varphi \leq \sV_{\kappa t}[a_t f] \leq (1+\varepsilon)\langle f, \tilde\varphi\rangle a_t\varphi\leq (1-\varepsilon)a_{k(t)}\varphi \leq \sV_{k(t)} .$$
	Then, from the semigroup property, we see that
	\[
	\sV_{(1-\kappa)t+n(t)} = \sV_{(1-\kappa)t}[\sV_{n(t)}] 
	\leq 
	\sV_{(1-\kappa)t}[\sV_{\kappa t}[a_t f]] = \sV_{t}[a_t f]
	\leq 
	\sV_{(1-\kappa)t}[\sV_{k(t)}] = \sV_{(1-\kappa)t+k(t)}
	\]
    and
    $$\frac{a_{(1-\kappa)t+n(t)}}{a_t} \leq \frac{a_{t}[a_t f]}{a_t} \leq \frac{a_{(1-\kappa)t+k(t)}}{a_t}.$$
	The rest of the details follow in a straightforward way from the branching Markov process setting.
\end{proof}

\section*{Appendix: Proof of Theorem \ref{superspine} in brief}\label{appn}

We give a brief proof  of the spine decomposition for non-local superprocesses given in Theorem \ref{superspine}. We don't bother giving a full proof as the experienced reader will know that, once the main components of the decomposition have been identified, the structural proof of how they fit together remains the same; see for example the proof of Theorem 11.1 in \cite{Horton2023} in the branching Markov process setting  as well as \cite{ren2020limit, LRSS2023} in the superprocess setting. We start with some preparatory material.

Recall that we write $(\xi_t, t\geq0)$, with probabilities $\mathbf{P}=(\mathbf{P}_x, x\in E)$ as the Markov process corresponding to the semigroup $({\sP}_t, t\geq0)$.
Recall also that the stochastic process $\hat{\xi} = (\hat{\xi}_t, t\geq0)$ with probabilities $\hat{\mathbf{P}} = (\hat{\mathbf{P}}_x, x\in E)$  has the same Markovian increments as  $(\xi, \mathbf{P})$, albeit  interlaced in time with additional {\it `special jumps'} which, when $\hat\xi$ is positioned at $x\in E$, arrive at instantaneous rate $\beta(x) {\sm}[\varphi](x)/\varphi(x)$ and are distributed according to the kernel $\kappa(x,\dd y)$, $y\in E$, where 
\[
\int_E f(y)\kappa(x, \dd y) =\frac{{\sm}[f\varphi](x)}{{\sm}[\varphi](x)}, \qquad f\in B^+(E).
\]
Moreover, define for $f\in B^+(E)$,
\[
\hat{\sP}_t^\varphi[f](x) = \frac{\sT_t\bra{f\varphi}(x)}{\varphi(x)}, \qquad t\geq0, \quad x \in E.
\] 

Below we prove two key results which form the main components of the  proof of Theorem \ref{superspine} and conclude the \hyperref[appn]{Appendix} with further discussion about additional components to the proof which we leave to the reader for the sake of brevity as they are otherwise standard given the references we have given above.

\begin{lemma}\label{superspinepre1} The semigroup $\hat{\emph{\sP}}^\varphi = (\hat{\emph{\sP}}_t^\varphi, t\geq0)$ is that of a conservative Markov process that corresponds to a Doob $h$-transform of $(\hat{\xi}, \hat{\mathbf{P}})$ via 
\begin{equation}
\hat{\emph{\sP}}_t^\varphi[f](x)= \hat{\mathbf{E}}_x\left[f(\hat{\xi}_t)\frac{\varphi(\hat{\xi}_t)}{\varphi(x)}{\rm e}^{\int_0^t \frac{\beta(\hat{\xi}_s)}{\varphi(\hat{\xi}_s)}(\emph{\sm}[\varphi](\hat{\xi}_s)-\varphi(\hat{\xi}_s)) + b(\hat{\xi}_s)\dd s}\right], \quad x\in E, \ t\geq0, \ f\in B^+(E).
\label{properMarkov}
\end{equation}
Moreover, $\hat{\emph{\sP}}^\varphi$ is equivalent to a conservative Markov process with semigroup $(\emph{\sP}^\varphi_t, t\geq0)$ given by 
\begin{equation}
\emph{\sP}^\varphi_t[f](x) = \mathbf{E}_x\left[f(\xi_t)\frac{\varphi(\xi_t)}{\varphi(x)}{\rm e}^{\int_0^t \frac{\beta(\xi_s)}{\varphi(\xi_s)}(\emph{\sm}[\varphi](\xi_s)-\varphi(\xi_s)) +
 b(\xi_s)  \dd s}\right], \quad x\in E, \ t\geq0, \ f\in B^+(E),
\label{Q}
\end{equation}
 interlaced  with additional {\it `special jumps'} with  instantaneous rate $\beta(x) \emph{\sm}[\varphi](x)/\varphi(x)$ and distribution given by $\kappa(x,\dd y)$, $y\in E$.
\end{lemma}
\begin{proof}
Throughout this proof we will appeal repeatedly to certain recursive integral equations. In particular we will identify and equate different solutions to these integral equations. Uniqueness to the recursive equations we use is therefore an implicit part of our reasoning. 
In all cases, uniqueness will follow thanks to standard reasoning of Gr\"onwall's Lemma; the details are left to the reader and referred to Chapter 2 and 8 of \cite{Horton2023} for further background.

Recall \eqref{superm21} and note that, with rearrangement, 
\begin{align}
\hat{\sP}_t^\varphi[f](x)&=\frac{\sP_t[f\varphi](x)}{\varphi(x)}
+\int_{0}^{t}\frac{1}{\varphi(x)}\sP_s\bra{\beta (\sm[\sT_{t-s}[f\varphi]]-\sT_{t-s}[f\varphi])+b\sT_{t-s}[f\varphi]} (x)\dd s\notag\\
&=\frac{\sP_t[f\varphi](x)}{\varphi(x)}
+\int_{0}^{t}\frac{1}{\varphi(x)}\sP_s\bra{\varphi \beta \frac{\sm[\varphi] }{\varphi}\Bigg(\frac{\sm[\varphi\hat{\sP}_{t-s}^\varphi[f]]}{\sm[\varphi]}-\hat{\sP}_{t-s}^\varphi[f]\Bigg)} (x)\dd s\notag\\
&\hspace{0.5cm}+ \int_0^t \frac{1}{\varphi(x)}\sP_s\bra{
\beta\Big({\sm[\varphi]}-{\varphi} \Big)\hat{\sP}_{t-s}^\varphi[f]
+\varphi b\hat{\sP}_{t-s}^\varphi[f]}(x)\dd s.
\label{introduceQ}
\end{align}
Using Theorem 2.1 of \cite{Horton2023} and the fact that $t\mapsto \sP_t[\varphi f]/\varphi$ is also a semigroup,  \eqref{introduceQ} can otherwise be written as 
\begin{align}
\hat{\sP}_t^\varphi[f](x)&=\sP^\varphi_t[f](x) + \int_0^t \sP^\varphi_s\bra{\beta \frac{\sm[\varphi] }{\varphi}\Bigg(\frac{\sm[\varphi\hat{\sP}_{t-s}^\varphi[f]]}{\sm[\varphi]}-\hat{\sP}_{t-s}^\varphi[f]\Bigg)}(x)\dd s.
\label{introduceQ2}
\end{align}
Next recall the  semigroup \eqref{Q},
then, appealing to Theorem 2.1 of \cite{Horton2023}, we note that, for $x\in E$, $t\geq0$,
\begin{equation}
\varphi(x)\sP^\varphi_t[1](x) =  \sP_t[\varphi](x) + \int_0^t \sP_s\left[\Big( \beta (\sm[\varphi] -\varphi ) + \varphi b \Big)\sP^\varphi_{t-s}[1]\right]\dd s.
\label{reflect1}
\end{equation}

Note also by setting $f = \varphi$ in \eqref{superm21}  and recalling that $\sT_t[\varphi] = \varphi$ for all $t\geq0$, we see that, for $x\in E$  and $t\geq0$,
\begin{equation}
\varphi(x) =  \sP_t[\varphi](x)  + \int_0^t \sP_s\left[\Big( \beta (\sm[\varphi] -\varphi ) + \varphi b \Big)\right]\dd s.
\label{reflect2}
\end{equation}
Comparing \eqref{reflect1} and \eqref{reflect2} we conclude that $\sP_t^\varphi[1] = 1$ and therefore that $(\sP^\varphi_t, t\geq0)$ is the semigroup of a conservative Markov process. 

Recalling the definition of $(\hat{\xi},\hat{\mathbf{P}})$,  we can further identify \eqref{introduceQ2} as a Markovian decomposition of \eqref{properMarkov} by splitting on the first special jump. 
Thanks to the eigen property of $\varphi$, $\hat{\sP}^\varphi_t[1] = \sT_t[\varphi]/\varphi =1$,  we conclude   that $\hat{\sP}^\varphi$ is the semigroup of a conservative Markov process.
\end{proof}

Next we move to the key component which identifies the nature of the spine decomposition (whilst not offering a complete proof of the spine decomposition).

\begin{theorem}\label{superspinepre2}Assume \ref{H2} and \ref{H6}. Suppose we denote  $\hat{X}^\varphi: = (\hat{X}^\varphi_t, t\geq0)$ the measure valued process whose evolution is described by the aggregated mass from  all but the first bullet points given in the decomposition Theorem \ref{superspine}. Denote by $\hat{\mathbb{P}}^\varphi_{\delta_x}$, $x\in E$, the associated probabilities. 
We have for $g\in B^+(E)$,
\begin{equation}
\mathbb{E}_{\delta_x}\left[ \frac{\langle \varphi, X_t\rangle}{\varphi(x)} {\rm e}^{-\langle g, X_t\rangle}\right]= {\rm e}^{- \emph{\sV}_t[g](x)}\hat{\mathbb{E}}^\varphi_{\delta_x}\left[{\rm e}^{-\langle g, \hat{X}^\varphi_t\rangle}\right], \qquad x\in E, \ t\geq0,\ g\in B^+(E).
\label{coreidentityforsuperspine}
\end{equation}
\end{theorem}
\begin{proof}
As with the previous proof, we will use throughout uniqueness of the integral equations discussed, which are left to the reader to verify. We point again to the use of  Gr\"onwall's Lemma and the style of reasoning in Chapter 2 and 8 of \cite{Horton2023} for further background. 

We start by noting that 
\begin{align*}
\mathbb{E}_{\delta_x}\left[ \frac{\langle \varphi, X_t\rangle}{\varphi(x)} {\rm e}^{-\langle g, X_t\rangle}\right] &=-\frac{1}{\varphi(x)} \frac{\partial}{\partial \theta}\mathbb{E}_{\delta_x}\left[  {\rm e}^{-\langle g+ \theta \varphi, X_t\rangle}\right]\Big|_{\theta =0}\notag\\
&=-\frac{1}{\varphi(x)} \frac{\partial}{\partial \theta}{\rm e}^{- \sV_t[g+\theta\varphi](x)}\Big|_{\theta=0}\notag\\
& = 
{\rm e}^{- \sV_t[g](x)}\frac{1}{\varphi(x)}\frac{\partial}{\partial \theta}\sV_t[g+ \theta\varphi](x)\Big|_{\theta =0}.
\end{align*}

Suppose now we define
\begin{equation}
\tilde{\sV}_t[g](x) = 
\frac{1}{\varphi(x)}\frac{\partial}{\partial \theta}\sV_t[g+ \theta\varphi](x)\Big|_{\theta =0}, \qquad x\in E,\ t\geq0,\ g\in B^+(E).
\label{differential}
\end{equation}
Note that as $\theta\to0$, 
\[
\sV_t[g+ \theta\varphi](x) \sim  \sV_t[g](x)+ \theta \varphi(x) \tilde{\sV}_t[g](x).
\]

We want to apply \eqref{differential} to \eqref{non-local-evolution-equation}, but first it is preferential to manipulate the latter. Indeed, considering the definition of the semigroup \eqref{Q}, we use Theorem 2.1 of \cite{Horton2023} to develop \eqref{non-local-evolution-equation} into 
\begin{align}
\frac{1}{\varphi(x)}\sV_t[g](x) &=\sP^\varphi_t\left[\varphi^{-1} g\right]\notag\\
&\hspace{0.5cm} - \int_{0}^{t}\sP^\varphi_s\Big[\varphi^{-1}\Big(\psi(\cdot,{\sV}_{t-s}[g])+\phi(\cdot,{\sV}_{t-s}[g])\notag\\
&\hspace{4cm}+\frac{\beta }{\varphi }({\sm}[\varphi] -\varphi ){\sV}_{t-s}[g] + b{\sV}_{t-s}[g] \Big)\Big](x)\dd s\notag\\
&=\sP^\varphi_t\left[\varphi^{-1} g\right]\notag\\
&\hspace{0.5cm}- \int_{0}^{t}\sP^\varphi_s\Bigg[\varphi^{-1}\Bigg(\sJ[\sV_{t-s}[g]]- \beta\sm[\varphi]   
\Big(\frac{\sm[\sV_{t-s}[g]]}{\sm[\varphi]}-\varphi^{-1}{\sV}_{t-s}[g] \Big) \Bigg)\Bigg](x)\dd s
\label{manipulated}
\end{align}
where we have used \eqref{local-branching-mechanism}, \eqref{non-local-branching-mechanism}, \eqref{supermh} and \eqref{J-H5}, so that
\[
\sJ[h](x) = \psi(x, h(x)) + \phi(x, h) + \beta(x)(\sm[h](x)-h(x)) + b(x)h(x), \qquad x\in E, \ h\in B^+(E).
\]

With \eqref{manipulated} in hand we have by dominated convergence, for $t\geq0$, $x\in E$ and $g\in B^+(E)$,
\begin{align}
\tilde{\sV}_t[g](x)&=\frac{\partial}{\partial \theta}\sP^\varphi_t[\varphi^{-1}g+\theta ](x)\Big|_{\theta =0}\notag\\
& 
\quad-\frac{\partial}{\partial \theta}\int_{0}^{t}\sP^\varphi_s\Bigg[\varphi^{-1}\Bigg(\sJ[\sV_{t-s}[g+\theta\varphi]]\notag\\
&\hspace{4cm}
-
 \beta\sm[\varphi]   
\Bigg(\frac{\sm[\sV_{t-s}[g+\theta\varphi]]}{\sm[\varphi]}-\varphi^{-1}{\sV}_{t-s}[g+\theta\varphi] \Bigg) \Bigg)\Bigg](x)\dd s\Bigg|_{\theta =0}\notag\\
&=1- \int_{0}^{t}\sP^\varphi_s\big[\varphi^{-1}\Phi[{\sV}_{t-s}[g], \varphi\tilde{\sV}_{t-s}[g]]\big](x)\dd s  \notag\\
&\hspace{4cm}+ \int_0^t\sP^\varphi_s\Bigg[ 
\varphi^{-1}\beta\sm[\varphi]  
\Bigg(\frac{\sm[\varphi \tilde{\sV}_{t-s}[g]]}{\sm[\varphi]}-\tilde{\sV}_{t-s}[g] \Bigg)  \Bigg](x)
 \dd s,\label{Vtilde}
\end{align}
where we have used the eigen property that $\sT_t[\varphi] = \varphi$, and for $h_1, h_2\in B^+(E)$,
\begin{align}
\Phi[h_1, h_2](x)&=2c(x)h_1(x)h_2(x) +h_2(x)\int_{(0,\infty)}( 1-{\rm e}^{-h_1(x) y})y\nu(x,\dd y)\notag\\
&\hspace{1cm} +\beta(x)\int_{M_{0}(E)}(1 -{\rm e}^{-\langle h_1,{\nu}\rangle} )\langle h_2, \nu\rangle\Gamma(x,\dd \nu).
\label{J'}
\end{align}

Recall we denoted by   $\hat{X}^\varphi: = (\hat{X}^\varphi_t, t\geq0)$ the measure valued process whose evolution is described by the aggregated mass from  all but the first bullet points given in the decomposition Theorem \ref{superspine}. Moreover, we denoted by $\hat{\mathbb{P}}^\varphi_{\delta_x}$, $x\in E$ its  associated probabilities. Next define
\[
 \hat{\sV}^\varphi_t[g](x):= \hat{\mathbb{E}}^\varphi_{\delta_x}\left[{\rm e}^{-\langle g, \hat{X}^\varphi_t\rangle}\right], \qquad x\in E,\ g\in B^+(E).
\]
We claim that $(\hat{\sV}^\varphi_t[g], t\geq0)$ also solves \eqref{Vtilde} and that the latter has a unique solution. 

Suppose we denote by $\xi^\varphi = (\xi^\varphi_t, t\geq0)$ with probabilities $\mathbf{P}^\varphi = (\mathbf{P}^\varphi_x, x\in E)$ as a version of the Markov process corresponding to the semigroup $(\sP^\varphi_t, t\geq0)$. From the description of the dressed spine in Theorem \ref{superspine} we can split according to the first `special jump'. In particular, conditioning on the spatial evolution of the spine and applying Campbell's formula (the reader is referred to the many calculations in e.g. \cite{RenSongYang2022, LRSS2023, Maren}), this gives us
\begin{align}
 \hat{\sV}^\varphi_t[g](x) &= \mathbf{E}^\varphi_x\left[{\rm e}^{-\int_0^t \beta(\xi^\varphi_s)\frac{\sm[\varphi](\xi^\varphi_s)}{\varphi(\xi^\varphi_s)}\dd s}
{\rm e}^{- \int_0^t \Psi_1[\sV_{t-s}[g]](\xi^\varphi_s)\dd s}
\right]\notag\\
& \quad+ \mathbf{E}^\varphi_x\Bigg[
\int_0^t \beta(\xi^\varphi_s)\frac{\sm[\varphi](\xi^\varphi_s)}{\varphi(\xi^\varphi_s)}
{\rm e}^{-\int_0^s\beta(\xi^\varphi_u) \frac{\sm[\varphi](\xi^\varphi_u)}{\varphi(\xi^\varphi_u)}\dd u}
{\rm e}^{- \int_0^s \Psi_1[\sV_{s-u}[g]](\xi^\varphi_u)\dd u}
\notag\\
&\hspace{1.75cm}
\frac{1}{\sm[\varphi](\xi^\varphi_s)}
\Bigg(\gamma(\xi^\varphi_s, \varphi \hat{\sV}^\varphi_t[g])+
 \int_{M_{0}(E)}{\rm e}^{-\langle \sV_{t-s}[g],{\nu}\rangle} \langle\varphi \hat{\sV}^\varphi_t[g], \nu\rangle\Gamma(\xi^\varphi_s,\dd \nu)
 \Bigg)
 \dd s\Bigg],
 \label{proposedspine}
\end{align}
where, for $h\in B^+(E)$, 
\[
\Psi_1[h](x)  = 2c(x)h(x) +\int_{(0,\infty)}( 1-{\rm e}^{-h(x) y})y\nu(x,\dd y).
\]
We have also used the fact that $\sV_t[f](x) = \mathbb{N}_x(1- {\rm e}^{-\langle f, X_t\rangle})$, which is a well known result due to \cite{DK04}.

Now let 
\[
\Psi_2[h_1,h_2](x) = \beta(x)\int_{M_{0}(E)}(1 -{\rm e}^{-\langle h_1,{\nu}\rangle} )\langle h_2, \nu\rangle\Gamma(x,\dd \nu) \qquad x\in E,\ h_1,h_2 \in B^+(E).
\]
so that 
\[
\Phi[h_1, h_2](x) = h_2(x)\Psi_1[h_1](x) + \Psi_2[h_1,h_2](x).
\]
Appealing to Lemma 2.1 of \cite{Horton2023}, and recalling \eqref{supermh}, we can transform \eqref{proposedspine} to take the form
\begin{align*}
 &\hat{\sV}^\varphi_t[g](x) \notag\\
 &=1 - 
  \mathbf{E}^\varphi_x\Bigg[
\int_0^t   
 \hat{\sV}^\varphi_{t-s}[g](x)\Psi_1[\sV_{t-s}[g]](\xi^\varphi_s)
 \dd s\Bigg]\notag\\
 &\quad -
  \mathbf{E}^\varphi_x\Bigg[
\int_0^t   
 \frac{\beta(\xi^\varphi_s)}{\varphi(\xi^\varphi_s)}
\Bigg(\int_{M_{0}(E)}\Big(1-{\rm e}^{-\langle \sV_{t-s}[g],{\nu}\rangle}\Big) \langle\varphi \hat{\sV}^\varphi_{t-s}[g], \nu\rangle\Gamma(\xi^\varphi_s,\dd \nu)\Bigg)
 \dd s\Bigg]\notag\\
 &\quad + \mathbf{E}^\varphi_x\Bigg[
\int_0^t   
 \frac{\beta(\xi^\varphi_s)}{\varphi(\xi^\varphi_s)}
 \Bigg(\gamma(\xi^\varphi_s, \varphi \hat{\sV}^\varphi_{t-s}[g])+\int_{M_{0}(E)}\langle\varphi \hat{\sV}^\varphi_t[g], \nu\rangle\Gamma(\xi^\varphi_s,\dd \nu)-\sm[\varphi](\xi^\varphi_s)\hat{\sV}^\varphi_{t-s}[g]
 \Bigg)
 \dd s\Bigg]\notag\\
 &=1- 
\int_0^t   \sP^\varphi_s\Bigg[
 \hat{\sV}^\varphi_{t-s}[g](x)\Psi_1[\sV_{t-s}[g]](\xi^\varphi_s) 
\Bigg] (x) \dd s - \int_0^t \sP^\varphi_s\Bigg[\varphi^{-1}\Psi_2[  {\sV}_{t-s}[g], \varphi\hat{\sV}_{t-s}^{\varphi}[g]]\Bigg](x)\dd s \notag\\
 &\quad 
 +\int_0^t \sP^\varphi_s\Bigg[
\varphi^{-1}\beta \sm[\varphi] \Bigg(\frac{\sm[ \hat{\sV}^\varphi_{t-s}[g]]}{\sm[\varphi]} -\hat{\sV}^\varphi_{t-s}[g]\Bigg)
 \Bigg](x)\dd s\notag\\
 &=1- \int_0^t \sP^\varphi_s\Bigg[\varphi^{-1}\Phi\left[  {\sV}_{t-s}[g], \varphi\hat{\sV}_{t-s}^{\varphi}[g]\right] +
 \varphi^{-1}\beta \sm[\varphi] \Bigg(\frac{\sm[ \hat{\sV}^\varphi_{t-s}[g]]}{\sm[\varphi]} -\hat{\sV}^\varphi_{t-s}[g]\Bigg)
 \Bigg](x)\dd s,
\end{align*}
and we see that $(\hat{\sV}^\varphi_t[g], t\geq0)$ solves \eqref{Vtilde}.

In conclusion we have that the distributional identity \eqref{coreidentityforsuperspine} holds as required.
\end{proof}

We conclude by noting that Lemma \ref{superspinepre1} and Theorem \ref{superspinepre2} are not sufficient to conclude the proof of Theorem \ref{superspine}. In order to do this, we note that the pair  $(\hat{\xi}^\varphi, \hat{X}^\varphi)$ are together Markovian by construction. We need to show that  the associated aggregate immigrated mass as per the statement of Theorem \ref{superspine} is alone Markovian. In that case, the conclusion of Theorem \ref{superspinepre2} tells us that it has the right transitions to match up with the Doob $h$-transform \eqref{COM}.  

In order to show that the aggregate immigrated mass from the statement of Theorem \ref{superspine} is Markovian, define 
\[
\hat{\chi}^\varphi = X' + \hat{X}^\varphi,
\]
where $X'$ is an indepndent copy of $(X, \mathbb{P}_{\delta_x})$. We need to show that 
\begin{equation}
\hat{\mathbb{E}}^\varphi_{\delta_x}\left[f(\hat{\xi}^\varphi_t){\rm e}^{-\langle g, \hat{\chi}^\varphi\rangle}\right] = \hat{\mathbb{E}}^\varphi_{\delta_x}\left[
\frac{\langle f\varphi, \hat{\chi}^\varphi\rangle}{\langle \varphi, \hat{\chi}^\varphi\rangle}
{\rm e}^{-\langle g, \hat{\chi}^\varphi\rangle}
\right] , \qquad x\in E, \ t\geq0,\ f,g\in B^+(E).
\end{equation}
This is Step 3 of the proof in the branching Markov process setting given in Theorem 11.1 of \cite{Horton2023}. We also leave this part of the proof to the reader as it is relatively straightforward and the reader can rely on the aforementioned reference as well as e.g. \cite{ren2020limit, LRSS2023}.

\section*{Acknowledgments} AEK and PM-C acknowledge the support of EPSRC programme grant EP/W026899/2. This work began in a research stay by NC-T at the University of Warwick. She gratefully acknowledges her host for the financial support, as well as for his kindness and  hospitality.

\bibliography{references}{}
\bibliographystyle{abbrv}

\end{document}